\documentclass{article}
\usepackage{fullpage}

\usepackage{amsmath,amsfonts}
\usepackage{mathrsfs}
\usepackage{color,xcolor}
\usepackage{ifpdf}
\usepackage{psfrag}
\usepackage{graphicx,graphics,subfigure,epsfig}
\usepackage{url,algorithm,algorithmic}
\usepackage{hyperref}
\usepackage{xspace}
\usepackage{nicefrac}
\usepackage{ulem}
\usepackage{epstopdf}
\DeclareMathOperator{\diag}{diag}

\newcommand{\dpar}[2]{\dfrac{\partial #1}{\partial #2}}

 \newcommand{\R}{\mathbb R}

\renewcommand{\P}{\mathbb P}

\newcommand{\bla}[1]{\textcolor{black}{#1}}

\newcommand{\bbf}{{\mathbf {f}}}
\newcommand{\bn}{{\mathbf {n}}}

\newcommand{\bu}{\mathbf{u}}
\newcommand{\bv}{\mathbf{v}}

\newcommand{\bbu}{\mathbf{u}}

\newcommand{\hbbf}{\hat{\mathbf{f}}}

\newcommand{\ba}{\mathbf{a}}

\newcommand{\bJ}{\mathbf{J}}
\newcommand{\bx}{\mathbf{x}}
\newcommand{\by}{\mathbf{y}}

\newcommand\eref[1]{(\ref{#1})}

\newcommand*\xbar[1]{%
  \hbox{%
    \vbox{%
      \hrule height 0.5pt 
      \kern0.4ex
      \hbox{%
        \kern-0.05em
        \ensuremath{#1}%
        \kern-0.00em
      }%
    }%
  }%
}

\begin{document}

\title{Active flux for triangular meshes for compressible flows problems}
\author{R\'emi Abgrall, Jianfang Lin and Yongle Liu\\
Institute of Mathematics, University of Z\"urich, Switzerland\\
\{remi.abgrall,jianfang.lin,yongle.liu\}@math.uzh.ch}
\date{}
\maketitle
\begin{abstract}
In this article, we show how to construct a numerical method for solving hyperbolic problems, whether linear or non-linear, using a continuous representation of the variables and their mean value in each triangular element. This type of approach has already been introduced by Roe, and others, in the multi-dimensional framework under the name of Active Flux, see \cite{AF1,AF2,AF3,AF4,AF5}. Here, the presentation is more general and follows \cite{Abgrall_AF,BarzukowAbgrall}. Various examples show the good behavior of the method in both linear and non-linear cases, including non-convex problems. The expected order of precision is obtained in both the linear and non-linear cases. This work represents a step towards the development of methods in the spirit of virtual finite elements for linear or non-linear hyperbolic problems, including the case where the solution is not regular.
\end{abstract}

\section{Introduction}
A new numerical method for solving hyperbolic problems was first introduced in \cite{AF1,AF2,AF3}. This new method uses the degrees of freedom (DoFs) of the quadratic polynomials, which are all on the boundary of the elements, and another degree of freedom that used is the cell average of the solution. This leads to a potentially third-order method, with a globally continuous approximation. The time integration relies on a back-tracing of characteristics. This has been extended to square elements, in a finite difference fashion, in \cite{AF4,AF5}. In \cite{Abgrall_AF}, a different, though very close, point of view is followed. The time stepping relies on a standard Runge-Kutta time stepping, and several versions of the same problem can be used without preventing local conservation. More precisely, a conservative version of the studied problem is used to update the average value and the local conservation property is guaranteed. In addition, we have more flexibility for the DoFs on the boundary of the elements. Namely, we can use several versions of the same problem, such as the primitive formulation, to update these DoFs. This is benefit from the fact that the point values are free of the local conservation and we do not need to restricted by the conservation formulation. These routes have been followed in \cite{BarzukowAbgrall} which gives several versions of higher order approximation, in one dimension.
 
The purpose of this paper is to extend the work of \cite{Abgrall_AF} to several dimensions, using simplex elements. We first describe the approximation space, then propose a discretization method, with some justifications, and then show numerical examples. We also show how to handle discontinuous solutions.
 
The hyperbolic problems we consider are written in a conservative form as:
\begin{equation}\label{eq:conservative}
\dpar{\bu}{t}+\text{ div }\bbf(\bbu)=0, \quad \bx \in \Omega\subset \R^d,
\end{equation}
with $d=2$, complemented by initial conditions
$$
\bbu(\bx,0)=\bbu_0(\bx), \qquad \bx\in \Omega,
$$
and boundary conditions that will be described case-by-case. The solution $\bu\in \R^m$ must stay in an open  subspace $\mathcal{D}\subset \R^d$ so that the flux $\bbf$ is defined. This is the invariant domain.

Two sets of problems are considered:
\begin{itemize}
\item the case of scalar conservation law with $\bbf$ linear taking the form of $\bbf(\bu,\bx)=\ba(\bx)\bu$ or non-linear. In that case, $\bu$ is a function with values in $\R$ and $\mathcal{D}=\R$ but the solution must stay in $[\min\limits_{\bx\in \R}\bu_0(\bx), \max\limits_{\bx\in \R}\bu_0(\bx)]$, due to Kruzhkov's theory: in these particular examples, we set $\mathcal{D}=[\min\limits_{\bx\in \R}\bu_0(\bx), \max\limits_{\bx\in \R}\bu_0(\bx)]$.
\item the case of the Euler equations, where $\bu=(\rho, \rho \bv, E)$, $\rho$ is the density, $\bv$ is the velocity, and $E$ is the total energy, i.e., the sum of the internal energy $e$ and the kinetic energy $\tfrac{1}{2}\rho \bv^2$. The flux is
$$
\bbf(\bu)=\begin{pmatrix}
\rho \bv\\
\rho \bv\otimes \bv+p\text{Id}_d\\
(E+p)\bv
\end{pmatrix}
$$
with $\text{Id}_d$ being a $d\times d$ identity matrix.
We have introduced the pressure that is related to the internal energy and the density by an equation of state. Here we make the choice of a perfect gas,
$$p=(\gamma-1)e$$
with $\gamma$ denoting the specific heat ratio. In all of the numerical examples tested in \S\ref{Euler_example}, we have taken $\gamma=1.4$.
In this example, $\mathcal{D}=\{\bu\in \R^4, \rho>0, e>0\}$.
\end{itemize}
 
We note that for smooth solutions, \eqref{eq:conservative} can be re-written in a non-conservative form,
$$\dpar{\bu}{t}+\bJ\cdot\nabla \bu=0,$$
where
$$\bJ\cdot\nabla \bu=A\dpar{\bu}{x}+B\dpar{\bu}{y},$$
with $A$ and $B$ being the Jacobians of the flux in the $x$- and $y$- direction, respectively. The system is hyperbolic, i.e., for any vector $\bn=(n_x, n_y)$, the matrix
$$\bJ\cdot\bn=A n_x+Bn_y$$ is diagonalisable in $\R$.

\section{Approximation space}
For now, we consider a triangle $E$ and the corresponding quadratic polynomial approximation. The degrees of freedom are the vertices $\sigma_i, i=1,\ldots 3$, the midpoints $\sigma_i, i=4,\ldots 6$, and the average, see Figure \ref{fig:normale2}. The first question whether we can find a polynomial space $V$ that contains at least the quadratic polynomials, and something more, in order to accommodate the average value as an independent variable, that is if we can find $7$ polynomials ($\varphi_i$) such that
\begin{equation}\label{phi_16}
  \varphi_i(\sigma_j)=\delta_i^j~\text{ for } ~1\leq i,j\leq 6, \quad\int_E\varphi_i(\bx)\; \bla{{\rm d}}\bx=0 ~\text{ for }~i=1, \ldots, 6,
\end{equation}
and
\begin{equation}\label{phi_7}
\varphi_7(\sigma_j)=0, \quad \int_E\varphi_7(\bx)\; \rm d\bx=|E|.
\end{equation}
An obvious choice is $\varphi_7=60 \lambda_1\lambda_2\lambda_3$ with $(\lambda_1,\lambda_2,\lambda_3)$ being the barycentric coordinates and $\varphi_i=$ the Lagrange basis function $+\beta_i \varphi_7$ with $\beta_i$ well chosen.
 
More precisely, the quadratic Lagrange polynomials are
$\varphi_i^L(\bx)=\lambda_i(2\lambda_i-1)$, $i=1,2,3$ and $\varphi^L_4=4\lambda_1\lambda_2$, $\varphi_5^L=4 \lambda_2\lambda_3$ and $\varphi_6^L=4\lambda_1\lambda_3$. Since 
$$\int_E\varphi_i^L\; {\rm d}\bx=0,$$ for $i=1,2,3$, we set $\varphi_i=\varphi_i^L$. Similarly, we see that we can set
$$\varphi_i=\varphi_i^L-\frac{1}{3}\varphi_7,$$ for $i=4,5,6.$ It is trivial to verify that the above-mentioned polynomials $\varphi_i$ satisfy the conditions \eref{phi_16} and \eref{phi_7}. Therefore, these polynomials serve as the quadratic basis functions of the polynomial space $V$.
 
The same idea can be extended to tetrahedrons and higher than quadratic degree. For example, we can get a cubic approximation by taking
$9$ Lagrange points with barycentric coordinates $(\tfrac{i}{3},\tfrac{j}{3},\tfrac{k}{3})$ with $i+j+k=\bla{3}$ and at least one  among the $\{i,j,k\}$ being set to $0$ so that we are on the boundary of the triangle $E$, and the average value. Hence the polynomials are:
$$\varphi_{10}=60\lambda_1\lambda_2\lambda_3,$$
and $$\varphi_{i}=\frac{1}{2}\lambda_i(3\lambda_i-1)(3\lambda_i-2)-2\lambda_1\lambda_2\lambda_3,\quad\text{ for }\quad i=1,2,3, $$ then
$$\varphi_4=\frac{9}{2}\lambda_1\lambda_{2}(3\lambda_1-1)-\frac{9}{2}\lambda_1\lambda_2\lambda_3,$$
$$\varphi_5=\frac{9}{2}\lambda_1\lambda_{2}(3\lambda_2-1)-\frac{9}{2}\lambda_1\lambda_2\lambda_3,$$
and the others $\{\varphi_6,\varphi_7,\varphi_8,\varphi_9\}$ are obtained by symmetry.

Our notations will be $\bu_\sigma$ for the point values and $\xbar \bu_E$ for the cell average. Finally, the solution $\bu$ can be approximated in each element by 
\begin{equation}\label{u1}
  \bu\approx\sum_{i=1}^{6}\bu_{\sigma_i}\varphi_i+\xbar \bu_E\varphi_7,
\end{equation}
for the third-order approximation and 
\begin{equation}\label{u2}
  \bu\approx\sum_{i=1}^{9}\bu_{\sigma_i}\varphi_i+\xbar \bu_E\varphi_{10},
\end{equation}
for the fourth-order approximation.

\section{Numerical schemes}\label{sec3}
In this section, we describe how to evolve in time the degrees of freedom. We use the method of line, with a SSP Runge-Kutta (third-order) in time. In the following, we describe the spatial discretization.
We first describe a high-order version, and then a first-order one.
\subsection{High order schemes}
The update of the average values is done by using the conservative formulation, i.e.,
\begin{equation}
\label{scheme:fv}\dpar{\bu}{t}+\text{div}\bbf(\bu)=0,
\end{equation}
both the conservative and a non-conservative version of the same system
\begin{equation}
\label{scheme:fd}\dpar{\bu}{t}+\bJ\cdot\nabla \bu=0.
\end{equation}
 
On a triangle $E$, $\bu$ is represented by point values at the Lagrange DoFs and its average on the element. The function  $\bu$ is a quadratic or cubic polynomial on each element, so that 
the update is done via:
\begin{equation}
\label{fv}\vert E\vert \dfrac{{\rm d}\bar\bu_E}{{\rm d}t}+\oint_{\partial E}\bbf(\bu)\cdot\bn\; \rm d\ell=0,
\end{equation} 
where the integral is obtained by using quadrature formula on each edge. In our simulations, we use a Gaussian quadrature formula with $3$ Gaussian points on each edge, so that there is no loss of accuracy.

The update of the boundary values is more involved and we describe several solutions.
Inspired by Residual Distribution schemes, and in particular the LDA (Low Diffusion advection A) scheme \cite{DeconinckMario}\footnote{The LDA scheme (LDA standing for Low Diffusion advection A scheme, there were a B and C version) is a scheme specialized to steady advection problems on triangular meshes for scalar problems. Point value residual are defined to update the solution in  a way similar as we do here. It is locally conservative, high-order and upwind but does not prevent spurious oscillations at discontinuities. It is the second order in its original version.  The accuracy is obtained by defining the residual at this point as the product of the upwind parameter with $J\cdot \nabla u$. The upwind property is obtained by looking at the sign of the scalar product of  local velocity (hence $J$ for a scalar problem) against the normal opposite to the vertex for which we evaluate the residual. Finally local conservation is guaranteed by multiplying the upwind parameter by a scaling factor, the so-called $N$ parameter, so that the sum or the fluctuations over the element is equal to $\int_E J\cdot \nabla u\, d\bx$. The $J$ Jacobians are evaluated so that this integral is equal to the integral over the boundary of the normal flux, by some generalization of the Roe average. An extension for systems was defined later in van der Weide's PhD thesis\cite{vanderWeide}, see also \cite{vanLeerCiCP}.}, one can define, for each of the polynomial DoFs $\sigma\in\{\sigma_1,\sigma_2,\ldots\}$ over the element $E$,
$$\Phi_\sigma^E(\bu)=\big (K^E_\sigma\big )^+N_\sigma\big ( \bJ(\bu_\sigma)\cdot\nabla \bu\big ), \quad N_\sigma^{-1}=\sum_{\sigma \in E}\big (K^E_\sigma\big )^+,$$ 
where, since the problem is hyperbolic and  whatever the vector $\bn_\sigma^E$, one can take the positive and negative parts of $K_\sigma^E$
$$
K^E_\sigma = (K^E_\sigma)^-+(K^E_\sigma)^+=R\Lambda R^{-1}=R\Lambda^- R^{-1}+R\Lambda^+R^{-1},
$$
with
\begin{align*}
    \Lambda^- &= \diag(\min(\lambda_1, 0), \dots, \min(\lambda_m, 0)), \\
    \Lambda^+ &= \diag(\max(\lambda_1, 0), \dots, \max(\lambda_m, 0)), \\
    K^E_\sigma &= \mathbf{J}(\bu_\sigma)\cdot \bn_\sigma^E.
\end{align*}
This normal vector, see Figure \ref{fig:normale2}, is defined as follows:
\begin{itemize}
\item for the quadratic approximation where we have used the DoFs of quadratic polynomials on the edges of the element, if $\sigma$ is a vertex, $\bn_\sigma$ is the inward normal to the opposite edge of $\sigma$, if $\sigma$ is a mid point, this is the outward normal to the edge where $\sigma$ is sitting.
\item for the cubic approximation where we have used the DoFs of cubic polynomials on the edges of the element, we apply the same procedure, described in the above case of quadratic approximation, to the vertices and the other DoFs on edges.
\end{itemize}
These two procedures can be reduced to the following. We notice that a degree of freedom belongs to two edges, the one sitting on its left and the one sitting on its right. These two edges can be identical. We then define the normal by adding the scaled normals of the outward normals to these edges.

Once this is done, we update $\bu$ by the following
$$\dfrac{\rm d\bu_\sigma}{\rm dt}+\sum\limits_{E, \sigma\in E} \Phi_\sigma^E=0.$$
The idea behind is to have some up-winding mechanism, 
however, the consistency is not clear: if we have a linear $\bu$, i.e., with $\nabla \bu$ constant, we would like to recover
$$\sum\limits_{E, \sigma\in E} \Phi_\sigma^E=\bJ(\bu_\sigma)\cdot \nabla \bu.$$
A priori, there is no reason for this. Hence a better idea is:
\begin{figure}[h]
\centerline{\includegraphics[trim=0.01cm 0.01cm 0.01cm 0.001cm,clip,width=4.5cm]{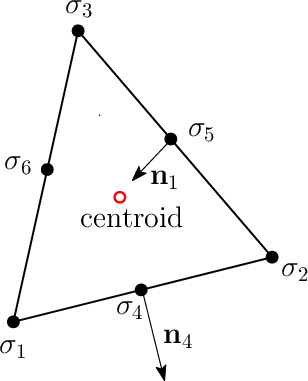}\hspace*{1.5cm}
\includegraphics[trim=0.01cm 0.01cm 0.01cm 0.001cm,clip,width=4.5cm]{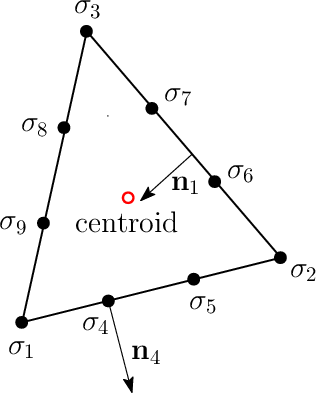}}
\caption{\label{fig:normale2} Normals for \eqref{scheme:pt}. Left: quadratic; Right: cubic.}
\end{figure}

\begin{subequations}\label{scheme:pt}
\begin{equation}\label{scheme:pt:1}\dfrac{{\rm d}\bu_\sigma}{{\rm d}t}+\sum\limits_{E, \sigma\in E} \Phi_\sigma^E(\bu)=0,\end{equation} with
\begin{equation}
\label{scheme:pt:2}\Phi_\sigma^E(\bu)=N_\sigma \big (K_\sigma^E)^+ \bJ(\bu_\sigma)\cdot \nabla \bu(\sigma),\end{equation}
\begin{equation}\label{scheme:pt:3}K_\sigma^E=\bJ(\bu_\sigma)\cdot \bn_\sigma^E \quad \text{and}\quad
N_\sigma^{-1}= \sum\limits_{E, \sigma\in E} \big (K_\sigma^E\big )^+.\end{equation}
\end{subequations}
Consistency is obvious.  From \cite{Abgrall99}, one can see that the matrix $N_\sigma$ exists for a symmetrizable system. A slightly different version is
$$\Phi_\sigma^E(\bu)=N_\sigma \text{ sign}\big (K_\sigma^E)\bJ(\bu_\sigma)\cdot \nabla \bu(\sigma), \qquad N_\sigma^{-1}=\sum\limits_{E, \sigma\in E} \text{sign}\big (K_\sigma^E\big ).$$
The numerical results are identical, and we will use the version \eqref{scheme:pt:3}.

\subsection{Low order schemes}\label{loworder}
For the average value, one can simply take a standard finite volume (so we do not use the boundary values). We have taken the Lax-Friedrichs numerical flux and the Roe scheme in our simulations.

The real problem is what to do with the point values.
A first-order version
\begin{figure}[h]
\centerline{\includegraphics[trim=0.02cm 0.01cm 0.01cm 0.01cm,clip,width=5.2cm]{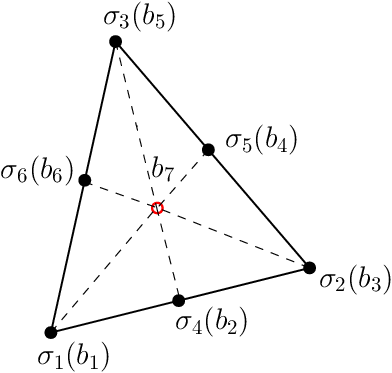}\hspace*{1.0cm}
\includegraphics[trim=0.02cm 0.01cm 0.01cm 0.01cm,clip,width=5.2cm]{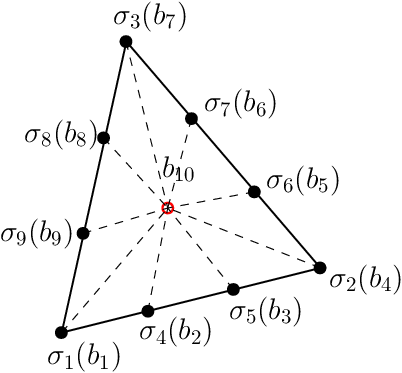}}
\caption{\label{fig:nscheme}  Geometry for the first-order scheme. Left: quadratic; Right: cubic.}
\end{figure}
can be defined as follows. We first, and temporarily for the definition of sub-elements only, re-number the boundary nodes. We pick one degree of freedom, call it $b_1$, and list the others in a counter-clockwise manner: we get the list $\{b_1, \ldots, b_{N-1}\}$ of all boundary degrees of freedom. The centroid is denoted by $b_N$, with $N=7$ for the quadratic approximation, and $N=10$ for the cubic one. Then we define the sub-elements of vertices $\{b_i, b_N, b_{i+1}\}_{i =1,\dots, N-2}$ and $\{b_{N-1}, b_N, b_1\}$ with a numbering $b_i$ modulo $N$. This defines $N-1$ sub-elements that we denote by $T^E_i$, then we get inward normals $\bn_{b_i}$. For quadratic elements, and using back the original numbering, this is illustrated in Figure \ref{fig:nscheme}. In what follows, we use the original numbering with this geometrical definition of sub-elements.

The update of $\bu_\sigma$ is done with \eqref{scheme:pt:1} where $\Phi_\sigma^E$ needs to be modified. We use
\begin{equation}
\label{schema:pt1:1}
\Phi_\sigma^E=\frac{1}{\vert C_\sigma\vert }\sum_{T^E_i, \sigma\in T^E_i}\Phi_{\sigma}^{T^E_i}(\bu).
\end{equation}
The residuals will be defined the following, which is done to get monotone first-order residual distribution schemes, see \cite{DeconinckMario,AbgrallEnergie} for example. We will use two types of residuals:
\begin{itemize}
\item a version that is inspired by the Local Lax-Friedrichs scheme
$$\Phi_\sigma^{T^E_i}=\frac{1}{3}\oint_{T^E_i} \bJ(\bu)\cdot \nabla \bu \; {\rm d}\bx+\alpha_{T^E_i}\big (\bu_\sigma-\overline{\bu}_{T^E_i}\big )$$
with
\begin{equation}\label{eq6}
\oint_{T^E_i} \bJ(\bu)\cdot \nabla \bu \; {\rm d}\bx\approx \vert T^E_i\vert \bJ(\overline{\bu}_{T^E_i})\nabla \bu(\bla{\xbar{\bu}_{T^E_i}}),
\end{equation}
where we use the $\P^1$ approximation on $T^E_i$  and
$\overline{\bu}_{T^E_i}$ is the arithmetic average of the $\bu$'s at  three vertices of $T^E_i$ (hence we use the average value here). Last, $$\alpha_{T^E_i}\geq \max\limits_{\substack{\bla{\sigma\in T^E_i,}\\\bn_{b_i} \text{normals of~} T^E_i}}\rho\big ( \bJ(\bla{\bu_\sigma})\cdot \bn_{b_i}\big ),$$ where $\rho(B)$ is the spectral radius of the matrix $B$.

\item a version that is inspired by the finite volume scheme re-written in the residual distribution framework, see \cite{ConservationRD} for details. We have used the Roe scheme with Harten-Yee entropy fix.
\end{itemize}
We still need to define $\vert C_\sigma\vert$. It is 
$$\vert C_\sigma\vert =\sum_{E, \sigma\in E}\sum_{T^E_i, \sigma\in T^E_i} \frac{\vert T^E_i\vert}{3},$$ and after some simple calculations, we see that
$$\vert C_\sigma\vert =\sum_{E, \sigma\in E}\frac{2\vert E\vert}{3\bla{(N-1)}}.$$

The low-order version can easily be shown to be invariant domain preserving under some CFL conditions. For example, the Local Lax-Friedrichs version needs that
$$3\frac{\Delta t}{\vert T^E_i\vert}\rho\big ( \bJ( \bu_{\sigma})\cdot \bn_{b_i}\big )\leq 1, \qquad \forall E, i$$
and this condition turns out to be much more restrictive than the one we can use in practice.

\subsection{Linear stability}
Linear stability is not easy to study when the mesh has no particular symmetry, and when the scheme under consideration does not satisfy a variational principle. In the one dimensional case, it is possible to show that the scheme developed in \cite{Abgrall_AF} is linearly stable: this has been done by W. Barsukow and reported in \cite{BarzukowAbgrall} using \cite{Miller}. However, the scheme \eqref{scheme:pt}, adapted in 1D, would not be exactly the same as in \cite{Abgrall_AF}. So we have conducted an experimental study on several {types} of meshes generated by GMSH \cite{zbMATH05643029} with several {types} of meshing option (frontal Delaunay, Delaunay). The all proved to be stable for the linear and non-linear problems tested above, with a CFL halved with respect to the standard value for the quadratic approximation and divided by 3 for the cubic approximation. This is more or less what is expected because for the quadratic (resp. cubic) approximation, the density of degrees of freedom is approximately doubled (resp. tripled).

For the first-order approximation, using the Rusanov scheme, it is easy to get a CFL condition under which the solution stays in its domain of dependence.
\subsection{Non-linear stabilization}\label{secMOOD}
When the solution develops discontinuities, the high-order scheme will be prone to numerical oscillations. To overcome this, we have used a simplified version of the MOOD paradigm in \cite{Mood1,Vilar}. Since point values and average are independent variables, we need to test both.

The idea is to work with several schemes ranging from order $p=3 \text{ or } 4$ to $1$, with the lowest order one able to provide results staying in the invariance domain. For the element $E$, we write the scheme for $\xbar\bu_E$ as $S_E(k)$ and for $\bu_\sigma$ as $S_\sigma(k)$.

We denote by $\xbar\bu^n_E$ and $\bu_\sigma^n$ the solution at the time $t_n$. After each Runge-Kutta cycle, the updated solution is denoted by $\tilde\bu_E$, $\tilde\bu_\sigma$. 

We first run the full order scheme, for each Runge-Kutta cycle.  Concerning the average $\xbar\bu_E$, we store the flux $\int_{\Gamma}\bbf(\bu^h)\cdot \bn$ for reasons of conservation. Then, for the average,
\begin{enumerate}
\item Computer Admissible Detector (CAD in short): We check if $\tilde\bu_E$ is a valid vector with real components: we check if each component is not NaN. Else, we flag the element and go to the next one in the list.
\item Physically Admissible Detector (PAD in short): We check if $\tilde\bu_E\in \mathcal{D}$. If this is false, the element is flagged, and go to the next one in the list.
\item 
Then we check if at $t_n$, the solution is not a constant in the elements used in the numerical stencil (so we check and compare between the average $\xbar\bu_E$ and point values $\bu_\sigma$ in $E$ with these in the elements sharing a face with $E$). This is done in order to avoid to detect a fake maximum principle. If the solution is declared locally constant, we move to the next element.
\item Discrete Maximum Principle (DMP in short): we check if $\tilde\bu_E$ is a local extrema. if we are dealing with the Euler equations, we compute the density and the pressure and perform this test on these two values only, even though for a system this is not really meaningful. We denote by $\xi$ the functional on which we perform the test (i.e., $\bu$ itself for scalar problems, and the density/pressure obtained from $V(E)$: the set of elements that share a face or a vertex with $E$, excluding $E$ itself).
We say we have a potential extrema if 
$$\xi_E^{n+1}\not\in [\big (\min\limits_{F\in V(E)}\xi_F^n\big )\bla{-}\varepsilon_E^n, \big (\max\limits_{F\in V(E)}\xi_F^n\big )\bla{+}\varepsilon_E^n],$$
where $\varepsilon_E^n$ is estimated as in \cite{Mood1}.
If the test is wrong, we go to the next element, else the element is flagged.
\end{enumerate}
If an element is flagged, then each of its faces are flagged, and re-compute the flux of the flagged faces using the first-order scheme. 

For the point values, the procedure is similar, and then for the flagged degrees of freedom, we re-compute
$\sum\limits_{E, \sigma\in E}\Phi_\sigma^E(\bu^n)$ with the first-order scheme.

\section{Numerical results}
In this section, we demonstrate the performance of the developed schemes in \S\ref{sec3} on several numerical examples.
When MOOD is activated, we use the local Rusanov scheme (with a small stencil) as described in \S\ref{loworder}. The CFL number (based on the elements) is always set to $0.4$.
\subsection{Scalar case}
We begin with the scalar case, 
\begin{equation}\label{scalar:test}\dpar{u}{t}+\mathbf{a}(\bx)\cdot
\nabla u=0,
\end{equation}
with, for $\bx=(x,y)$, $\ba(\bx)=\big (a^{(1)},a^{(2)}\big )$ representing the velocity field that the quantity $u$ is moving with.
\subsubsection{Convection problem}
In the first example, we consider
$\ba(\bx)=2\pi\big (-(y-y_0),x-x_0\big )$ \bla{and} note that $\text{ div
}\ba=0$, so that \eqref{scalar:test} can be put in a conservative form.
The domain is $\Omega=[-20,20]^2$, $(x_0,y_0)=(0,0)$ and the initial
condition is
\begin{equation}\label{test411}
    u_0(\bx)=e^{-\Vert\mathbf{x}-\mathbf{\tilde{x}}_0\Vert^2}, ~~\mathbf{\tilde{x}}_0=(5, 5).
\end{equation}
The solution is checked after one rotation ($T=1$) with the exact one.
The initial mesh is obtained with GMSH, using the front Delaunay option.
Then we sub-divide the mesh by cutting the edges into two equal or three equal segments,
and a quadratic approximation or a cubic approximation is employed. The discrete $L^1$-, $L^2$-, and $L^{\infty}$-errors of the point values $u$ computed with quadratic and cubic approximations are shown in
Tables \ref{Table:1} and \ref{Table:2}, respectively. As one can clearly see that, the expected third- and fourth-order accuracy are achieved.
\begin{table}[ht!]
 \caption{
\label{Table:1} Errors and rates on the point values for the initial condition \eqref{test411} with the quadratic approximation of \eqref{scheme:fv}-\eqref{scheme:pt}.}
     \begin{center}
         \begin{tabular}{|c||c|c|c|c|c|c|c|}
             \hline
             $h $& $L^1$&$\text{rate}$& $L^2$ &$\text{rate}$&$L^\infty$ & $\text{rate}$\\
             \hline
             $2.00\times10^{-1}$   &     $ 4.309\times10^{-1}$ &-       &   $6.703\times 10^{-2}$  &-      & $3.270\times10^{-2}$  &-\\
             $1.91\times 10^{-1}$   &     $3.749\times10^{-1}$  &2.85  &   $5.857\times 10^{-2}$  &2.77 & $2.871\times10^{-2}$  &2.66\\
             $1.82\times10^{-1}$   &     $3.241\times10^{-1}$  &3.13  &   $5.079\times10^{-2}$  &3.06 & $2.512\times10^{-2}$  &2.87\\
             $1.74\times10^{-1}$   &     $2.850\times10^{-1}$  &2.89 &   $4.485\times10^{-2}$  &2.80 & $2.215\times10^{-2}$  &2.83\\
             $1.67\times10^{-1}$   &     $2.553\times10^{-1}$  &2.59  &   $4.026\times10^{-2}$  &2.54 & $2.007\times10^{-2}$  &2.39\\
             $1.60\times10^{-1}$   &     $2.269\times10^{-1}$  &2.89  &   $3.584\times10^{-2}$  &2.84 & $1.786\times10^{-2}$  &2.78\\          
             $1.33\times10^{-1}$   &     $1.317\times10^{-1}$  &2.98  &   $2.100\times10^{-2}$  &2.93 & $1.056\times10^{-2}$  &2.88\\
             \hline
         \end{tabular}
     \end{center}
\end{table}

\begin{table}[ht!]
\caption{
\label{Table:2} Errors and rates on the point values for the initial condition \eqref{test411} with the cubic approximation of \eqref{scheme:fv}-\eqref{scheme:pt}.}
	\begin{center}
		\begin{tabular}{|c||c|c|c|c|c|c|c|}
			\hline
			$h $& $L^1$&$\text{rate}$& $L^2$ &$\text{rate}$& $L^\infty$ & $\text{rate}$\\
			\hline
			$4.00\times10^{-1}$   &     $5.219\times10^{-1}$ &-       &    $6.800\times10^{-2}$  &-      & $2.206\times10^{-2}$  &-\\
			$3.33\times10^{-1}$   &     $2.440\times10^{-1}$ &4.17  &    $3.268\times10^{-2}$  &4.02 & $1.121\times10^{-2}$  &3.71\\
			$2.86\times10^{-1}$   &     $1.358\times10^{-1}$ &3.80  &    $1.839\times10^{-2}$  &3.73 & $6.339\times10^{-3}$  &3.70\\
			$2.50\times10^{-1}$   &     $7.978\times10^{-2}$ &3.99  &    $1.083\times10^{-2}$  &3.96 & $3.705\times10^{-3}$  &4.02\\
			$2.22\times10^{-1}$   &     $4.868\times10^{-2}$ &4.19  &    $6.634\times10^{-3}$  &4.16 & $2.321\times10^{-3}$  &3.97\\			
			$2.00\times10^{-1}$   &     $3.267\times10^{-2}$ &4.19  &    $4.458\times10^{-3}$  &4.16 & $1.542\times10^{-3}$  &3.97\\
			\hline
		\end{tabular}
	\end{center}
\end{table}

\subsubsection{Zalesak test case}
In the second example, we consider a Zalesak’s problem, which involves the solid body rotation of a notched disk.
   Here the domain is $[0,1]^2$ and the rotation field is $\ba(\bx)$ also defined to be  a rotation with respect to $(x_0,y_0)=(0.5,0.5)$ and angular speed of $2\pi$.
   The initial condition is 
   $$u_0(\bx)=\left \{
   \begin{array}{ll}
   0.25\big ( 1+\cos(\frac{\pi r_{1}}{0.15})) & \text{ if }r_{1}=\Vert\bx-(0.25, 0.5)\Vert\leq 0.15,\\
  1-\frac{r_2}{0.15}&\text{ if } r_{2}=\Vert\bx-(0.5,0.25)\Vert\leq 0.15, \\
  1&\text{ if } r_{3}=\Vert\bx-(0.5,0.75)\Vert\leq 0.15,
  \\
  0& \text{ if }|x-0.5| \leq 0.025 ~\text{ and }~y\in[0.6,0.85]. 
  \end{array}\right .$$
  
\begin{figure}[ht!]
\centerline{\subfigure[point values $\bu_\sigma$]{\includegraphics[trim=1.2cm 2.2cm 3.5cm 1.8cm,clip,width=3.2cm]{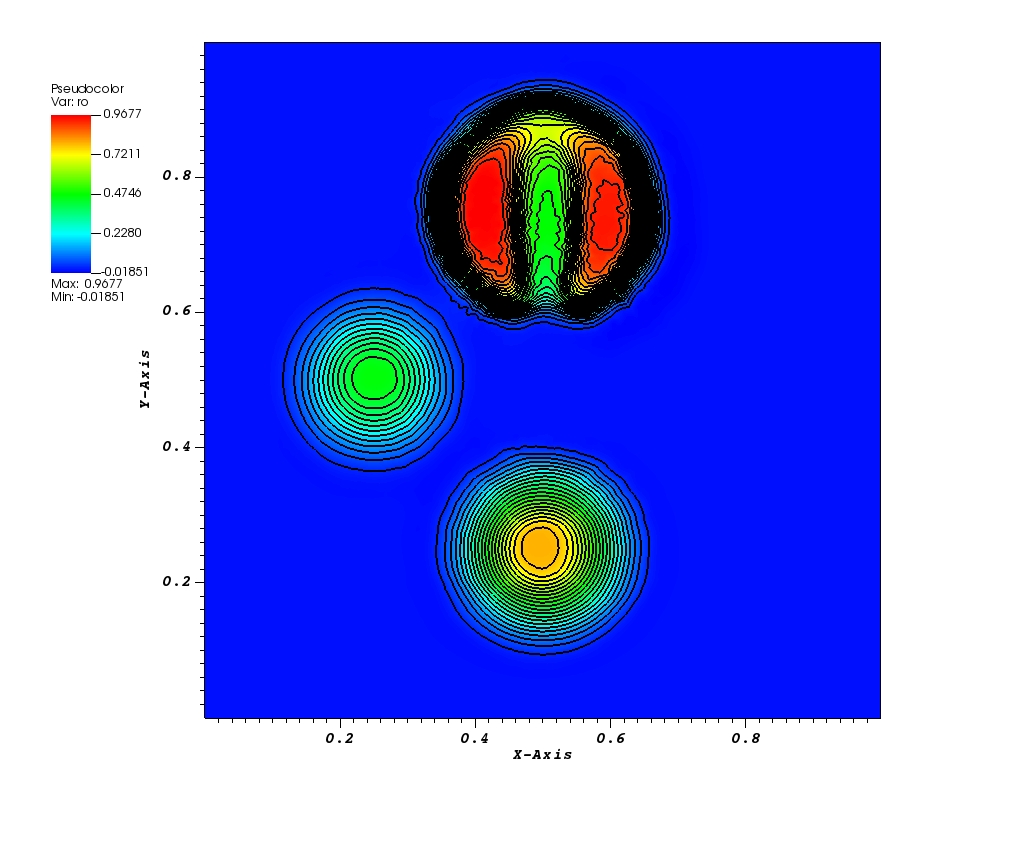}}\hspace*{0.002cm}
\subfigure[averages $\xbar{\bu}_E$]{\includegraphics[trim=1.2cm 2.2cm 3.5cm 1.8cm,clip,width=3.2cm]{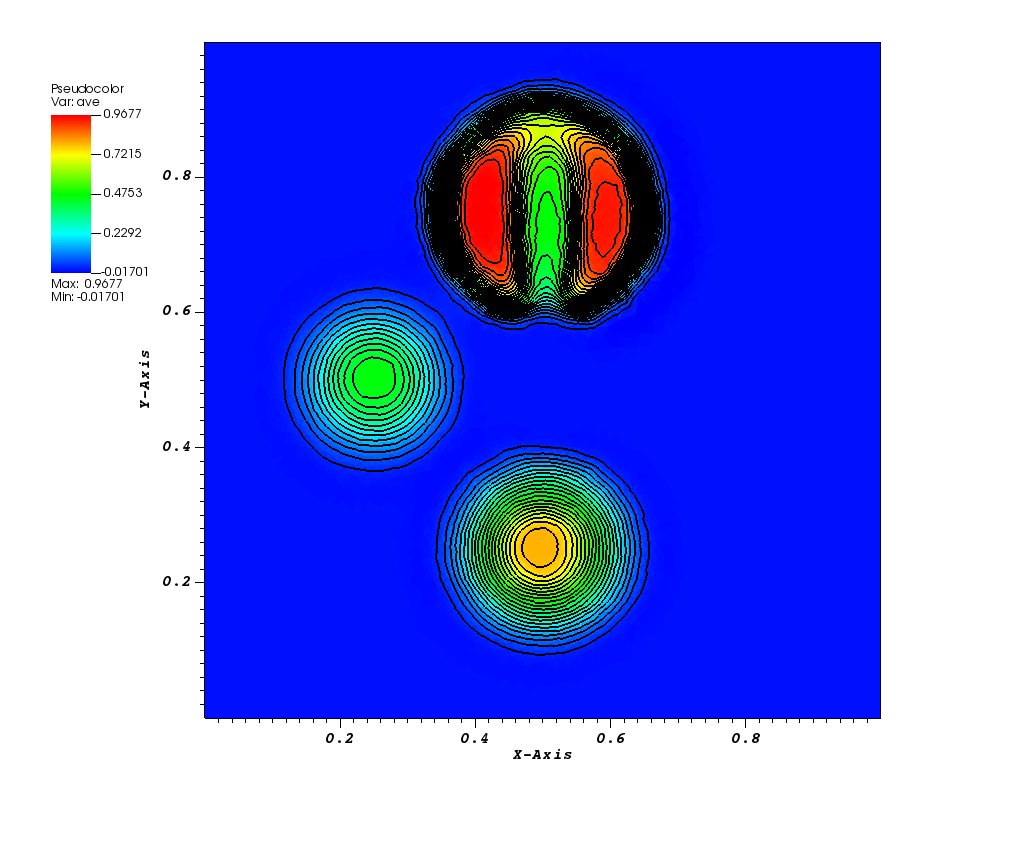}}\hspace*{0.002cm}
\subfigure[flag on $\bu_\sigma$]{\includegraphics[trim=1.2cm 2.2cm 3.3cm 1.8cm,clip,width=3.2cm]{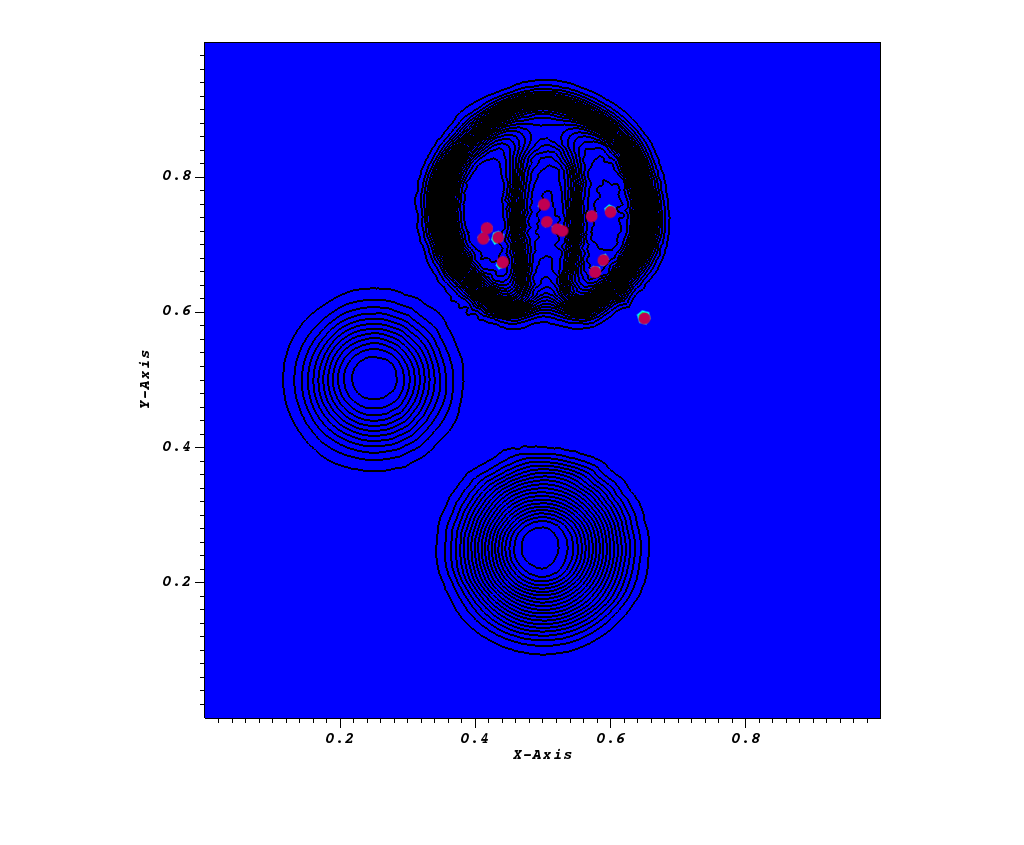}}\hspace*{0.02cm}
\subfigure[flag on $\xbar{\bu}_E$]{\includegraphics[trim=1.2cm 2.2cm 3.3cm 1.8cm,clip,width=3.2cm]{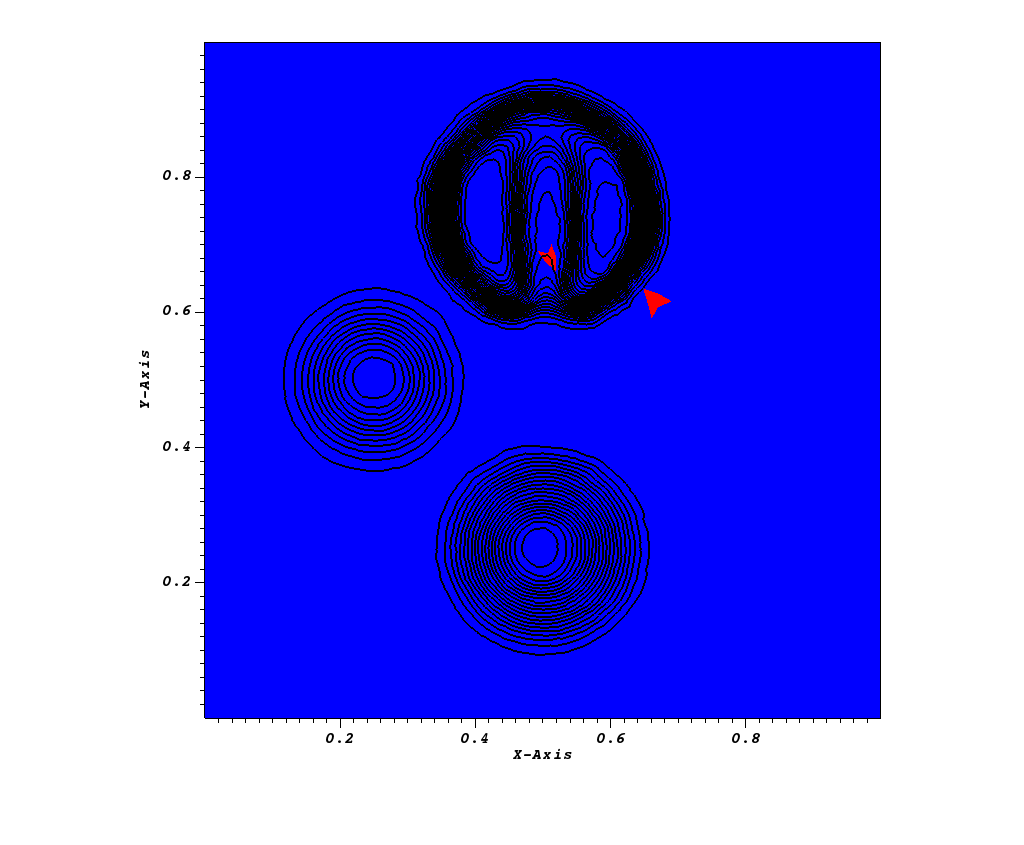}}}
\vskip5pt
\centerline{\subfigure[point values $\bu_\sigma$]{\includegraphics[trim=1.2cm 2.2cm 2.5cm 1.8cm,clip,width=3.1cm]{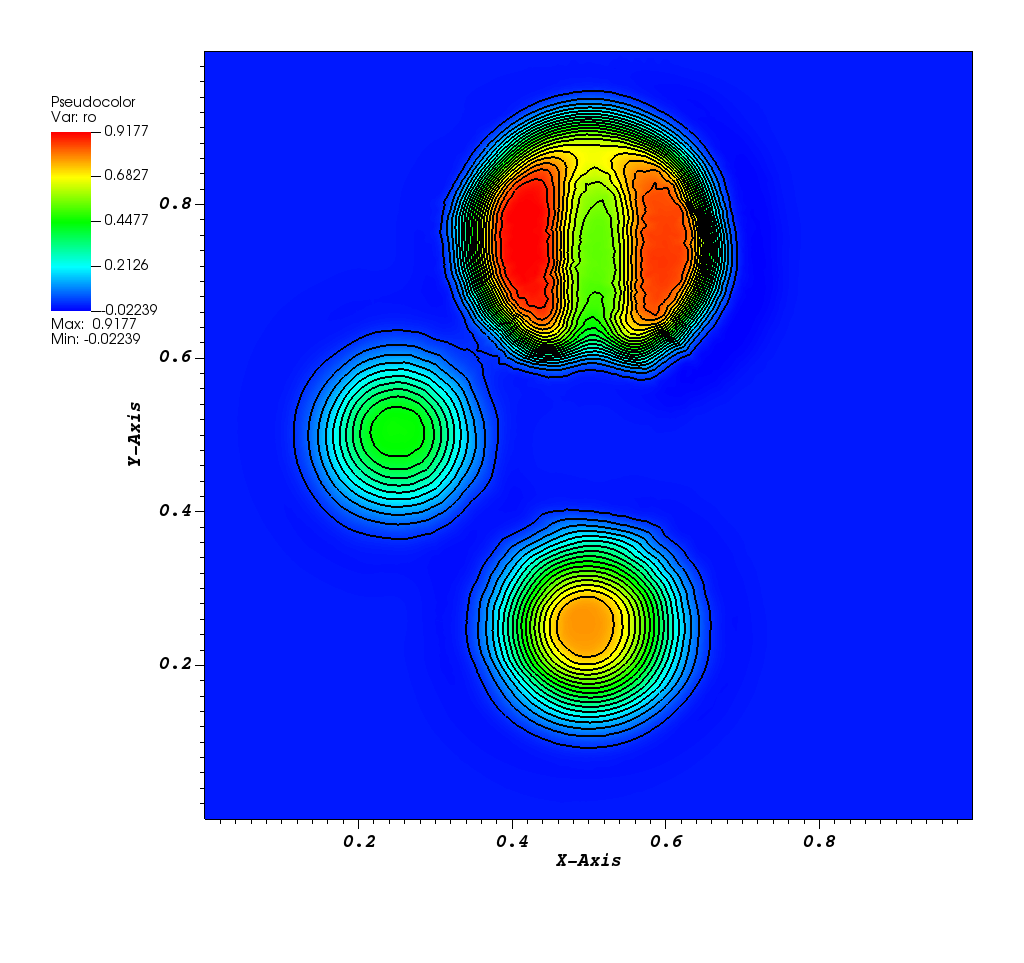}}\hspace*{0.002cm}
\subfigure[averages $\xbar{\bu}_E$]{\includegraphics[trim=1.2cm 2.2cm 2.5cm 1.8cm,clip,width=3.1cm]{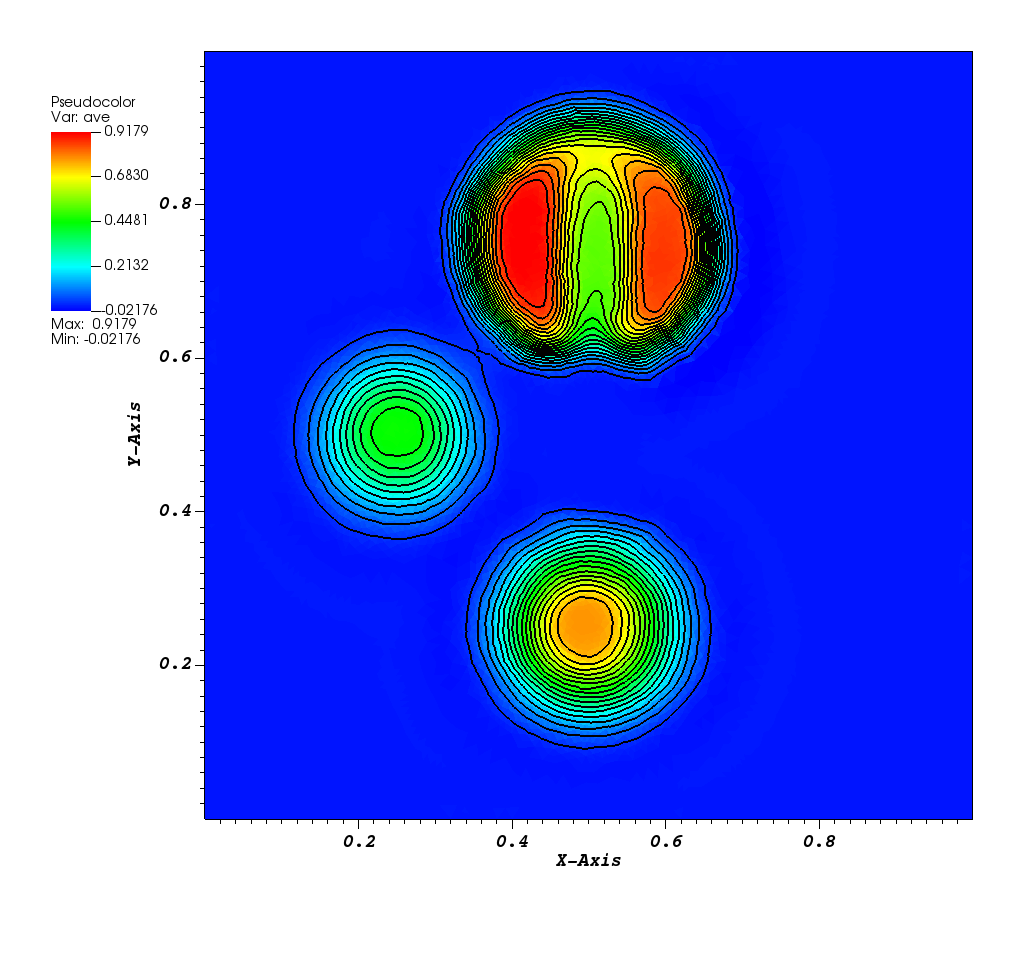}}\hspace*{0.002cm}
\subfigure[flag on $\bu_\sigma$]{\includegraphics[trim=1.2cm 2.2cm 2.5cm 1.8cm,clip,width=3.1cm]{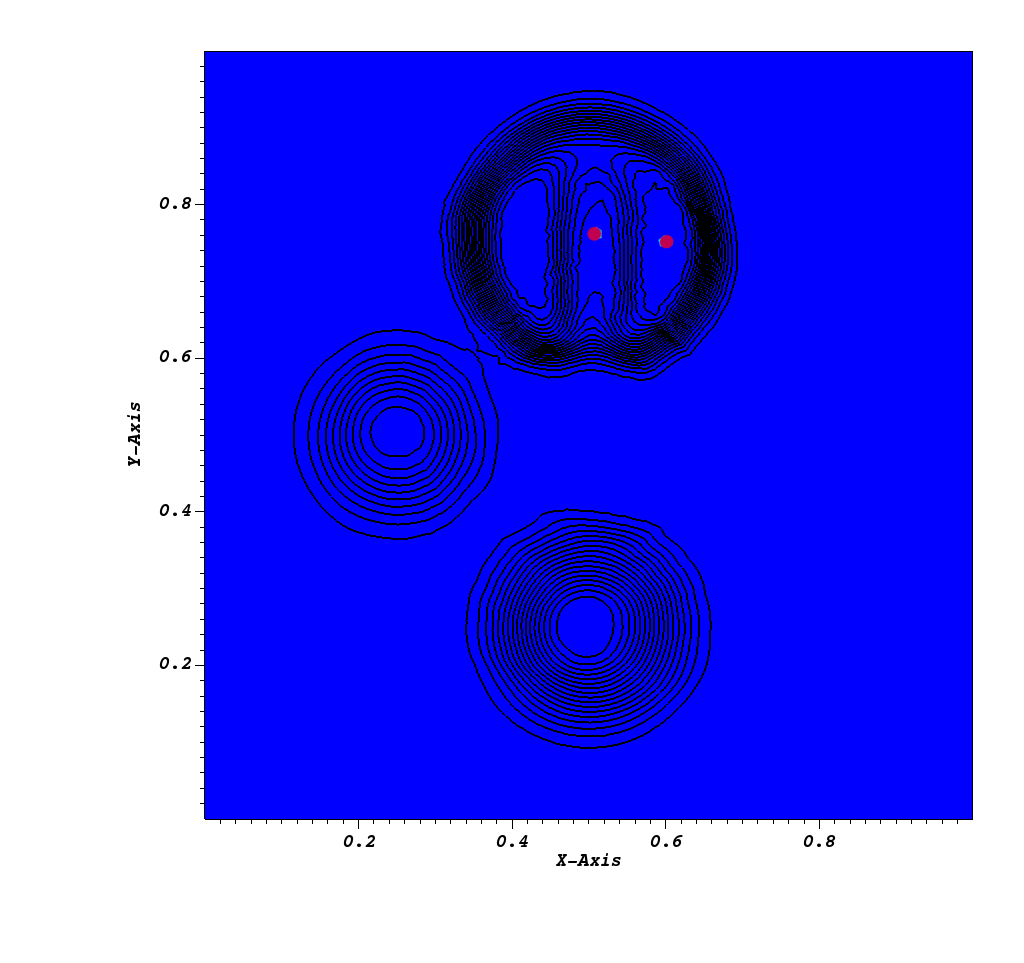}}\hspace*{0.002cm}
\subfigure[flag on $\xbar{\bu}_E$]{\includegraphics[trim=1.2cm 2.2cm 2.5cm 1.8cm,clip,width=3.1cm]{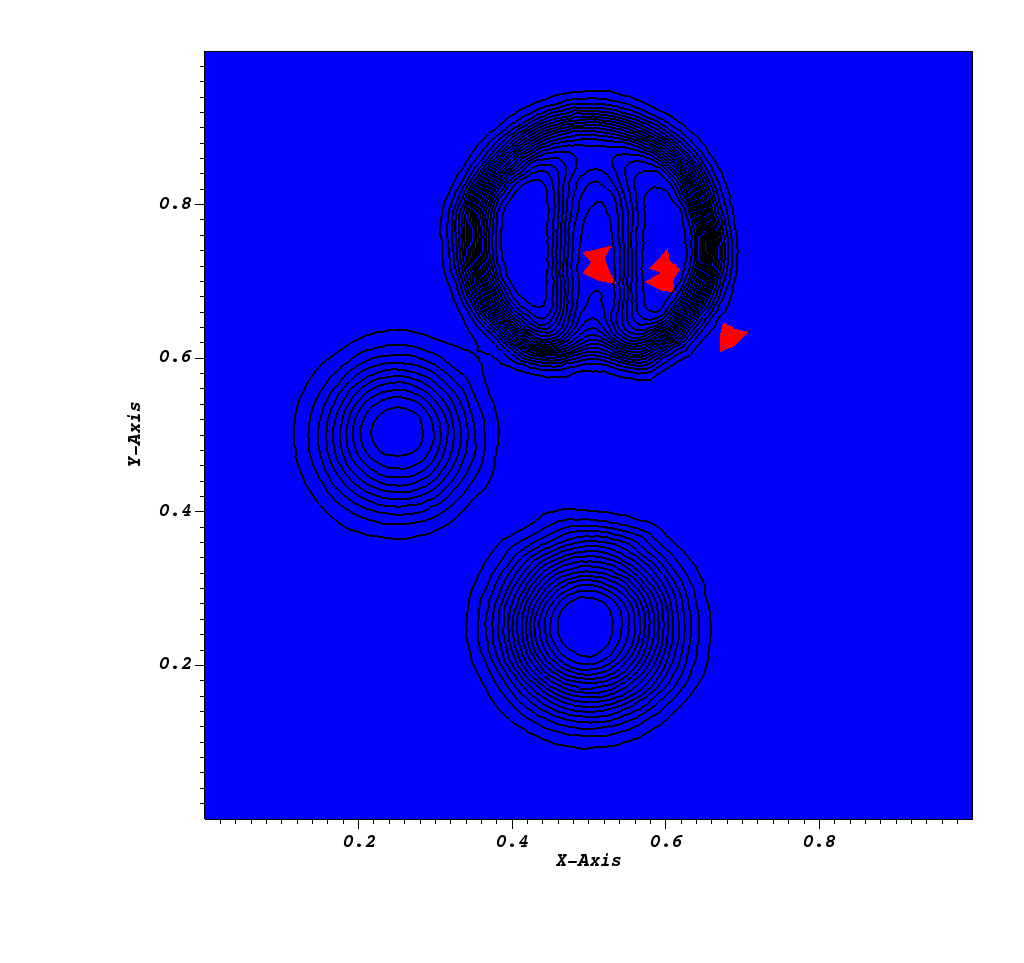}}}
\vskip5pt
\centerline{\subfigure[point values $\bu_\sigma$]{\includegraphics[trim=1.2cm 2.2cm 2.5cm 1.8cm,clip,width=3.1cm]{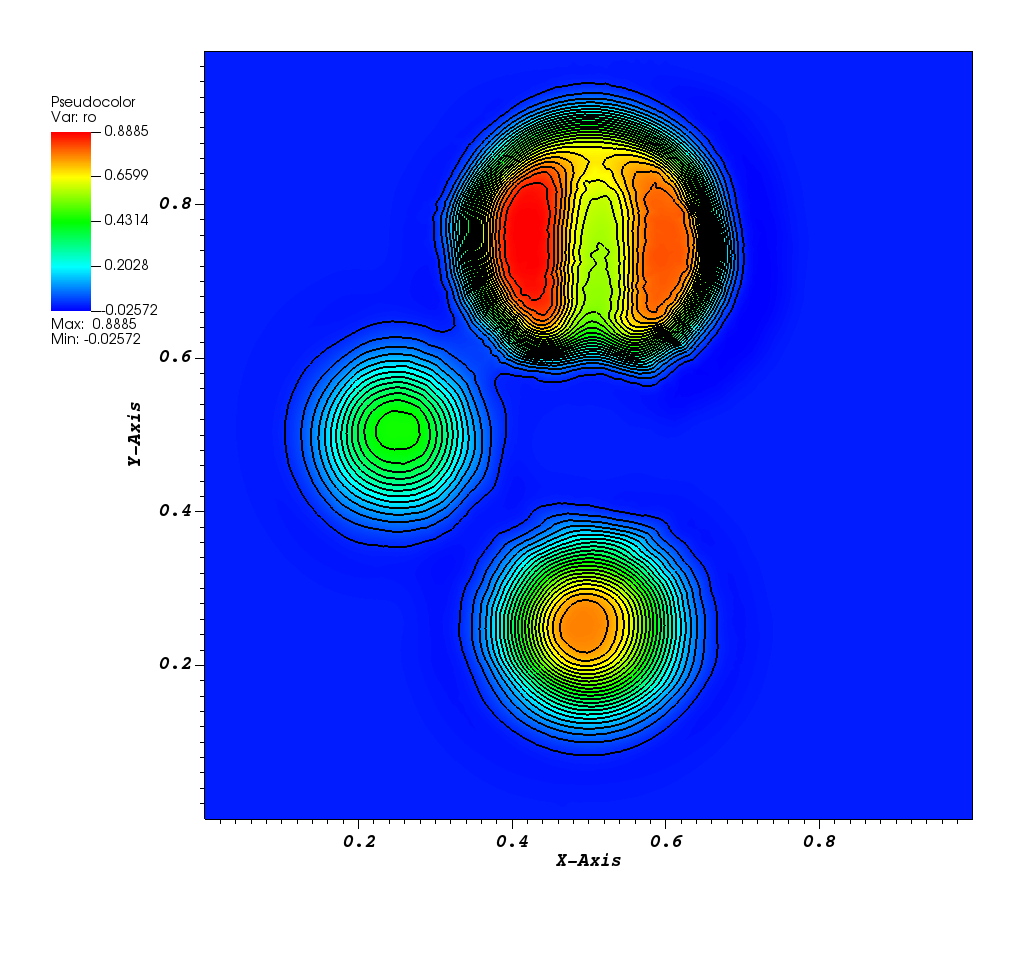}}\hspace*{0.002cm}
\subfigure[averages $\xbar{\bu}_E$]{\includegraphics[trim=1.2cm 2.2cm 2.5cm 1.8cm,clip,width=3.1cm]{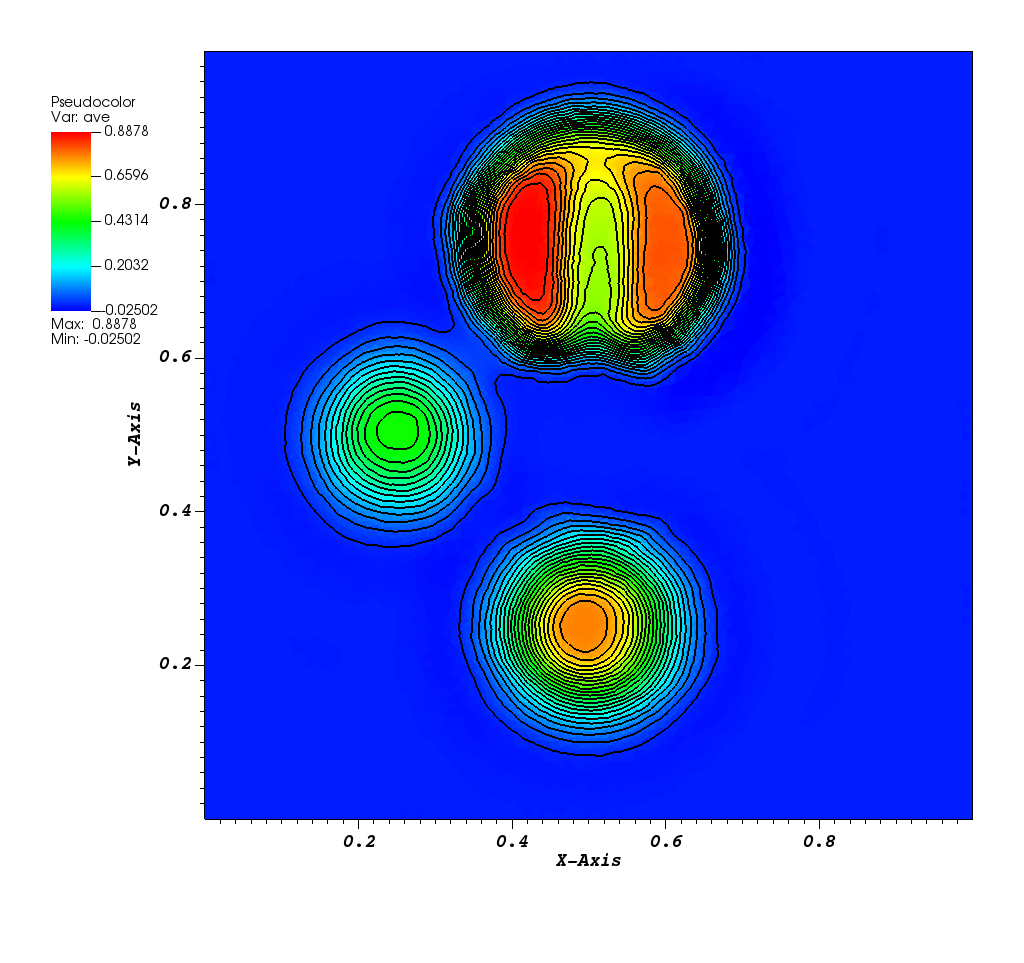}}\hspace*{0.002cm}
\subfigure[flag on $\bu_\sigma$]{\includegraphics[trim=1.2cm 2.2cm 2.5cm 1.8cm,clip,width=3.1cm]{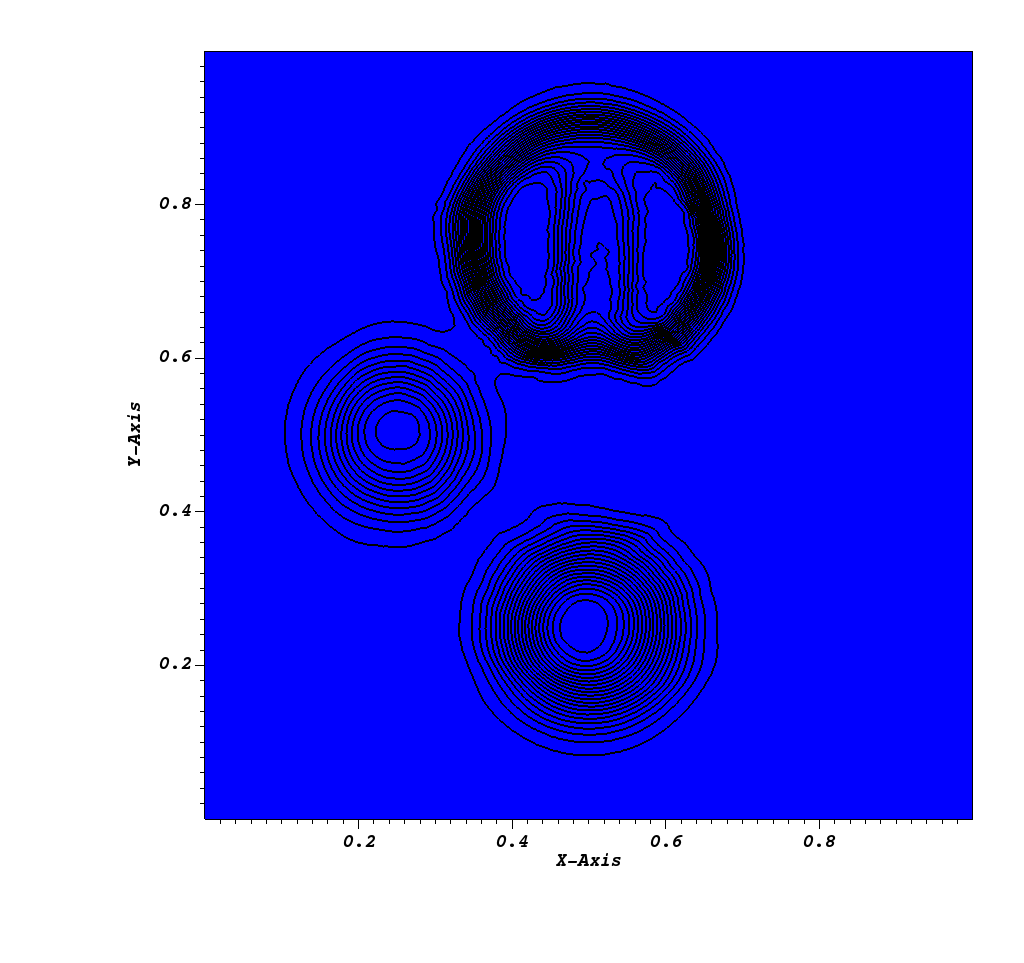}}\hspace*{0.002cm}
\subfigure[flag on $\xbar{\bu}_E$]{\includegraphics[trim=1.2cm 2.2cm 2.5cm 1.8cm,clip,width=3.1cm]{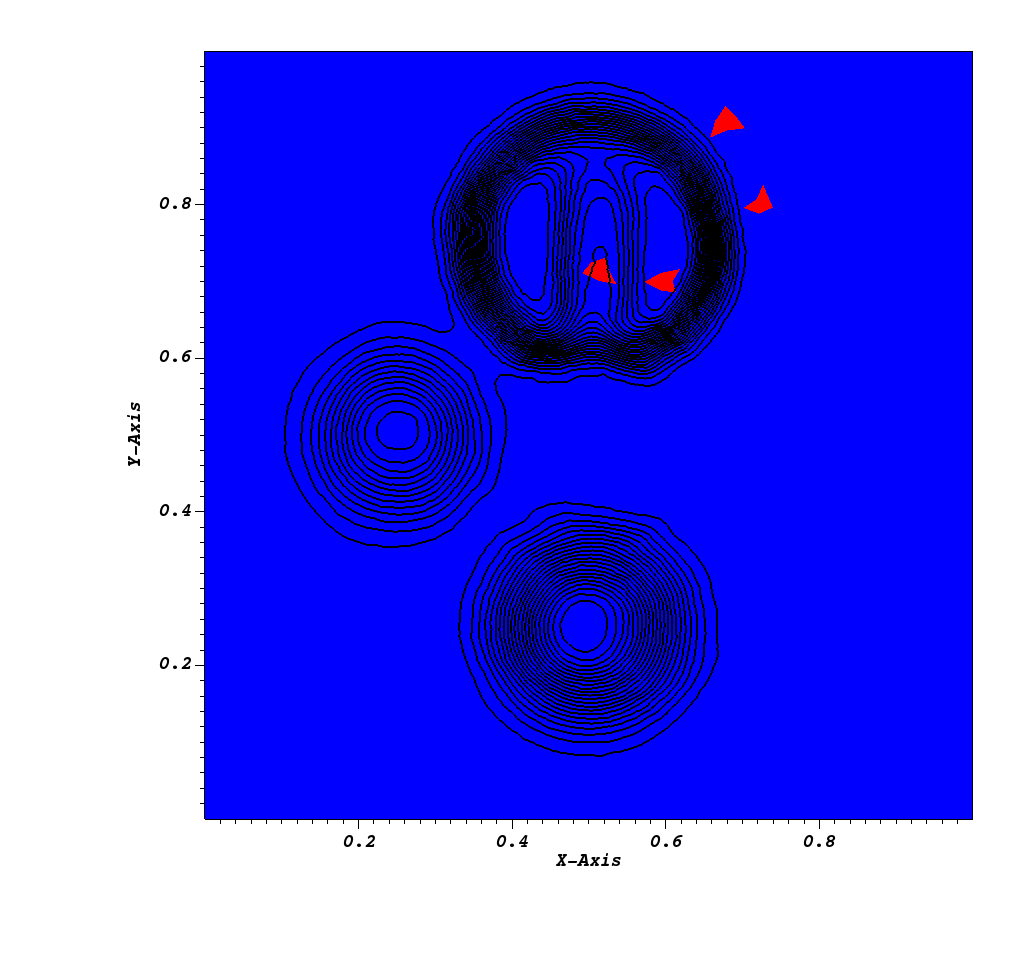}}}
\caption{Zalesak problem: quadratic approximation and MOOD criteria are checked expect for PAD. Top row: $T=1$; Middle row: $T=2$; Bottom row: $T=3$.\label{fig:zakesak_p2_mood1}}
\end{figure}

\begin{figure}[ht!]
\centerline{\subfigure[point values $\bu_\sigma$]{\includegraphics[trim=1.2cm 2.2cm 2.5cm 1.8cm,clip,width=3.2cm]{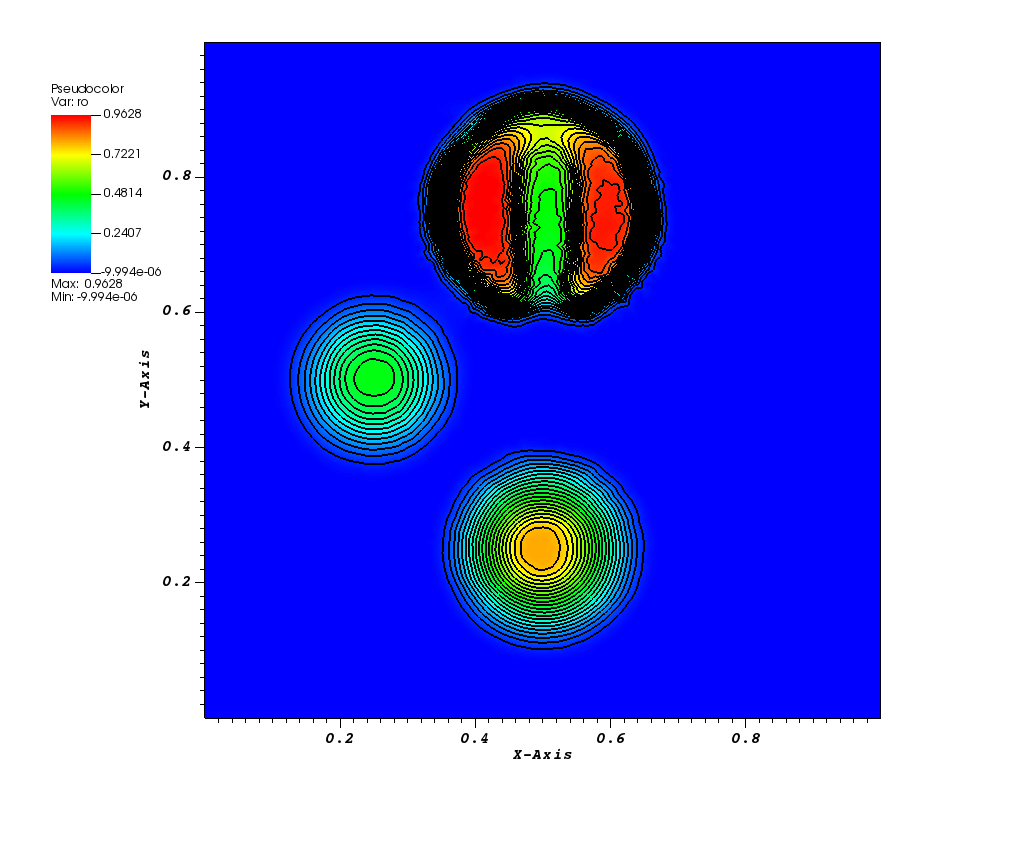}}\hspace*{0.02cm}
\subfigure[averages $\xbar{\bu}_E$]{\includegraphics[trim=1.2cm 2.2cm 2.5cm 1.8cm,clip,width=3.2cm]{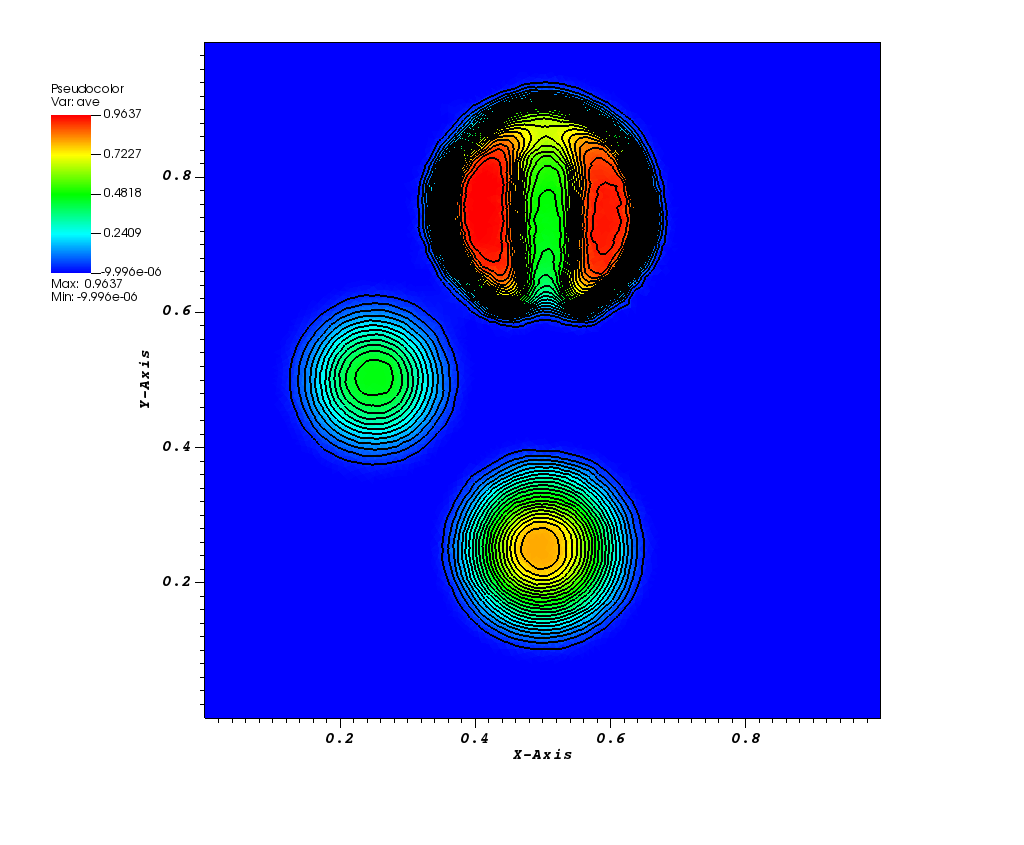}}\hspace*{0.02cm}
\subfigure[flag on $\bu_\sigma$]{\includegraphics[trim=1.2cm 2.2cm 2.5cm 1.8cm,clip,width=3.2cm]{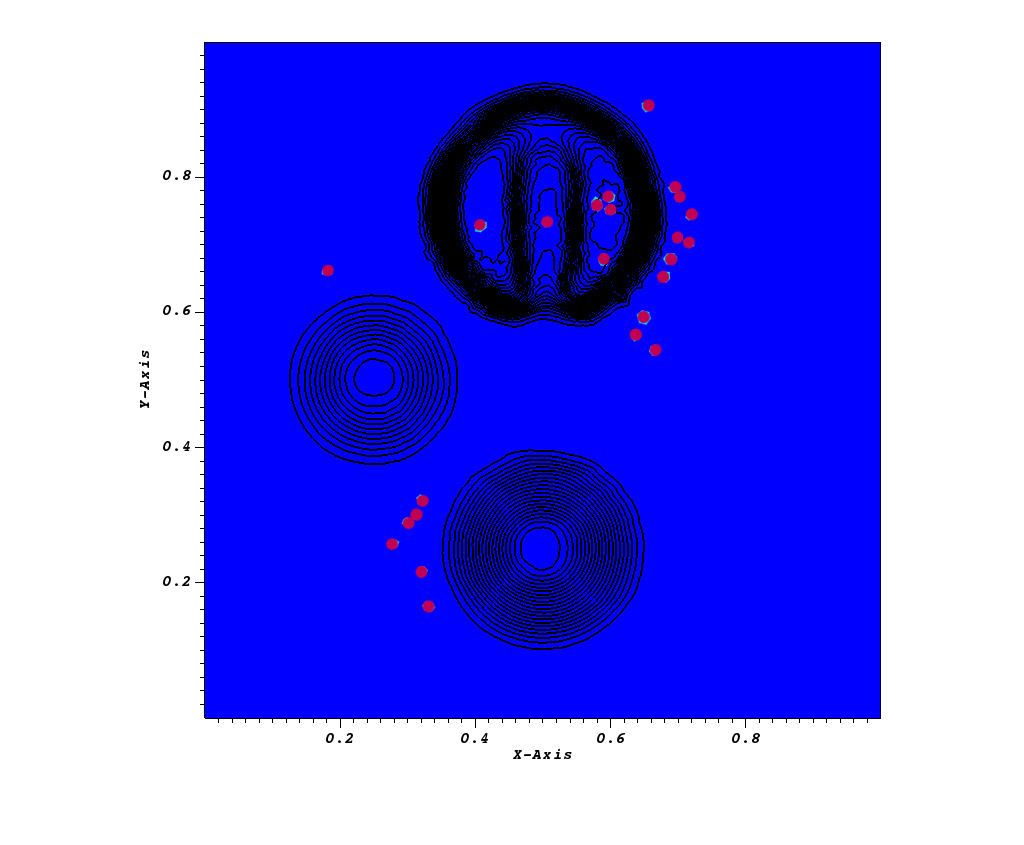}}\hspace*{0.02cm}
\subfigure[flag on $\xbar{\bu}_E$]{\includegraphics[trim=1.2cm 2.2cm 2.5cm 1.8cm,clip,width=3.2cm]{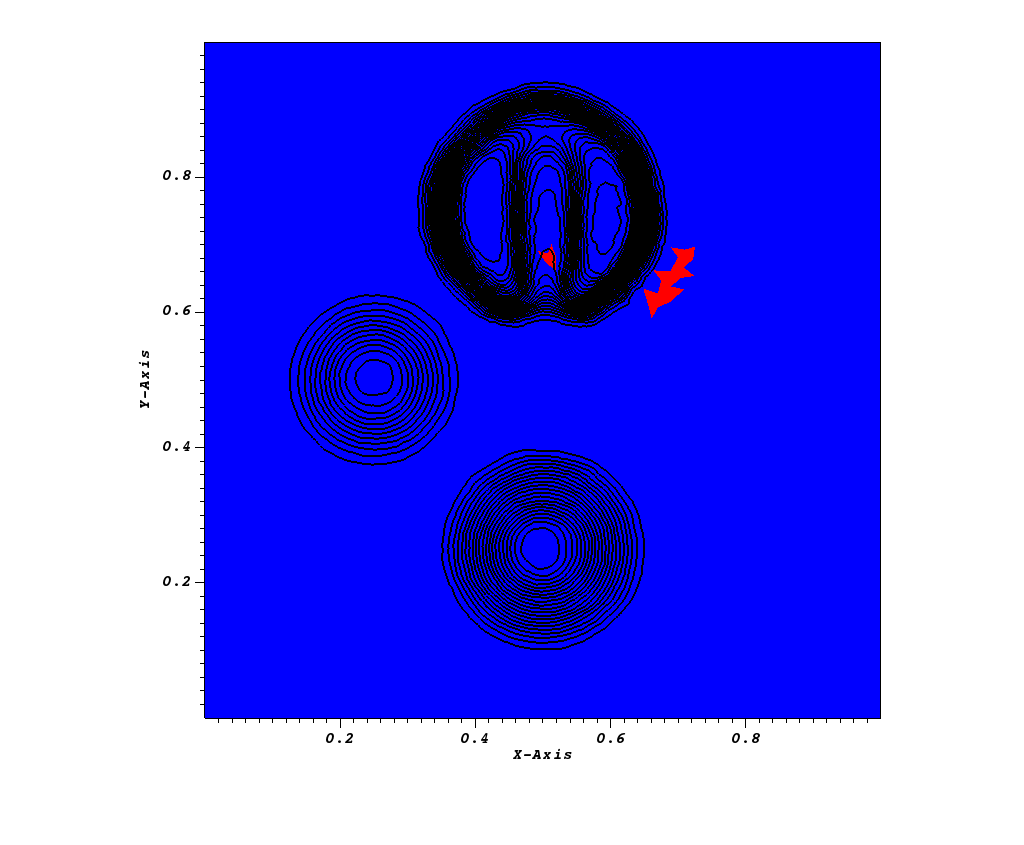}}}
\vskip5pt
\centerline{\subfigure[point values $\bu_\sigma$]{\includegraphics[trim=1.2cm 2.2cm 2.5cm 1.8cm,clip,width=3.2cm]{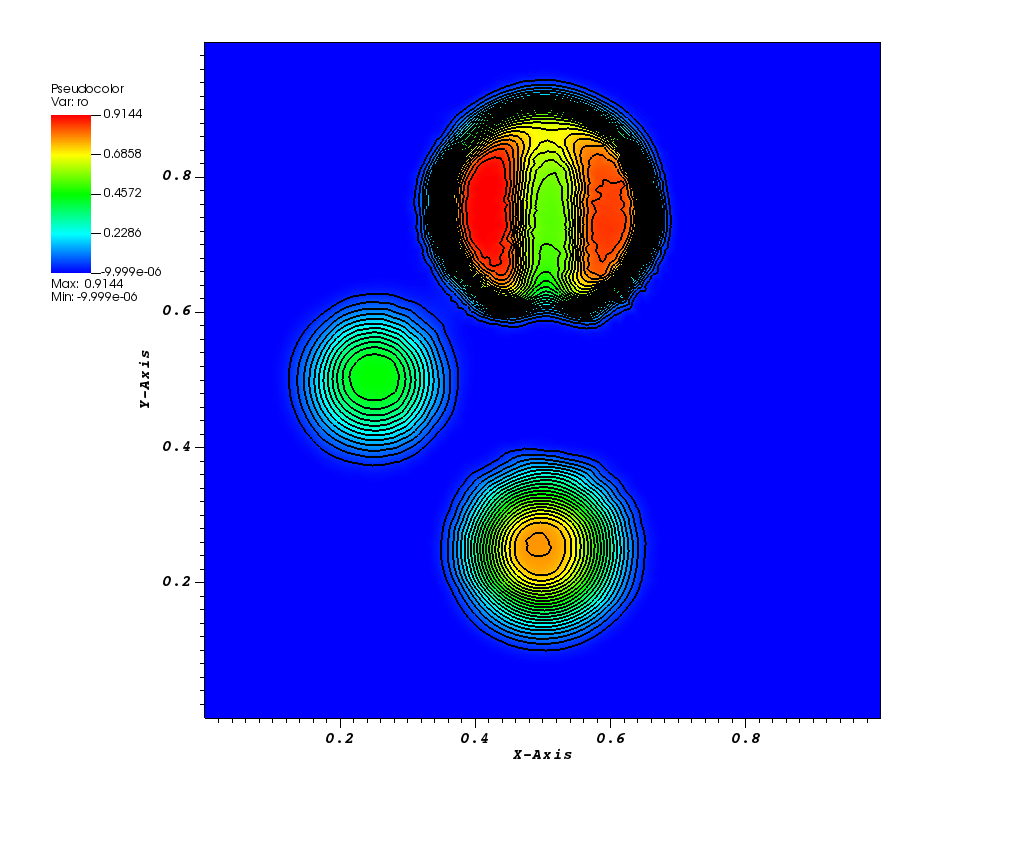}}\hspace*{0.02cm}
\subfigure[averages $\xbar{\bu}_E$]{\includegraphics[trim=1.2cm 2.2cm 2.5cm 1.8cm,clip,width=3.2cm]{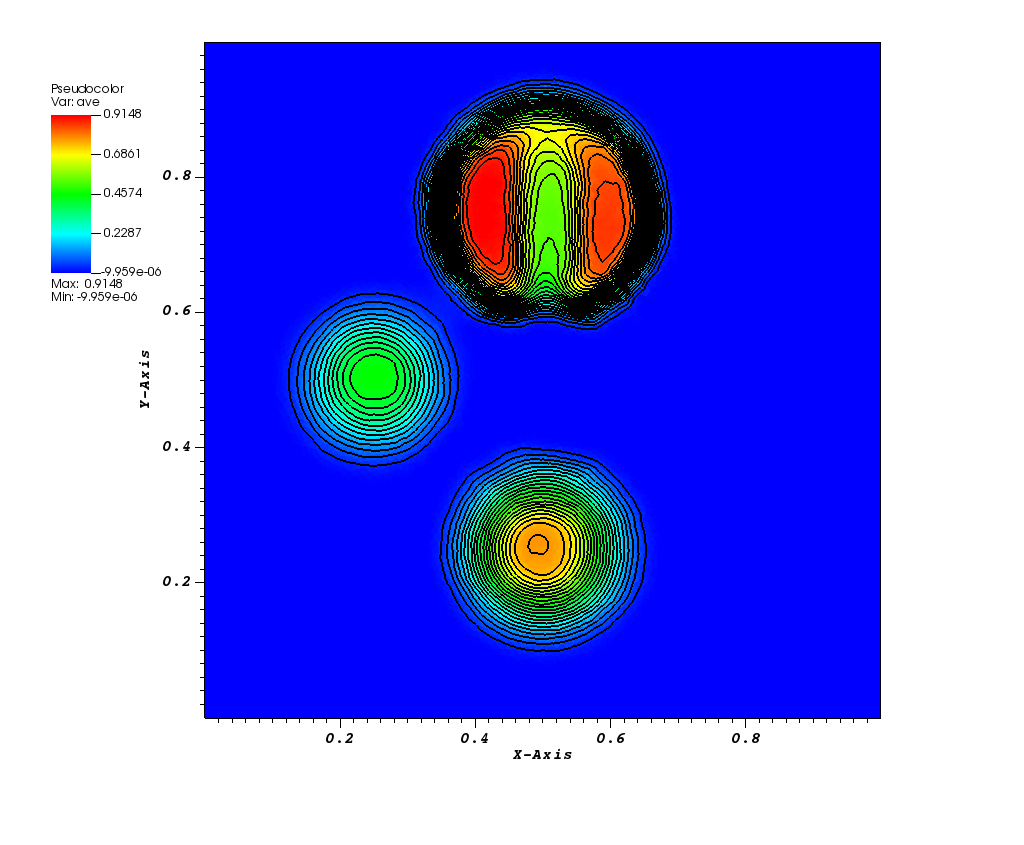}}\hspace*{0.02cm}
\subfigure[flag on $\bu_\sigma$]{\includegraphics[trim=1.2cm 2.2cm 2.5cm 1.8cm,clip,width=3.2cm]{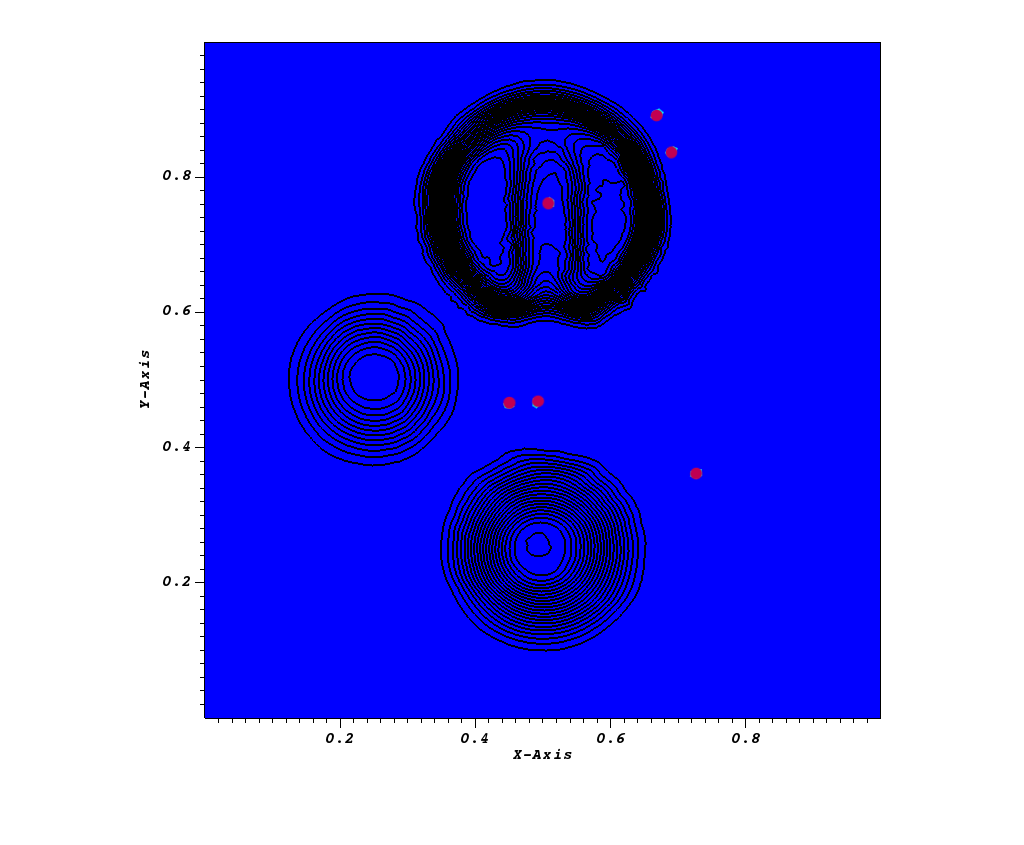}}\hspace*{0.02cm}
\subfigure[flag on $\xbar{\bu}_E$]{\includegraphics[trim=1.2cm 2.2cm 2.5cm 1.8cm,clip,width=3.2cm]{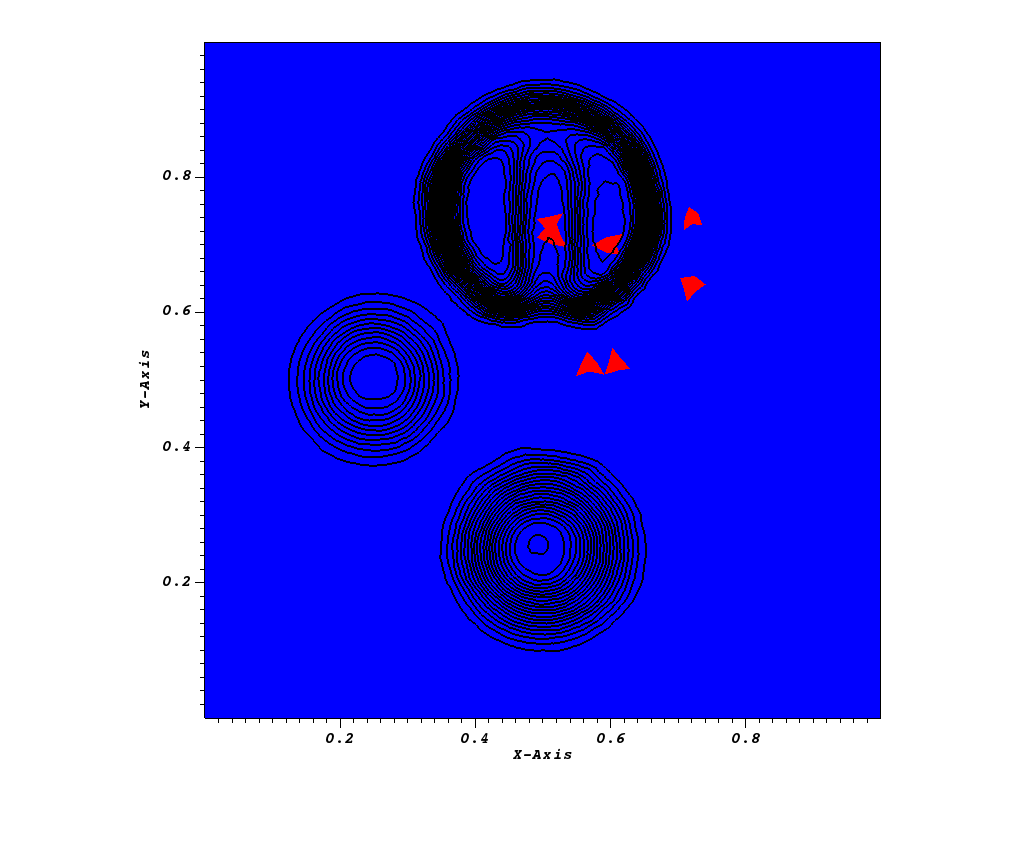}}}
\vskip5pt
\centerline{\subfigure[point values $\bu_\sigma$]{\includegraphics[trim=1.2cm 2.2cm 2.5cm 1.8cm,clip,width=3.2cm]{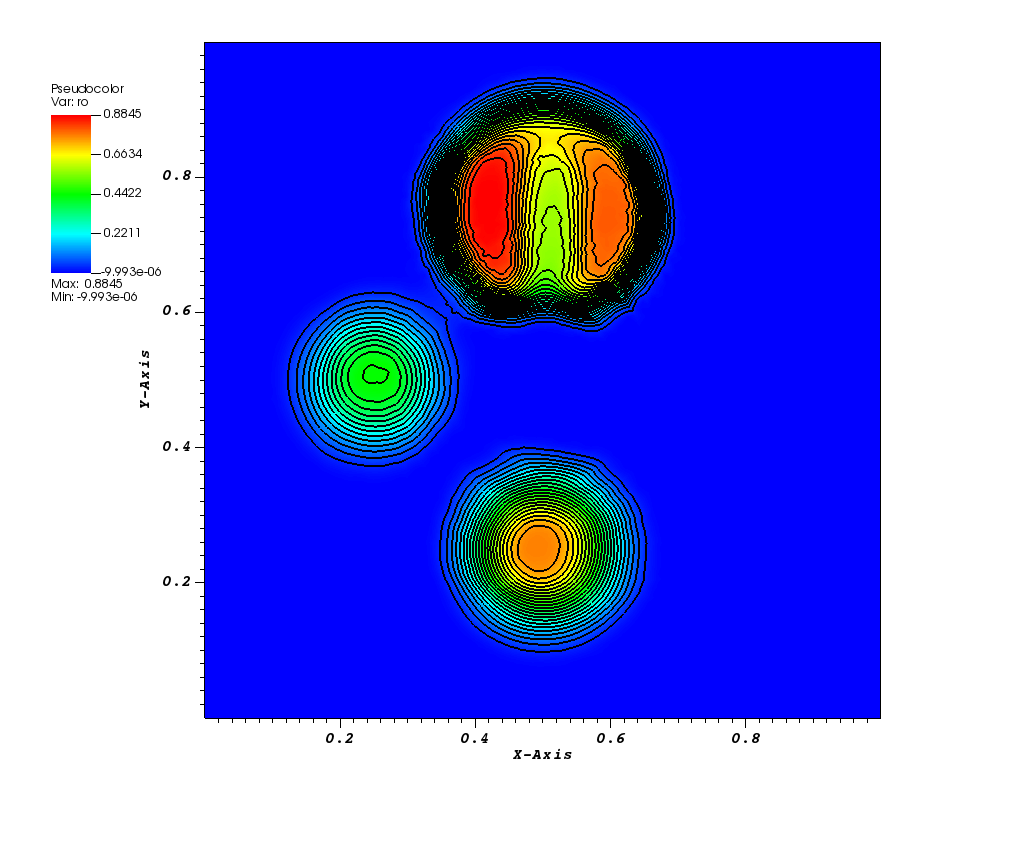}}\hspace*{0.02cm}
\subfigure[averages $\xbar{\bu}_E$]{\includegraphics[trim=1.2cm 2.2cm 2.5cm 1.8cm,clip,width=3.2cm]{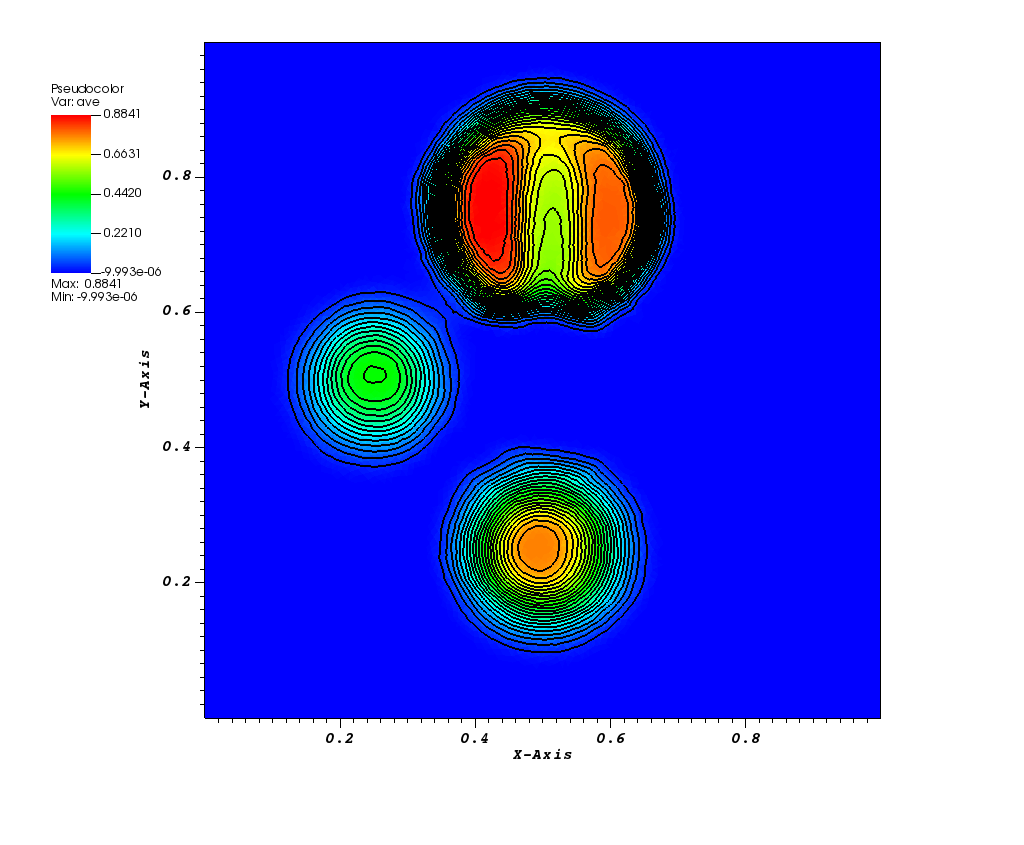}}\hspace*{0.02cm}
\subfigure[flag on $\bu_\sigma$]{\includegraphics[trim=1.2cm 2.2cm 2.5cm 1.8cm,clip,width=3.2cm]{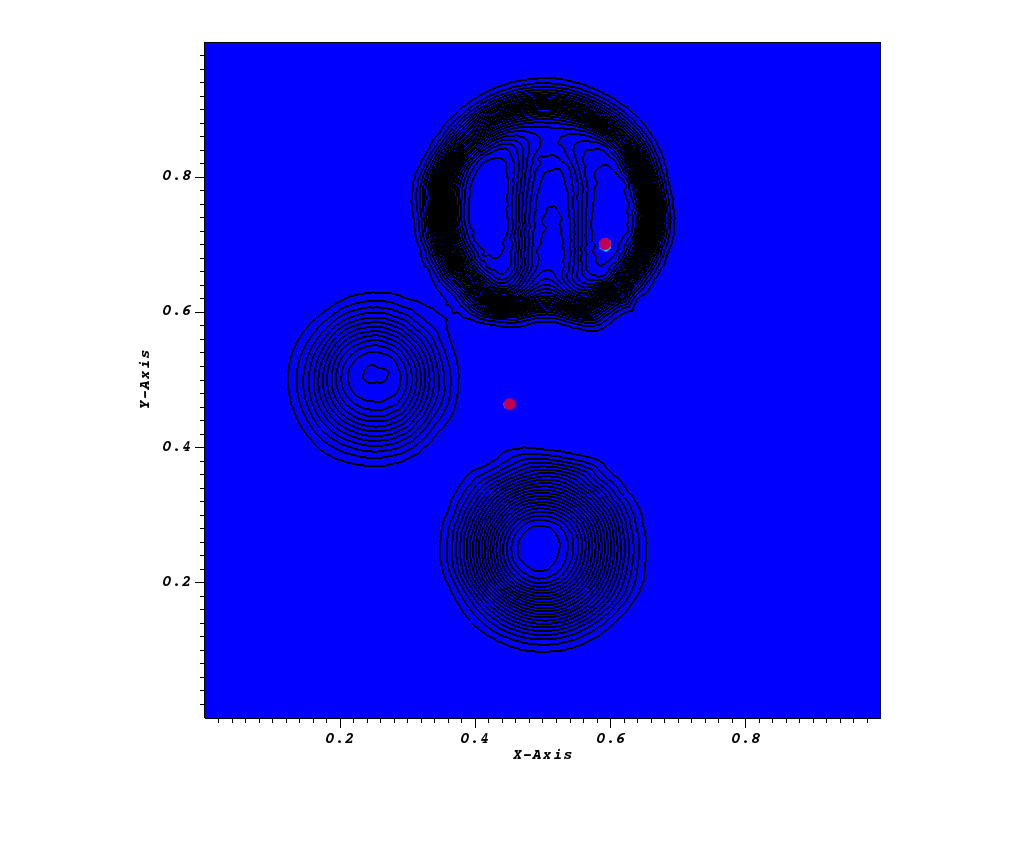}}\hspace*{0.02cm}
\subfigure[flag on $\xbar{\bu}_E$]{\includegraphics[trim=1.2cm 2.2cm 2.5cm 1.8cm,clip,width=3.2cm]{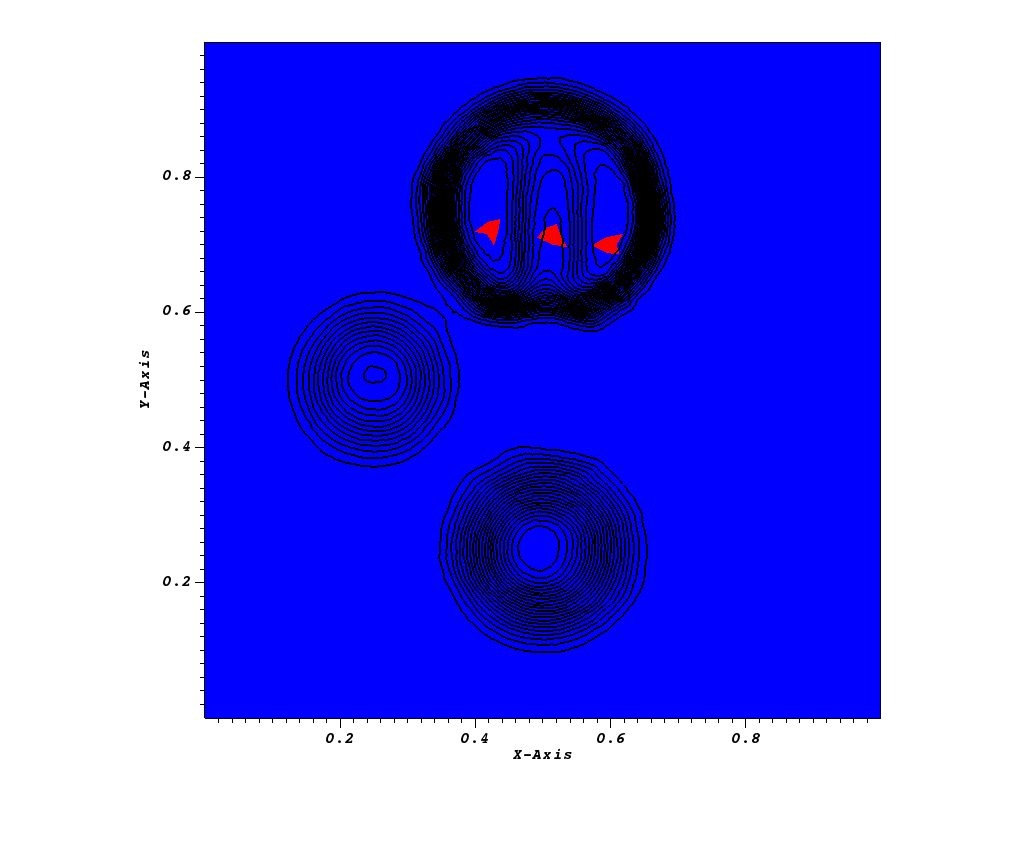}}}
\caption{Zalesak problem: quadratic approximation and all MOOD criteria are checked. Top row: $T=1$; Middle row: $T=2$; Bottom row: $T=3$.\label{fig:zakesak_p2_mood2}}
\end{figure}

We use the MOOD procedure and compute the numerical solution with a mesh that has 3442 vertices, 10123 edges, and 16881 elements (hence 13565 (resp. 40569) DoFs in total for the quadratic (resp. cubic) approximation, thanks to the Euler relation). The results computed by the quadratic approximation after 1, 2, and 3 rotations are displayed in Figures \ref{fig:zakesak_p2_mood1}--\ref{fig:zakesak_p2_mood2}.
From the sub-plots in the first two columns of these figures, one can observe that, both the point values and the average values can be correctly captured. In the last two columns of these figures, we also indicate in red where the MOOD criteria are violated. We first check the criteria expect for PAD. From the sub-plots in the last two columns of Figure \ref{fig:zakesak_p2_mood1}, we see that, as expected, the first-order scheme is only used around the region where the solution is non-smooth. We also observe that the minimum values of $\bu_\sigma$ and $\xbar{\bu}_E$ are $\approx-10^{-2}$. This is due to the fact that the PAD is not checked. 

We then check all the MOOD criteria including PAD, and say that PAD is violated if $\bu_\sigma<-10^{-5}$ or $\xbar{\bu}_E<-10^{-5}$. From the sub-plots in the last two columns of Figure \ref{fig:zakesak_p2_mood2}, we clearly see that the first-order scheme will be also used in the region where the solution is smooth but PAD may not be satisfied. At the same time, we can see that the minimum values of $\bu_\sigma$ and $\xbar{\bu}_E$ are $>-10^{-5}$, which satisfy the assumed PAD criterion.   

\begin{figure}[ht!]
\centerline{\subfigure[point values $\bu_\sigma$]{\includegraphics[trim=1.2cm 2.2cm 2.5cm 1.8cm,clip,width=3.2cm]{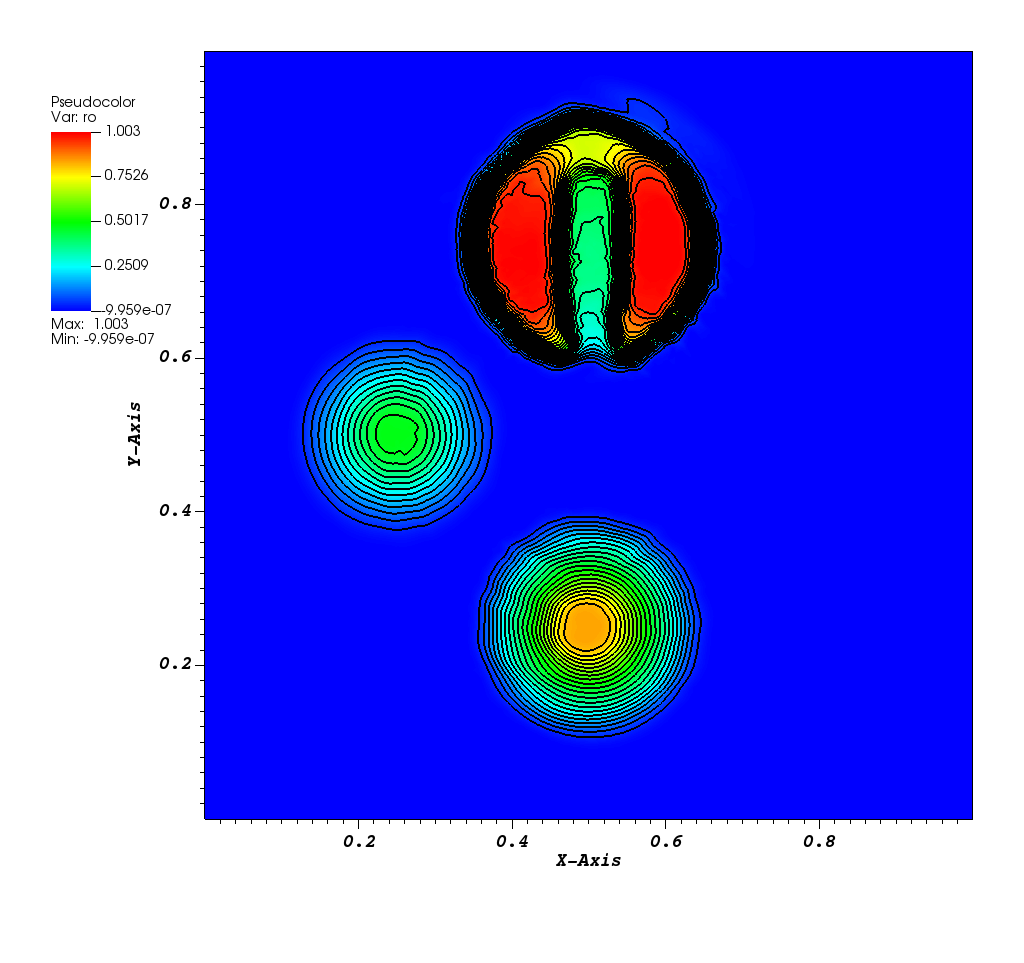}}\hspace*{0.02cm}
\subfigure[averages $\xbar{\bu}_E$]{\includegraphics[trim=1.2cm 2.2cm 2.5cm 1.8cm,clip,width=3.2cm]{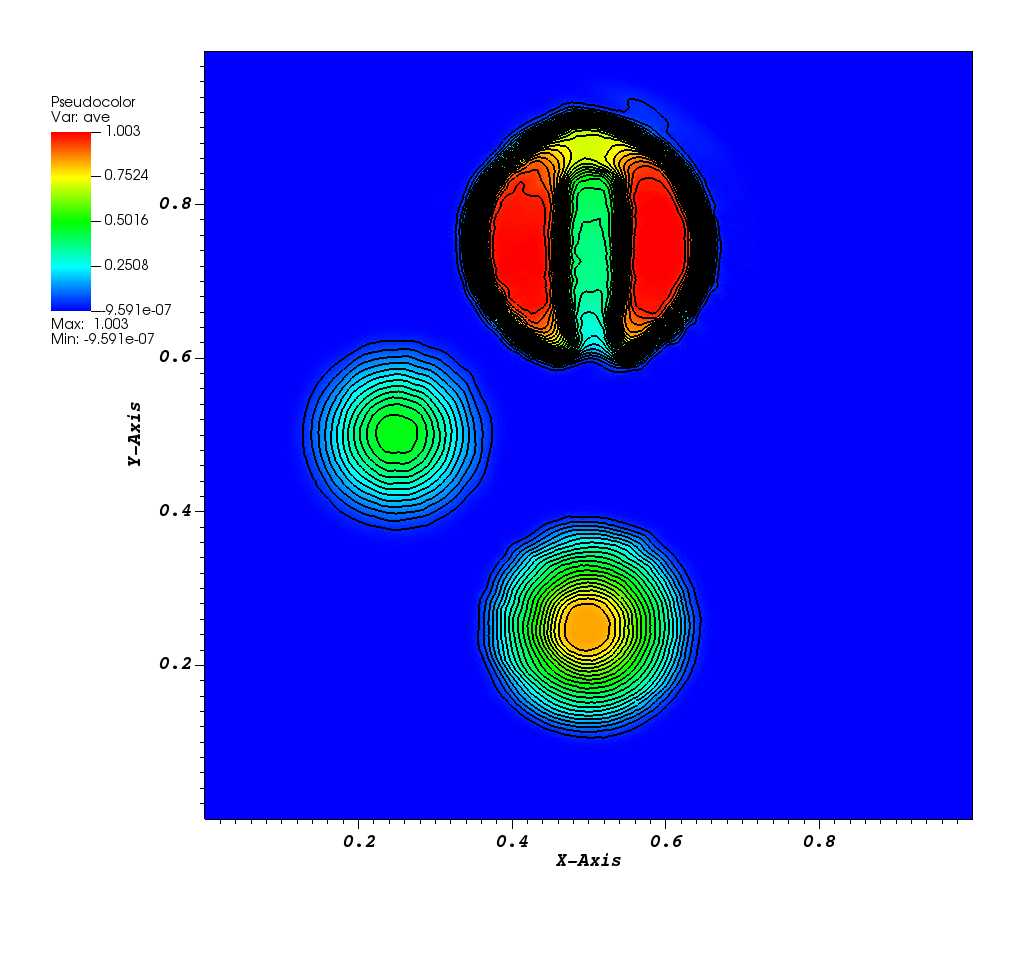}}\hspace*{0.02cm}
\subfigure[flag on $\bu_\sigma$]{\includegraphics[trim=1.2cm 2.2cm 2.5cm 1.8cm,clip,width=3.2cm]{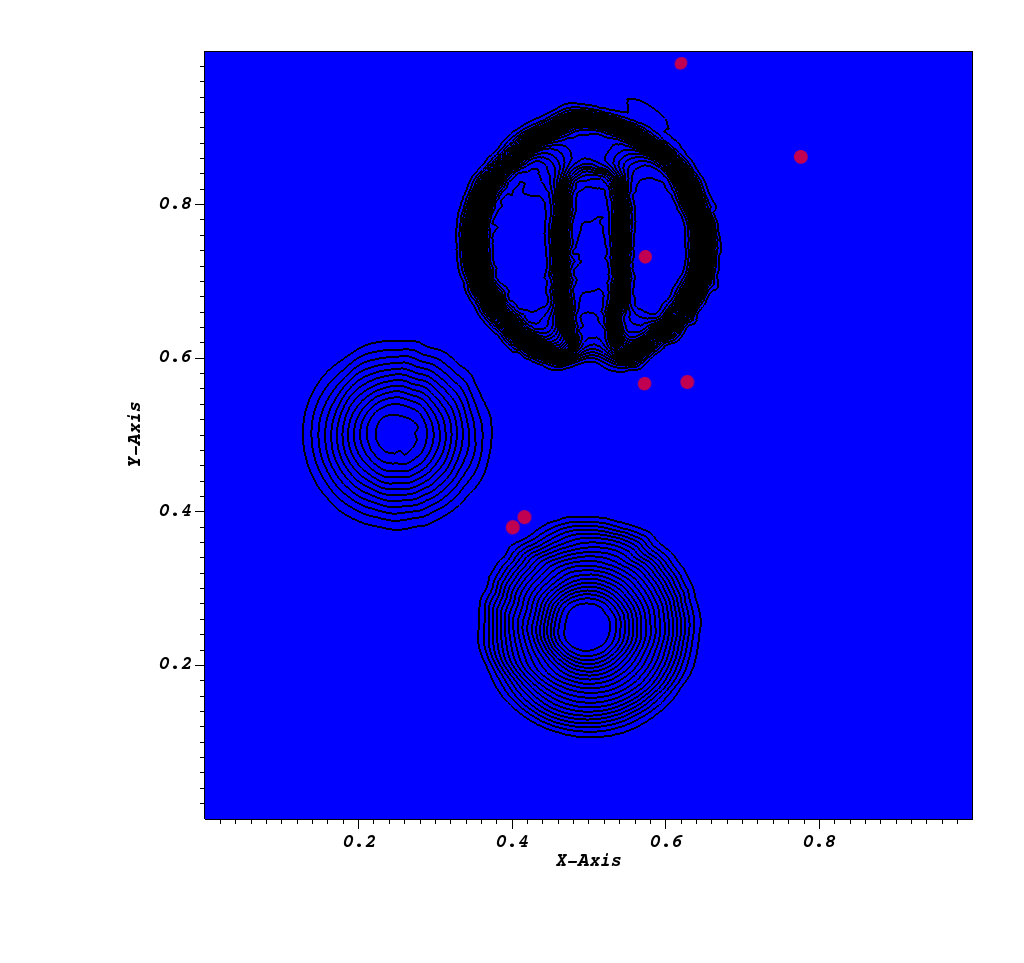}}\hspace*{0.02cm}
\subfigure[flag on $\xbar{\bu}_E$]{\includegraphics[trim=1.2cm 2.2cm 2.5cm 1.8cm,clip,width=3.2cm]{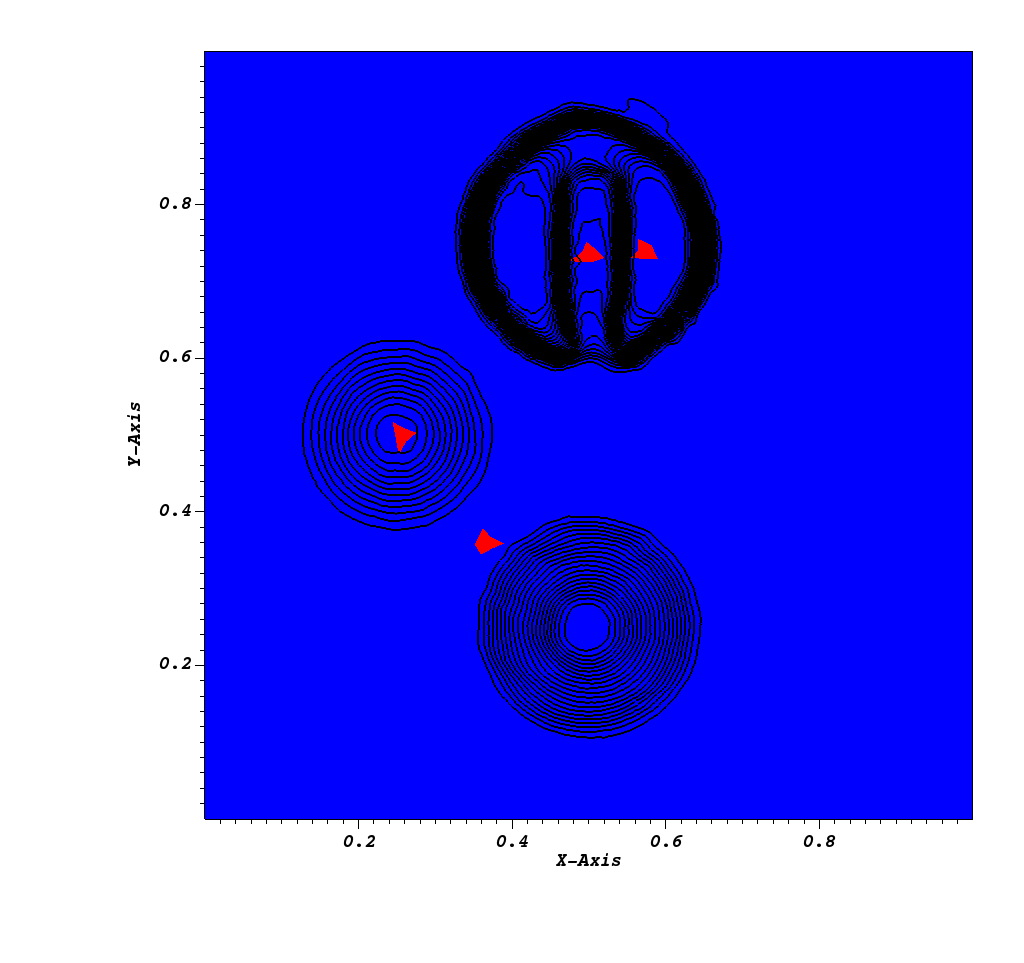}}}
\vskip5pt
\centerline{\subfigure[point values $\bu_\sigma$]{\includegraphics[trim=1.2cm 2.2cm 2.5cm 1.8cm,clip,width=3.2cm]{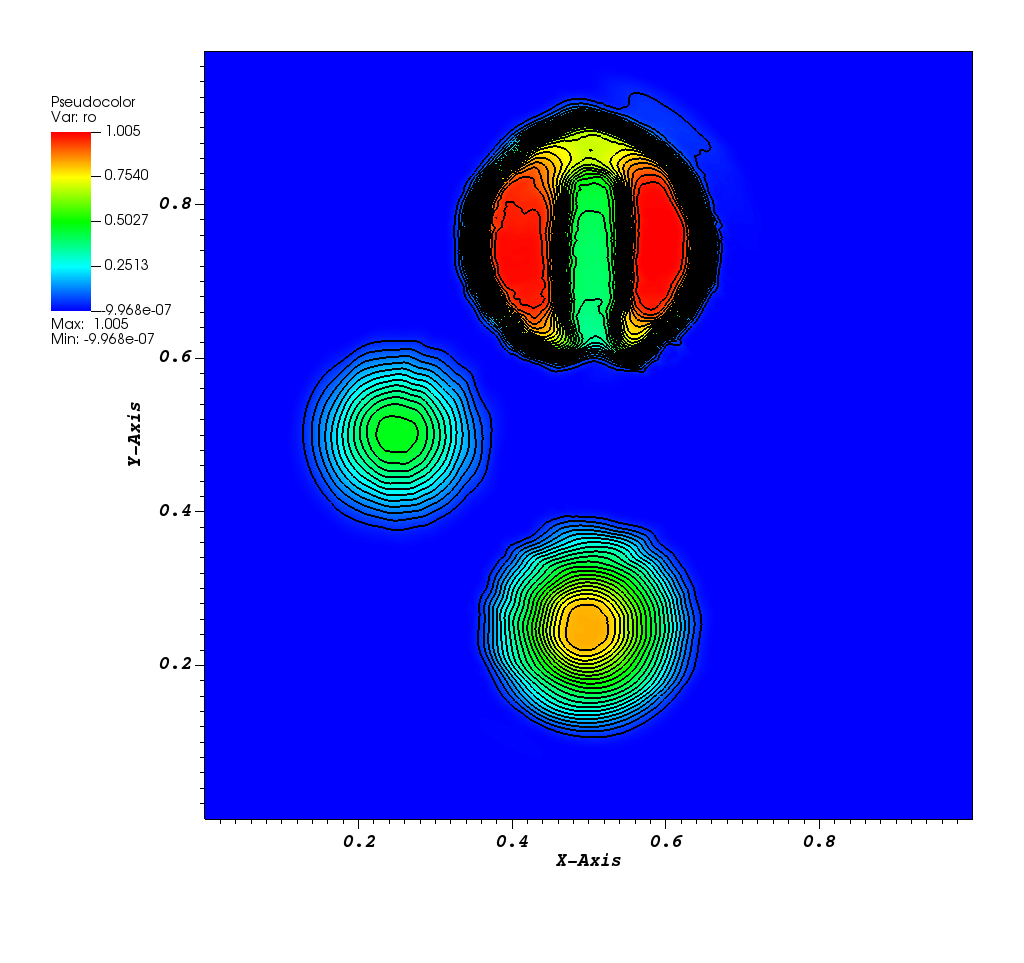}}\hspace*{0.02cm}
\subfigure[averages $\xbar{\bu}_E$]{\includegraphics[trim=1.2cm 2.2cm 2.5cm 1.8cm,clip,width=3.2cm]{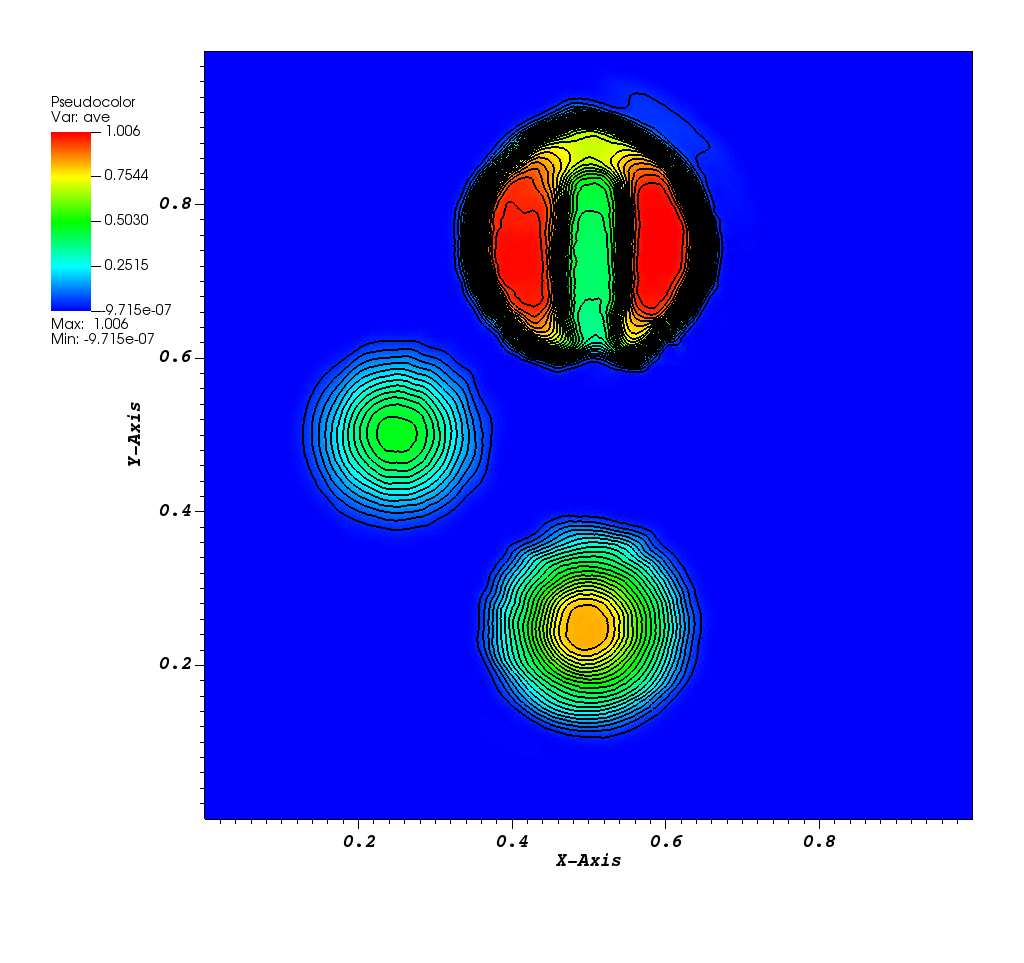}}\hspace*{0.02cm}
\subfigure[flag on $\bu_\sigma$]{\includegraphics[trim=1.2cm 2.2cm 2.5cm 1.8cm,clip,width=3.2cm]{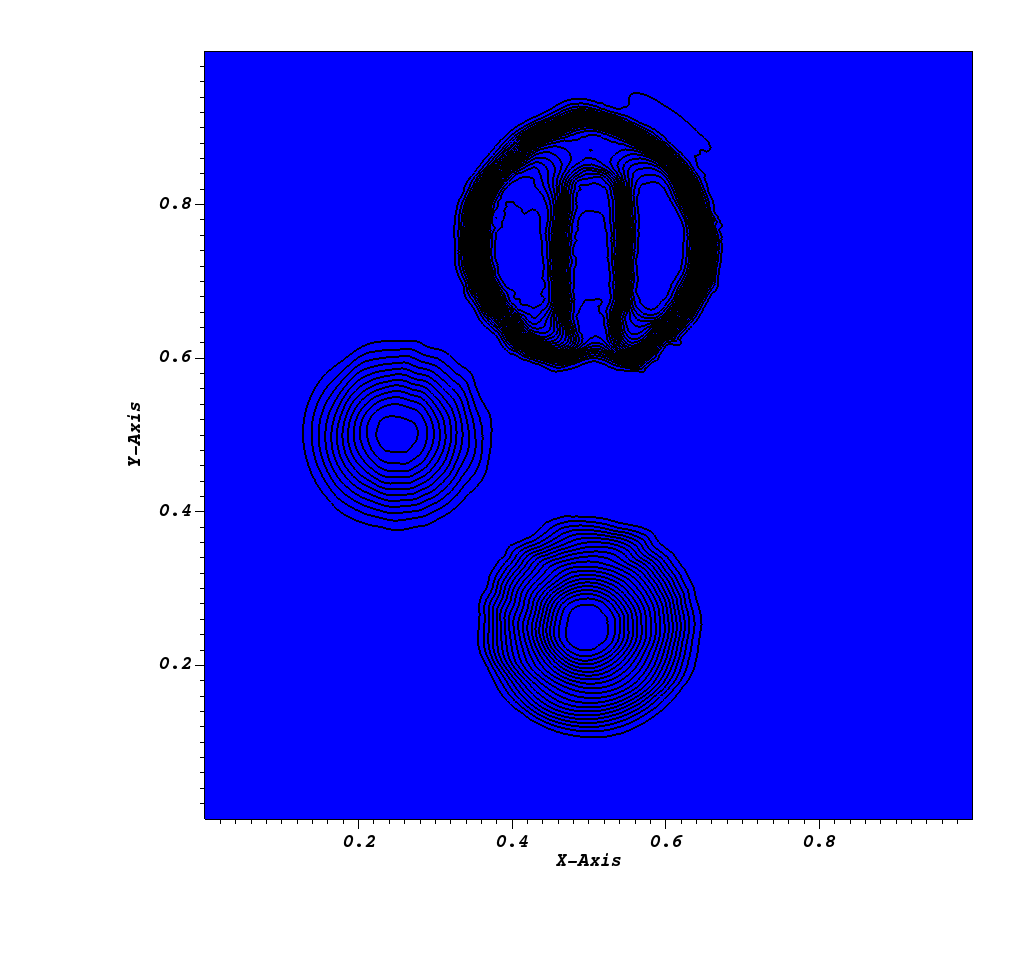}}\hspace*{0.02cm}
\subfigure[flag on $\xbar{\bu}_E$]{\includegraphics[trim=1.2cm 2.2cm 2.5cm 1.8cm,clip,width=3.2cm]{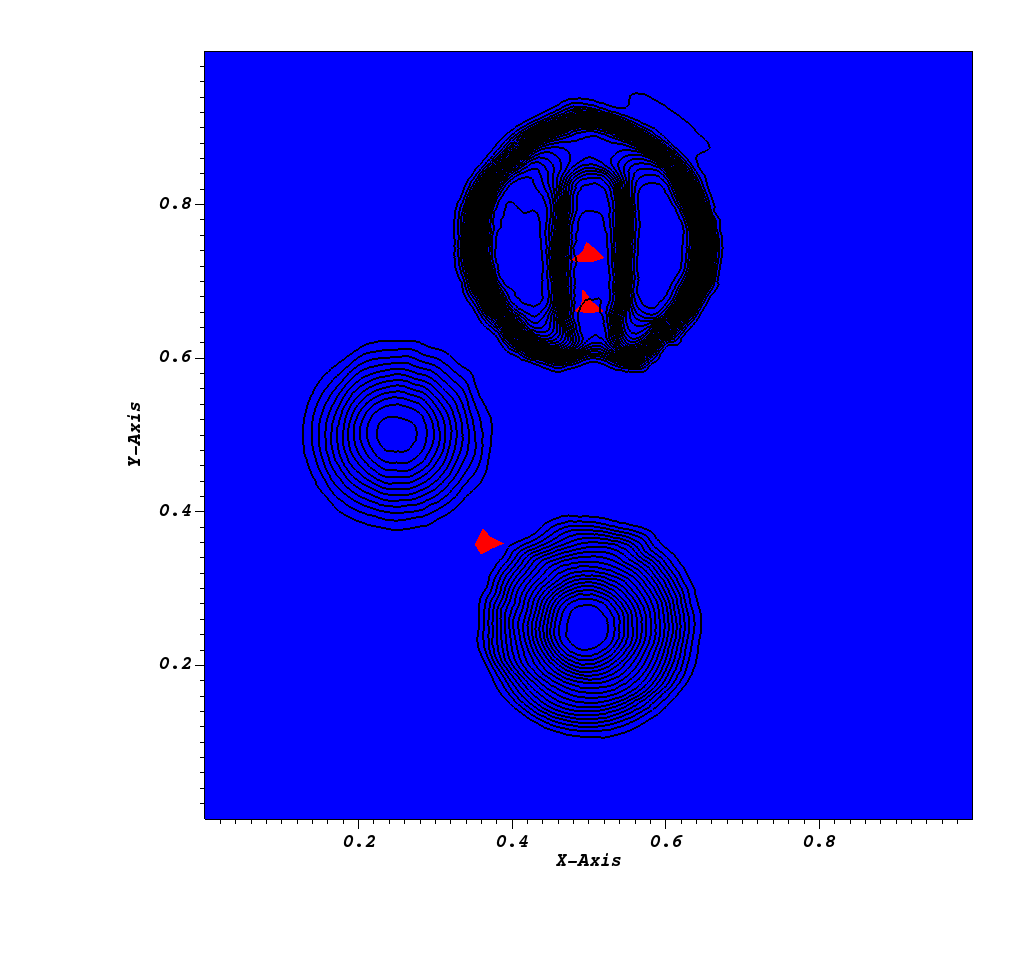}}}
\vskip5pt
\centerline{\subfigure[point values $\bu_\sigma$]{\includegraphics[trim=1.2cm 2.2cm 2.5cm 1.8cm,clip,width=3.2cm]{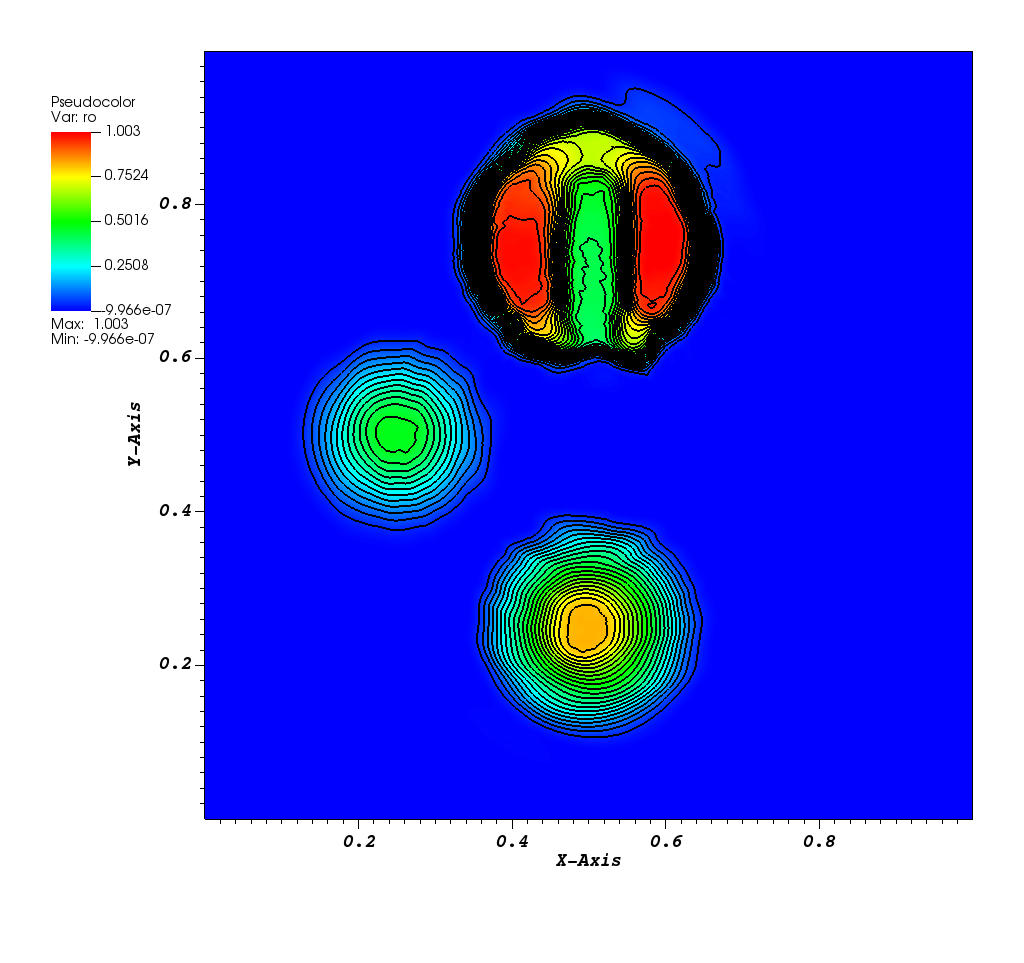}}\hspace*{0.02cm}
\subfigure[averages $\xbar{\bu}_E$]{\includegraphics[trim=1.2cm 2.2cm 2.5cm 1.8cm,clip,width=3.2cm]{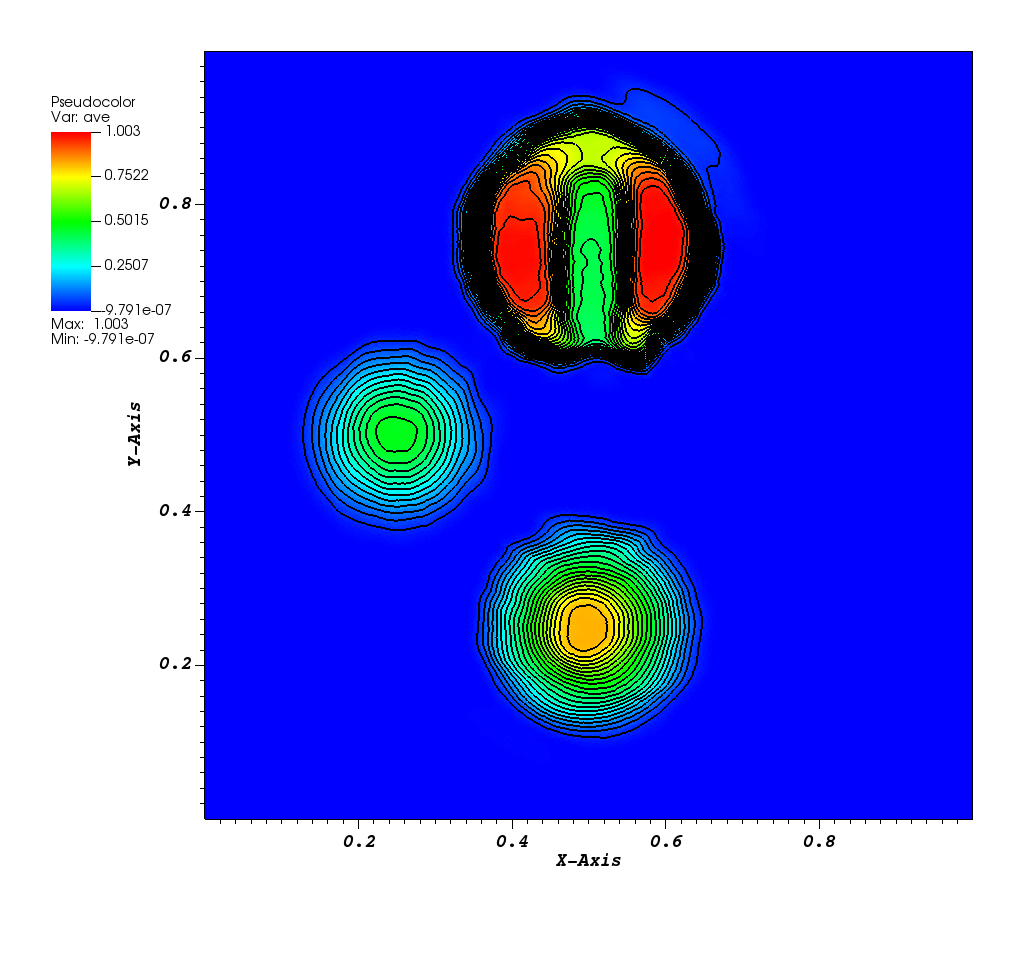}}\hspace*{0.02cm}
\subfigure[flag on $\bu_\sigma$]{\includegraphics[trim=1.2cm 2.2cm 2.5cm 1.8cm,clip,width=3.2cm]{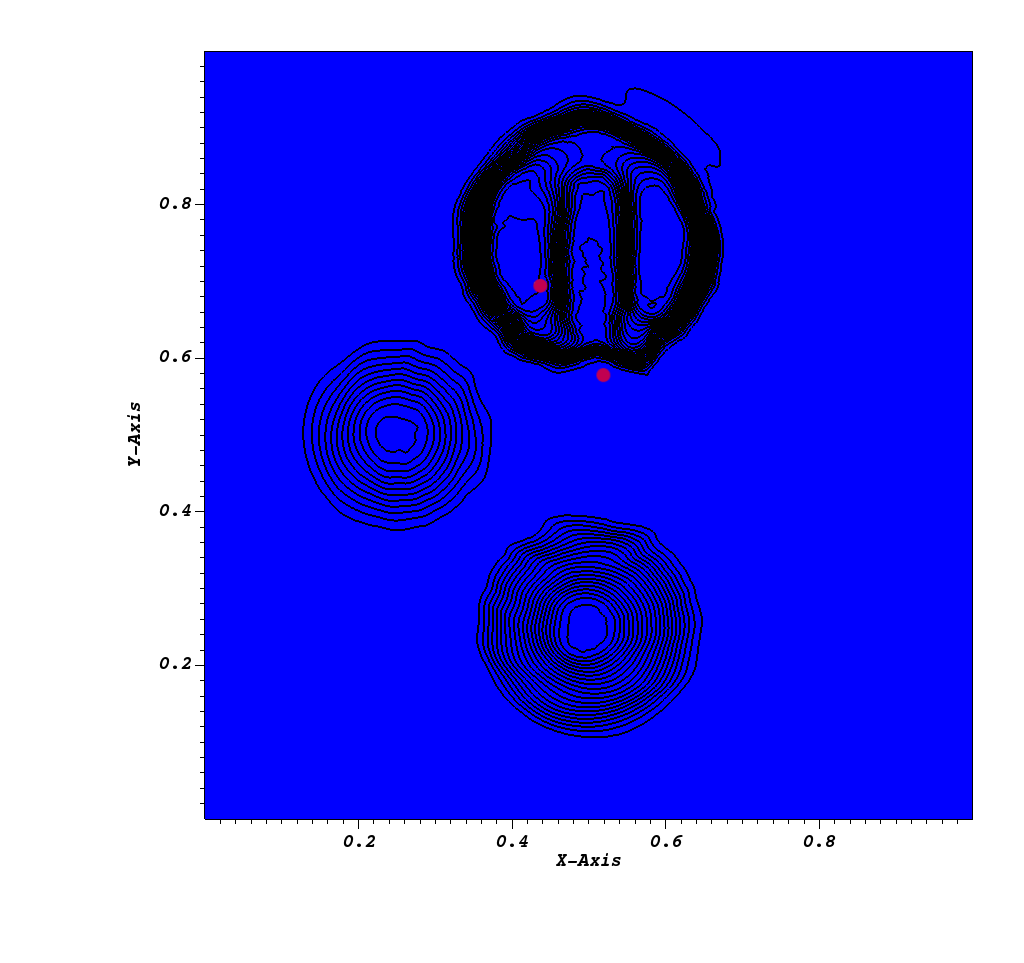}}\hspace*{0.02cm}
\subfigure[flag on $\xbar{\bu}_E$]{\includegraphics[trim=1.2cm 2.2cm 2.5cm 1.8cm,clip,width=3.2cm]{figures/Zalesak/zalesak_p3_flage2_t2.png}}}
\caption{Zalesak problem: cubic approximation and all MOOD criteria are checked. Top row: $T=1$; Middle row: $T=2$; Bottom row: $T=3$. \label{fig:zakesak_p3_mood2}}
\end{figure}

Finally, we compute the solutions using the cubic approximation. For saving space, only the results obtained with all MOOD criteria checked are presented in Figure \ref{fig:zakesak_p3_mood2}. Comparing the sensor results with these computed by the quadratic approximation, we observe that the cubic approximation leads to less activation of the MOOD paradigm.

\subsubsection{KPP test case}
\begin{figure}[ht!]
\centerline{\subfigure[point values $\bu_\sigma$]{\includegraphics[trim=1.2cm 2.2cm 2.5cm 2.5cm,clip,width=4.25cm]{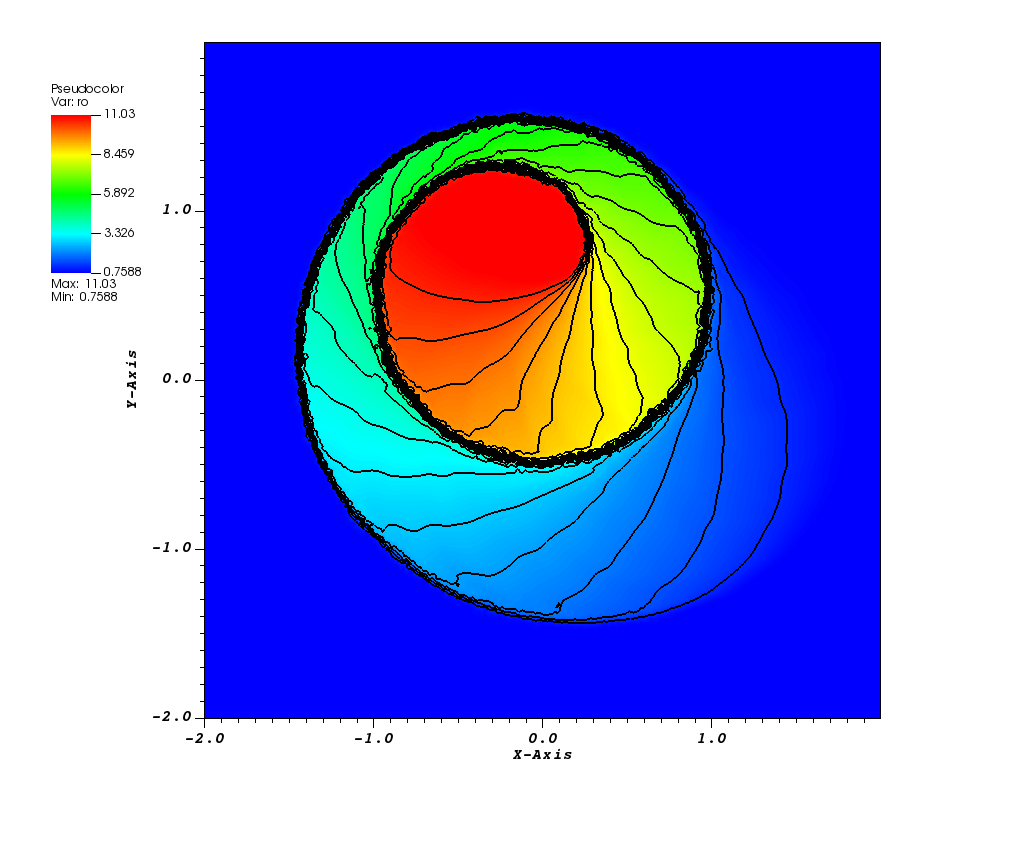}}\hspace*{0.05cm}
\subfigure[average values $\xbar{\bu}_E$]{\includegraphics[trim=1.2cm 2.2cm 2.5cm 2.5cm,clip,width=4.25cm]{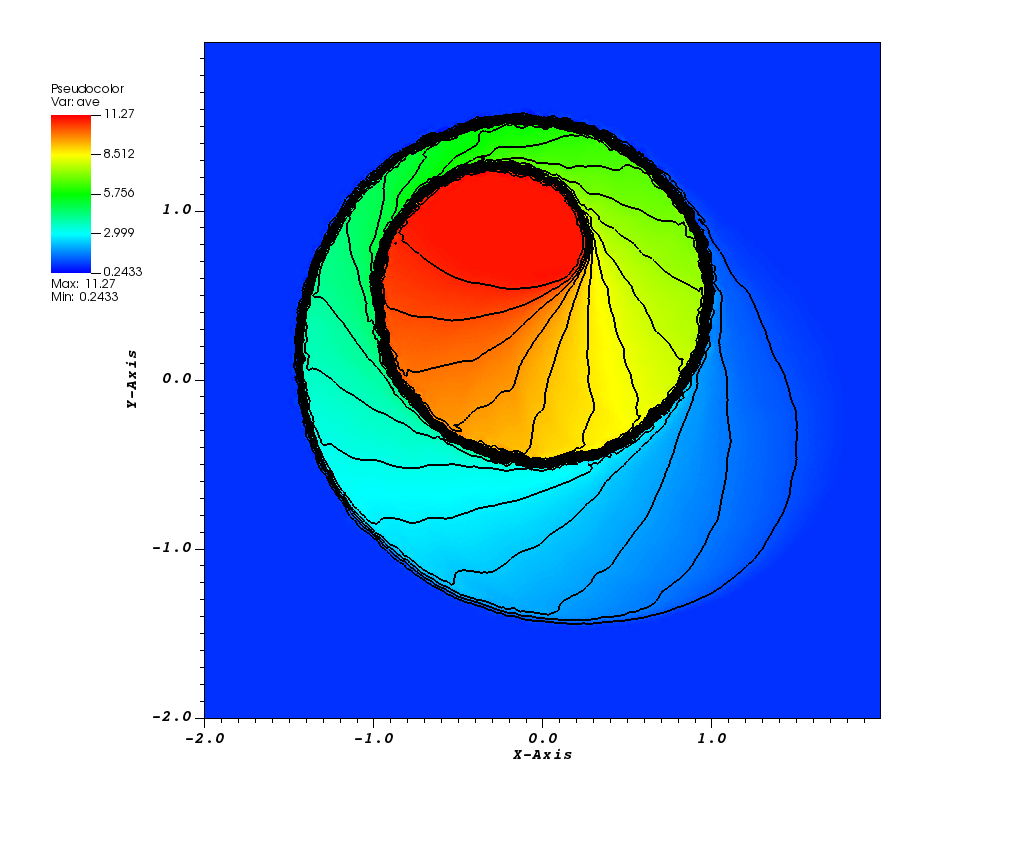}}\hspace*{0.05cm}
\subfigure[zoom of average values]{\includegraphics[trim=1.2cm 2.2cm 2.5cm 2.5cm,clip,width=4.25cm]{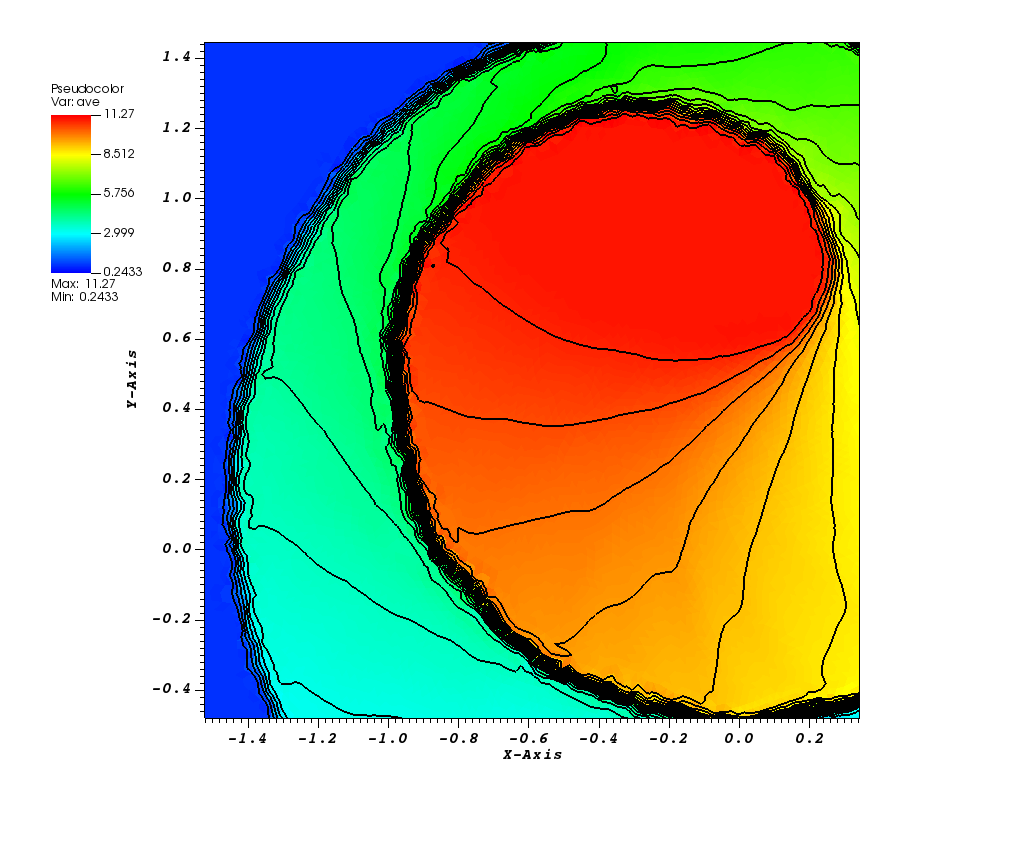}}}
\vskip5pt
\centerline{\subfigure[flag on point values $\bu_\sigma$]{\includegraphics[trim=1.2cm 2.2cm 2.5cm 2.5cm,clip,width=4.25cm]{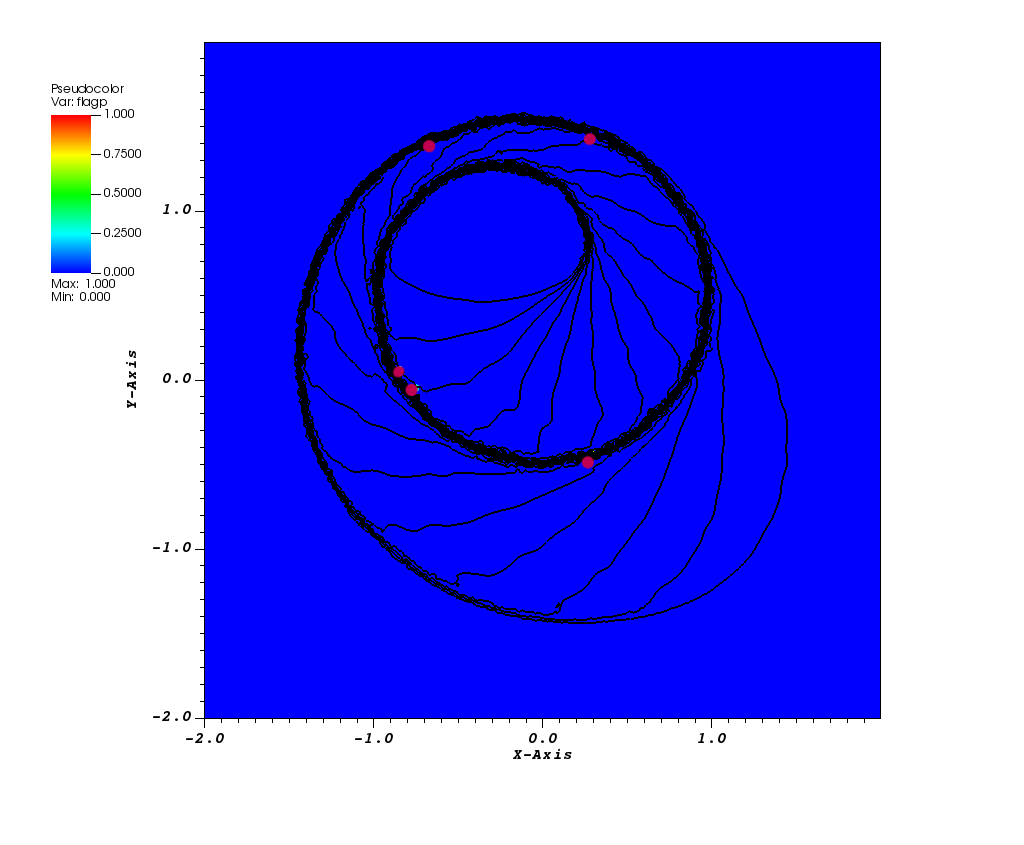}}\hspace*{0.15cm}
\subfigure[flag on average values $\xbar{\bu}_E$]{\includegraphics[trim=1.2cm 2.2cm 2.5cm 2.5cm,clip,width=4.25cm]{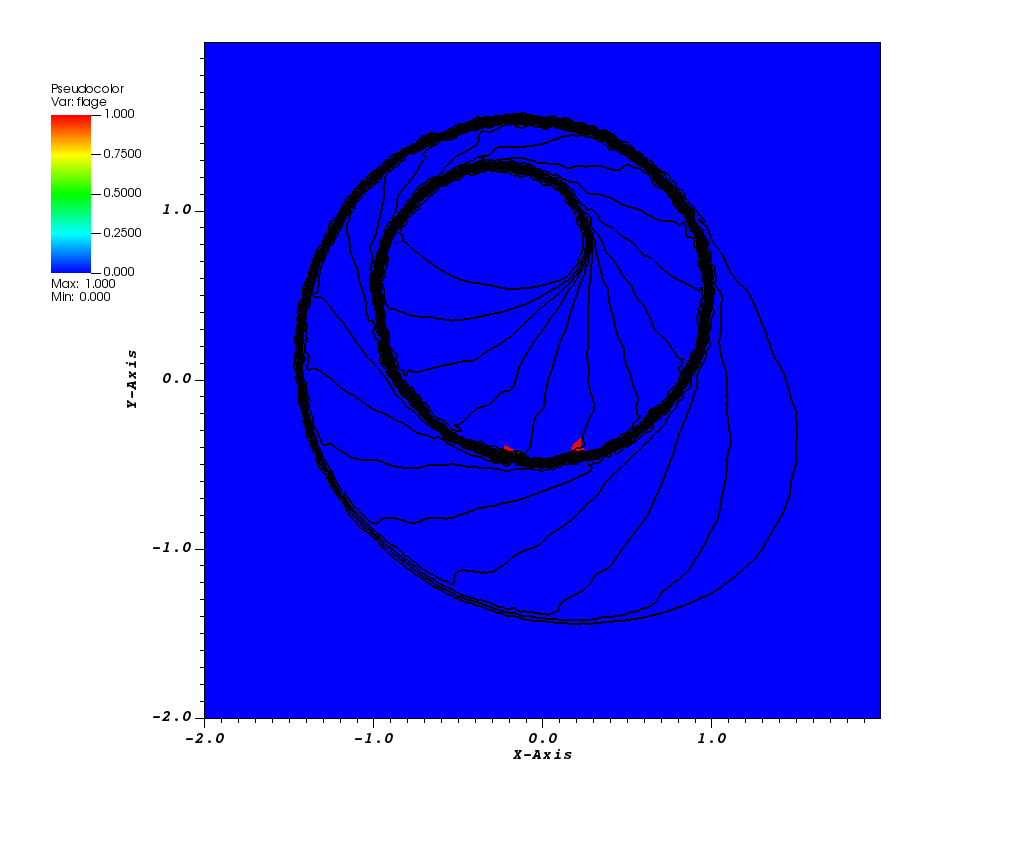}}}
\caption{KPP problem: using quadratic approximation.\label{fig:kpp_p2}}
\end{figure}

\begin{figure}[ht!]
\centerline{\subfigure[point values $\bu_\sigma$]{\includegraphics[trim=1.2cm 2.2cm 2.5cm 2.5cm,clip,width=4.25cm]{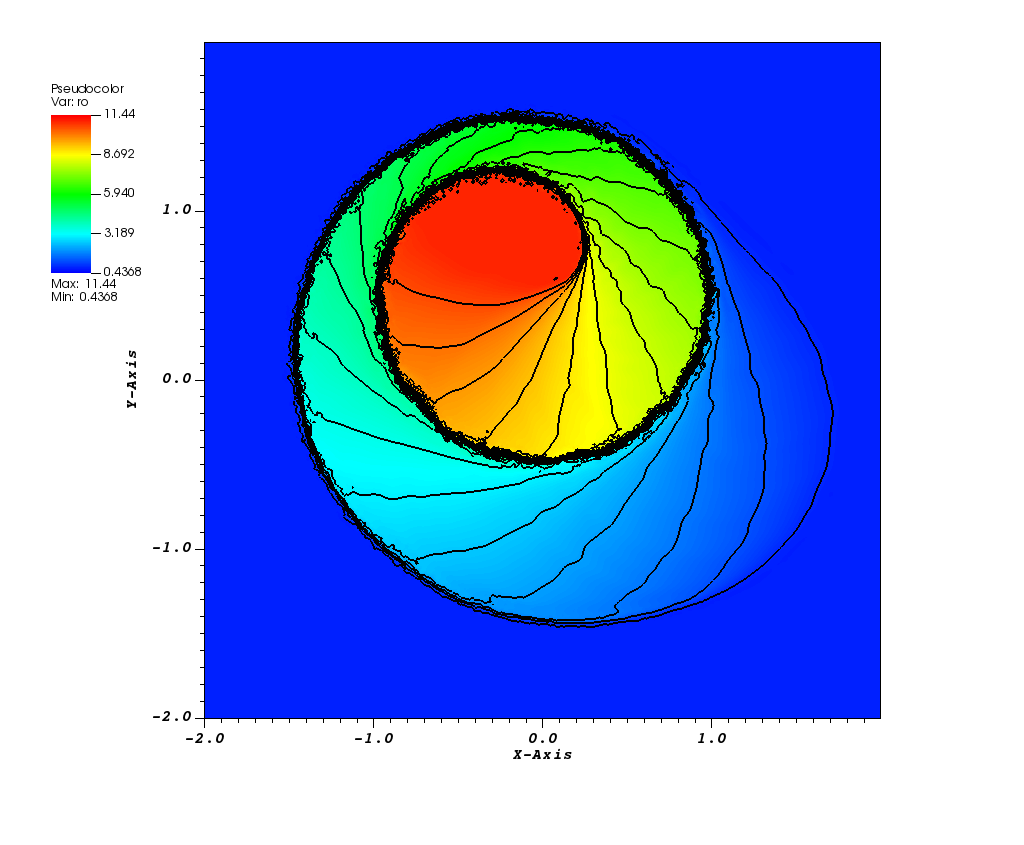}}\hspace*{0.05cm}
\subfigure[average values $\xbar{\bu}_E$]{\includegraphics[trim=1.2cm 2.2cm 2.5cm 2.5cm,clip,width=4.25cm]{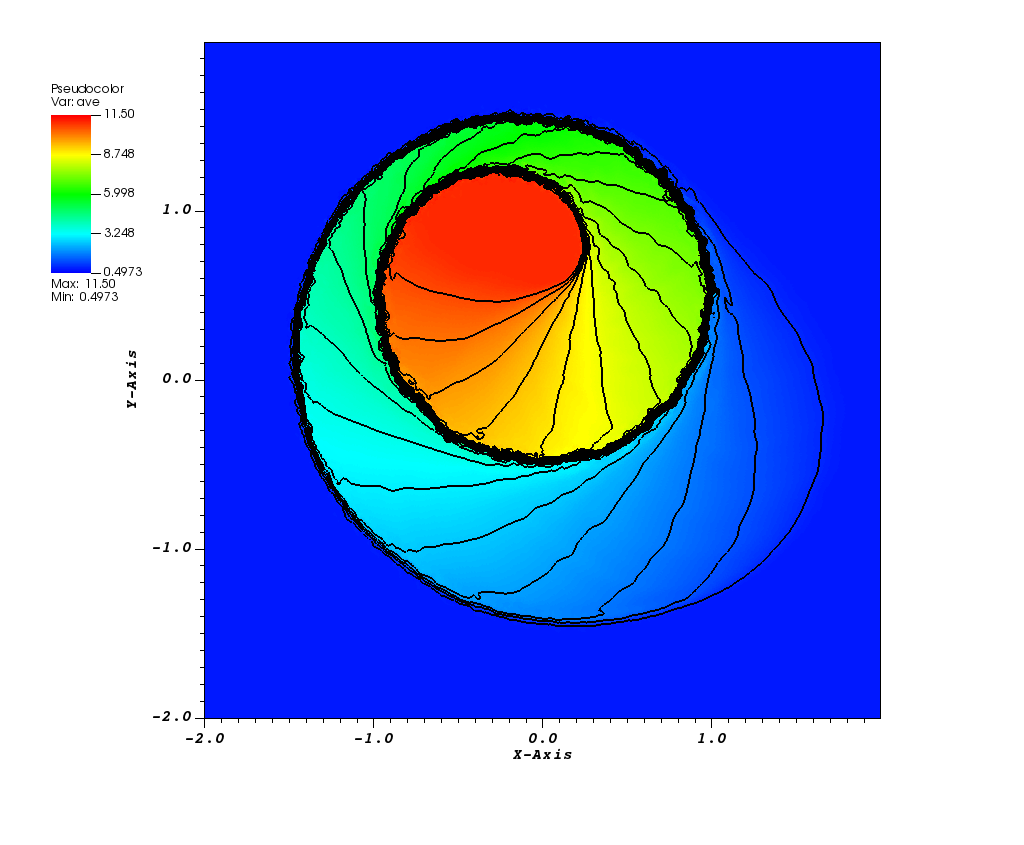}}\hspace*{0.05cm}
\subfigure[zoom of average values]{\includegraphics[trim=1.2cm 2.2cm 2.5cm 2.5cm,clip,width=4.25cm]{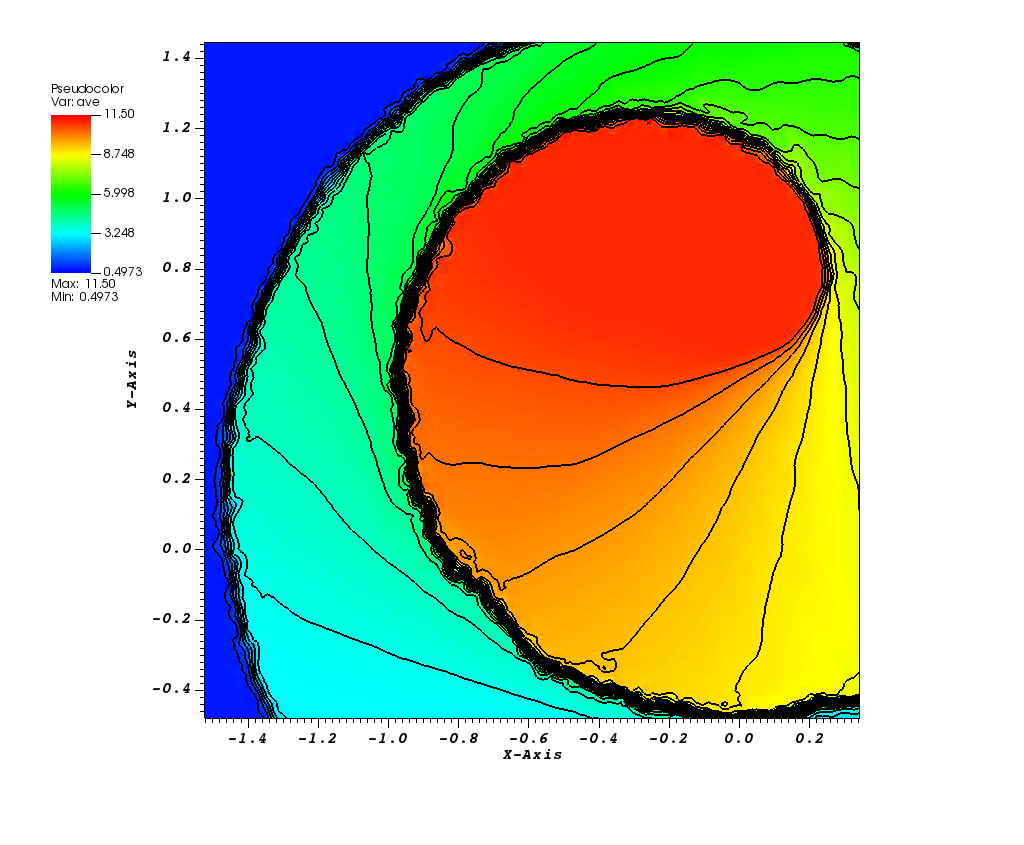}}}
\vskip5pt
\centerline{\subfigure[flag on point values $\bu_\sigma$]{\includegraphics[trim=1.2cm 2.2cm 2.5cm 2.5cm,clip,width=4.25cm]{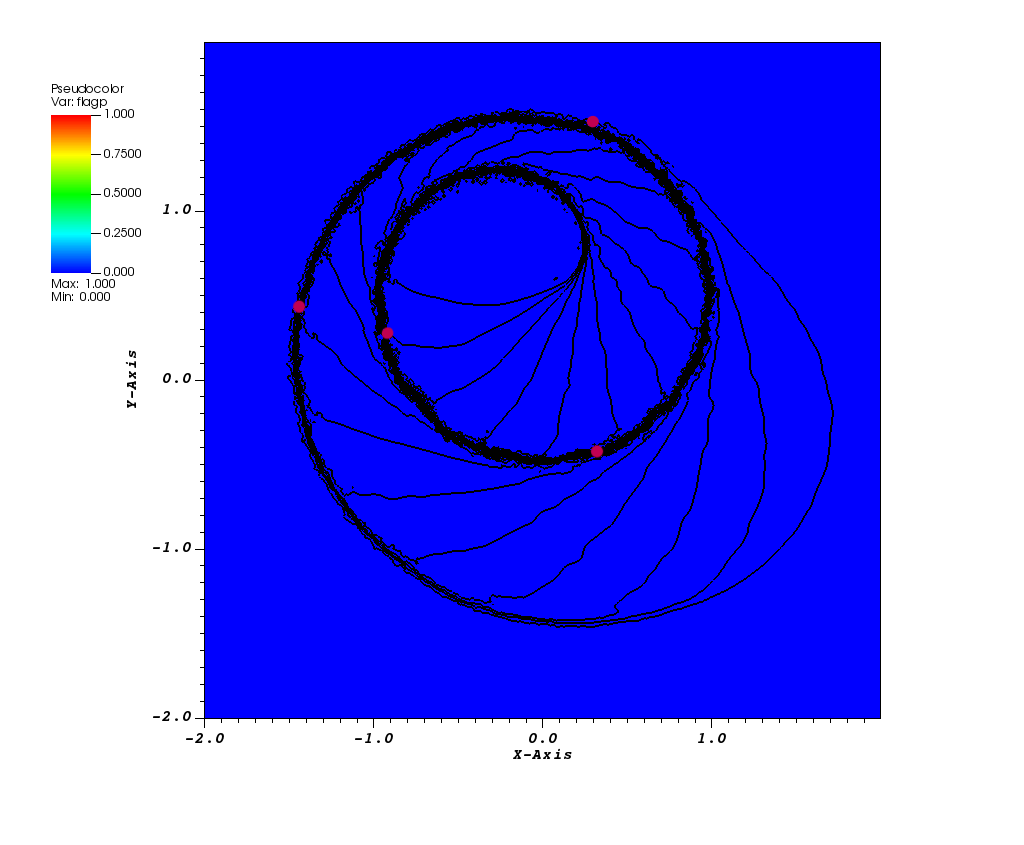}}\hspace*{0.15cm}
\subfigure[flag on average values $\xbar{\bu}_E$]{\includegraphics[trim=1.2cm 2.2cm 2.5cm 2.5cm,clip,width=4.25cm]{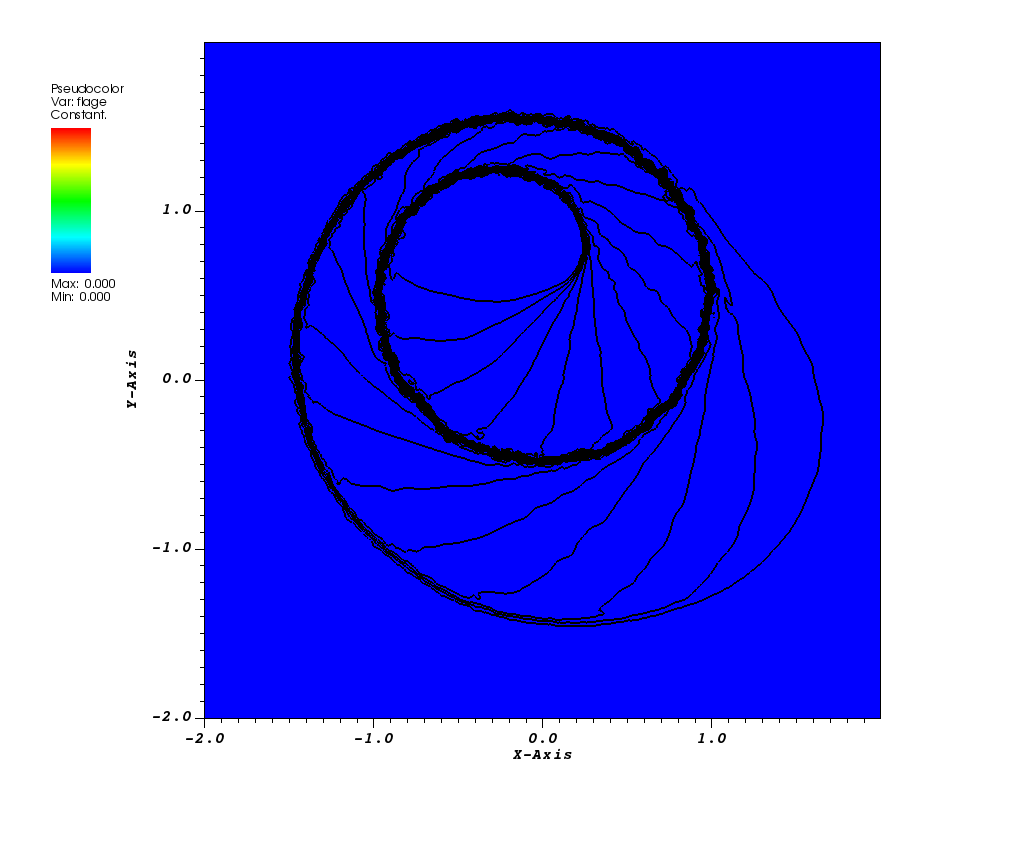}}}
\caption{KPP problem: using cubic approximation.\label{fig:kpp_p3}}
\end{figure}

In the third example, we consider the KPP problem which has been first considered by Kurganov, Popov and Petrova in \cite{zbMATH05380178}. The governing PDE is
\begin{equation}
\label{eq:KPP}
\dpar{u}{t}+\dpar{\big (\bla{\sin} \; u\big )}{x}+\dpar{\big (\bla{\cos} \; u\big )}{y}=0,
\end{equation}
in a domain $[-2,2]^2$
with the initial condition
$$u_0(\bx)=\left \{\begin{array}{ll}
\frac{7}{2}\pi& \text{ if } \Vert \bx-(0,0.5)\Vert^2 \leq 1,\\
\frac{\pi}{4}& \text{ else.}
\end{array}
\right .
$$
The problem \eqref{eq:KPP} is non-convex, in the sense that compound wave may exist. Here we have used the MOOD paradigm with the Local Lax-Friedrichs scheme as a first-order scheme. We compute the numerical solutions until the final time $T=1$ using both the quadratic and cubic approximations. For the quadratic (resp. cubic) approximation, the mesh has 8601 vertices, 25480 edges, 16880 elements, and 34081 (resp. 76441) degrees of freedom. We present the obtained numerical results in Figures \ref{fig:kpp_p2} and \ref{fig:kpp_p3}.

If the method is not dissipative enough, the two shocks in the northwest direction at around $22^\circ$ will be attached which is not correct. Here this is not the case, and the solution looks correct, at least in comparison with the existing literature on the subject. We are not aware of any form of exact analytical solution. We also indicate in red where the computed quantities ($\bu_\sigma$ or $\xbar{\bu}_E$) violate MOOD criteria (PAD is included) given in \S\ref{secMOOD}. As one can see, the low-order scheme is only activated at a limited number of locations. When the cubic approximation is used, the average values do not violate the MOOD criteria at the final time.

\subsection{Euler equations}\label{Euler_example}
We proceed with Euler equations of gas dynamics and consider the following six cases:
\begin{itemize}
\item a stationary vortex case in $[-10,10]^2$ for $T=3$.
\item a moving vortex case in $[-20,20]^2$ for $T=20$.
\item the Sod problem in $[-1,1]^2$ for $T=0.16$.
\item the 13th case studied by Liu and Lax \cite{zbMATH01929393}, where the solution is shown at time $T=1$.
\item the 3rd case of \cite{zbMATH01872721} by Kurganov and Tadmor for $T=1$.
\item the double Mach reflection case of \bla{Woodward and Collela, \cite{collela}}.
\end{itemize}

\subsubsection{Accuracy test}
In the first example of non-linear system cases, we verify the order of convergence of the proposed scheme. To this end, we consider a smooth stationary isentropic vortex flow. The computational domain is $[-10,10]\times[-10,10]$ with Dirichelet boundary conditions everywhere. The initial condition is given by 
\begin{equation}\label{Ex1_cond}
    \rho=(1+\delta T)^{\frac{1}{\gamma-1}}, \quad u=-\frac{\varepsilon}{2\pi}x_2e^{\frac{1-r^2}{2}}, \quad v=\frac{\varepsilon}{2\pi}x_1e^{\frac{1-r^2}{2}},\quad p=\rho^{\gamma},
\end{equation}
where the vortex strength is $\varepsilon=5$ and $\delta T=-\frac{(\gamma-1)\varepsilon^2}{8\gamma\pi^2}e^{1-r^2}$ with $r=\|\mathbf{x}\|=\sqrt{x_1^2+x_2^2}$. This is a stationary equilibrium of the system so the exact solution coincides with the initial condition at any time. The proposed scheme is not able to exactly preserve this equilibrium state and the truncation error is $\mathcal{O}(h^p)$, where $h$ is the characteristic size of the triangular mesh and $p$ is the order of the scheme. Thus, it severs as a good numerical example to verify the order of convergence, see also, e.g. in \cite{GBCKSD,ABT}.

Tables \ref{vortex:error:P2ave}--\ref{vortex:error:P3point} report the convergence rates of quadratic and cubic approximations for the vortex test
problem run on a sequence of successively refined meshes. The meshes are obtained from the Cartesian $40\times40$, $80\times80$, $160\times160$, $320\times320$, and $640\times640$ meshes which are all cut by the diagonal. From the obtained results, we can clearly see that the expected third- and fourth-accuracy are achieved by the proposed schemes with quadratic and cubic approximations, respectively.

\begin{table}[ht!]
\caption{\label{vortex:error:P2ave} Errors and rates on the average values of density, quadratic approximation, $T=3$.}
\begin{center}
\begin{tabular}{|c||c|c|c|c|c|c|}
\hline
$h$& $L^1$ &rate&  $L^2$& rate&$L^\infty$& rate \\
\hline
$3.6262\times{10^{-1}}$     &  $1.262\times10^{-4}$&- &  $2.850\times10^{-5}$& -     &     $8.406\times10^{-3}$  & - \\
$1.7901\times{10^{-1}} $ &  $2.220\times10^{-5}$&$2.46$&  $4.851\times10^{-6}$& $2.51$&       $1.364 \times10^{-3}$  & $2.58$\\
$8.8944\times{10^{-2}}$ &  $3.179\times10^{-6}$& $2.78$&  $6.865\times10^{-7}$& $2.80$&       $2.058\times10^{-4}$  & $2.70$\\
$4.4332\times{10^{-2}}$ &  $4.172\times10^{-7}$& $2.92$&  $9.093\times10^{-8}$& $2.90$&       $2.843\times10^{-5}$  & $2.84$\\
$2.2132\times{10^{-2}}$ &  $5.298\times10^{-8}$& $2.97$&  $1.160\times10^{-8}$& $2.96$&       $3.660\times10^{-6}$  & $2.95$\\
\hline
\end{tabular}
\end{center}
\end{table}

\begin{table}[ht!]
\caption{\label{vortex:error:P2point} Errors and rates on the point values of density, quadratic approximation, $T=3$.}
\begin{center}
\begin{tabular}{|c||c|c|c|c|c|c|}
\hline
$h$ & $L^1$ &rate&  $L^2$& rate&$L^\infty$& rate\\
\hline
$2.9608\times{10^{-1}}$   &  $9.100\times10^{-5}$&-&  $2.696\times10^{-5}$& -     &       $1.979\times10^{-2}$  & -   \\
$1.4616\times{10^{-1}} $&  $1.526\times10^{-5}$&$2.53$ &  $4.083\times10^{-6}$& $2.68$ &       $1.960 \times10^{-3}$  & $3.28$\\
$7.2623\times{10^{-2}}$ &  $2.182\times10^{-6}$& $2.78$&  $5.799\times10^{-7}$& $2.79$&       $2.355\times10^{-4}$  & $3.03$\\
$3.6198\times{10^{-2}}$ &  $2.860\times10^{-7}$& $2.92$&  $7.685\times10^{-8}$& $2.90$&       $3.323\times10^{-5}$  & $2.81$\\
$1.8070\times{10^{-2}}$&  $3.630\times10^{-8}$& $2.97$ &  $9.796\times10^{-9}$& $2.96$&       $4.356\times10^{-6}$  & $2.92$\\
\hline
\end{tabular}
\end{center}
\end{table}

\begin{table}[ht!]
\caption{\label{vortex:error:P3ave} Errors and rates on the average values of density, cubic approximation, $T=3$.}
\begin{center}
\begin{tabular}{|c||c|c|c|c|c|c|}
\hline
$h$& $L^1$& rate &  $L^2$ &rate&   $L^\infty$& rate\\
\hline
$3.6262\times{10^{-1}}$ & $1.956 \times10^{-5}$        & -     &  $4.651 \times10^{-6}$& -     & $2.051 \times10^{-3}$ &-\\

$1.7901\times{10^{-1}} $ &       $8.975 \times10^{-7}$  & $4.37$&  $1.603 \times10^{-7}$& $4.77$&  $5.134\times10^{-5}$ &$5.22$\\

$8.8944\times{10^{-2}}$ &       $4.970 \times10^{-8}$  & $4.14$&  $8.579 \times10^{-9}$& $4.19$&  $3.373 \times10^{-6}$& $3.89$\\

$4.4332\times{10^{-2}}$ &       $1.867 \times10^{-7}$  & $4.08$&  $4.829 \times10^{-10}$& $4.13$&  $2.895 \times10^{-9}$& $4.16$\\

$2.2132\times{10^{-2}}$ &  $1.755 \times10^{-10}$& $4.03$&  $2.863 \times10^{-11}$& $4.07$&       $1.128 \times10^{-8}$  & $4.04$\\
\hline
\end{tabular}
\end{center}
\end{table}

\begin{table}[ht!]
\caption{\label{vortex:error:P3point} Errors and rates on the point values of density, cubic approximation, $T=3$.}
\begin{center}
\begin{tabular}{|c||c|c|c|c|c|c|}
\hline
$h$& $L^1$ &rate&  $L^2$& rate &$L^\infty$ & rate \\
\hline
$2.4175\times{10^{-1}}$ &       $1.462 \times10^{-5}$  & -     &  $4.990  \times10^{-6}$& -     &  $2.169 \times10^{-3}$&-\\

$1.1934 \times{10^{-1}} $ &       $6.647\times10^{-7}$  & $4.38$&  $1.713 \times10^{-7}$& $4.78$&  $9.512\times10^{-5}$&$4.43$\\

$5.9296 \times{10^{-2}}$ &       $ 3.733 \times10^{-8}$  & $4.12$&  $9.757 \times10^{-9}$& $4.10$&  $5.200 \times10^{-6}$& $4.16$\\

$ 2.9555 \times{10^{-2}}$ &       $2.253 \times10^{-9}$  & $4.03$&   $5.953 \times10^{-10}$& $4.02$&  $3.018 \times10^{-7}$& $4.09$\\

$ 1.4754 \times{10^{-2}}$ &       $1.399 \times10^{-10}$  & $4.00$&   $3.798 \times10^{-11}$& $3.96$&  $ 2.254 \times10^{-8}$& $3.73$\\
\hline
\end{tabular}
\end{center}
\end{table}

\subsubsection{Moving vortex problem} 
\begin{figure}[ht!]
\centerline{\subfigure[point values $\rho_\sigma$, quadratic]{\includegraphics[trim=1.2cm 2.2cm 0.8cm 0.8cm,clip,width=6.0cm]{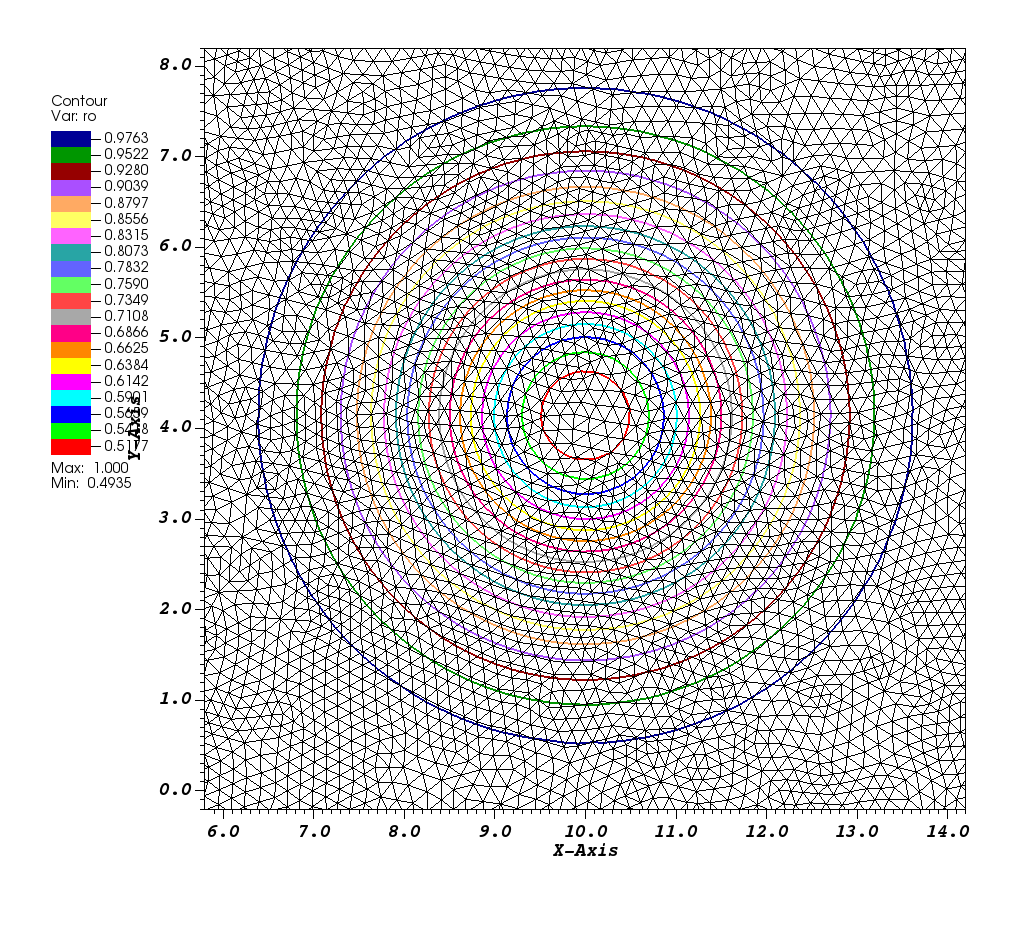}}\hspace*{0.05cm}
\subfigure[zoom of $\rho_\sigma$, quadratic]{\includegraphics[trim=1.2cm 2.2cm 0.8cm 0.8cm,clip,width=6.0cm]{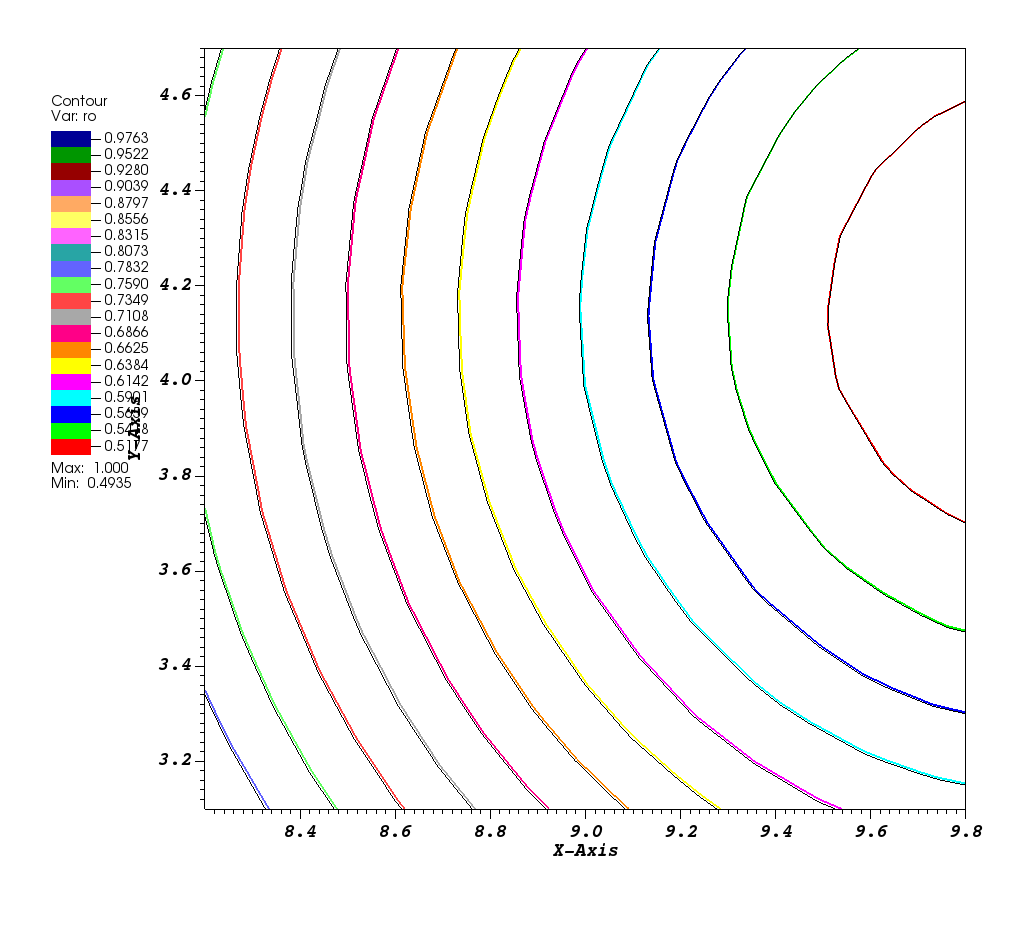}}}
\vskip5pt
\centerline{\subfigure[average values $\xbar{\rho}_E$, quadratic]{\includegraphics[trim=1.2cm 2.2cm 0.8cm 0.8cm,clip,width=6.0cm]{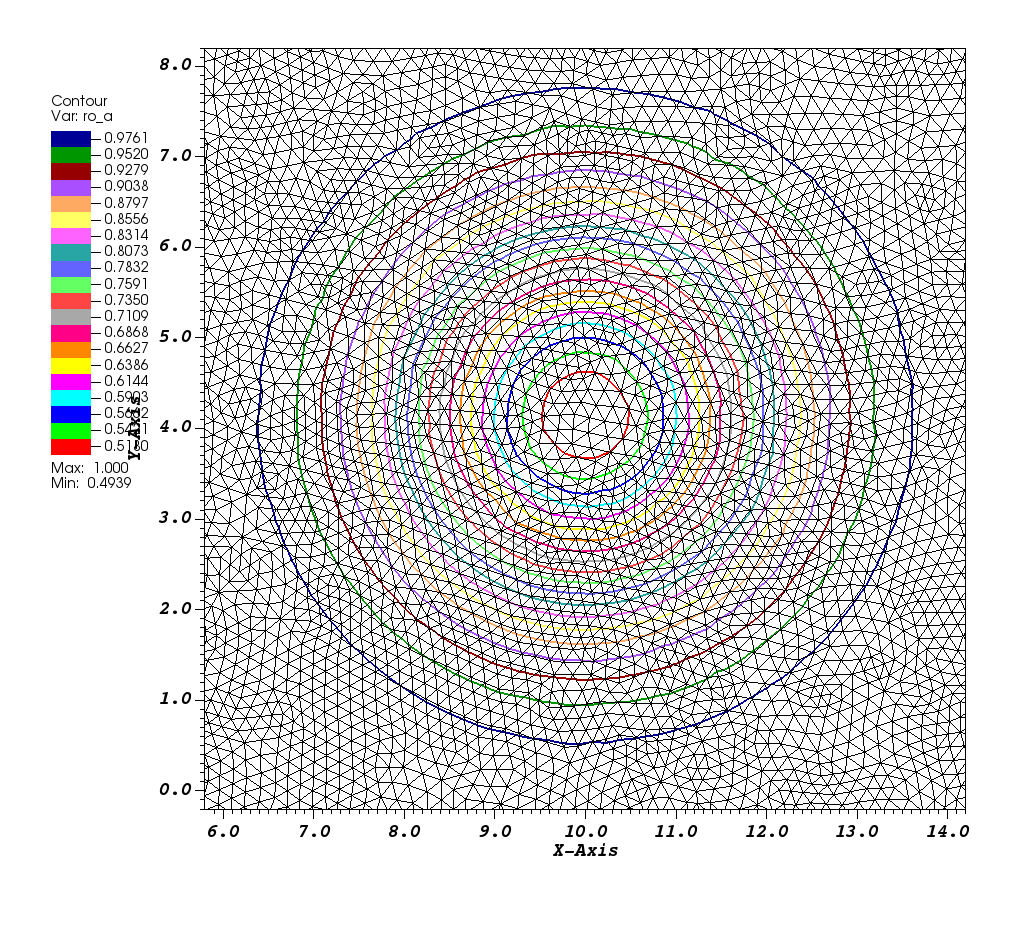}}\hspace*{0.05cm}
\subfigure[zoom of $\xbar{\rho}_E$, quadratic]{\includegraphics[trim=1.2cm 2.2cm 0.8cm 0.8cm,clip,width=6.0cm]{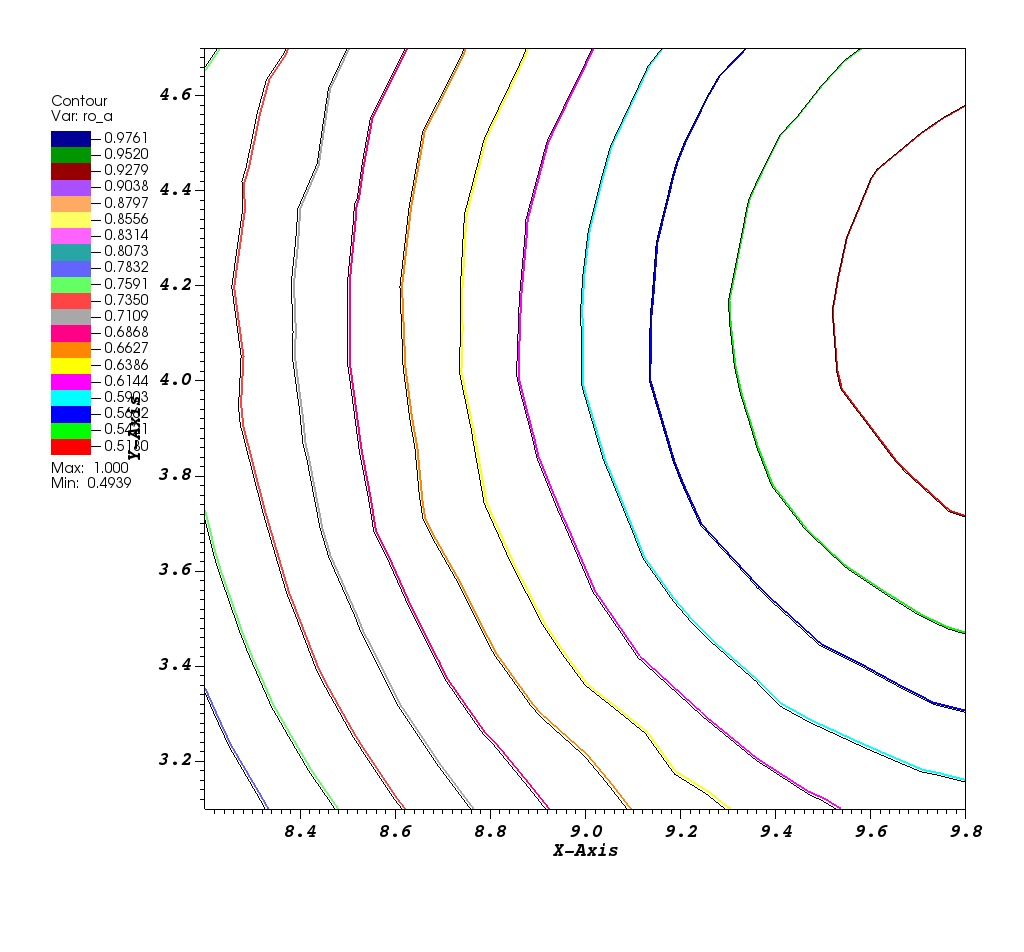}}}
\caption{\label{vortex:1} Moving vortex problem: density, zoom of the solution and the mesh. The exact solution is in black lines, the numerical solution with colored isolines. The mesh has 69705 DoFs, 35082 elements, third-order scheme.}
\end{figure}

\begin{figure}[ht!]
\centerline{\subfigure[point values $\rho_\sigma$, cubic]{\includegraphics[trim=1.2cm 2.2cm 0.8cm 0.8cm,clip,width=6.0cm]{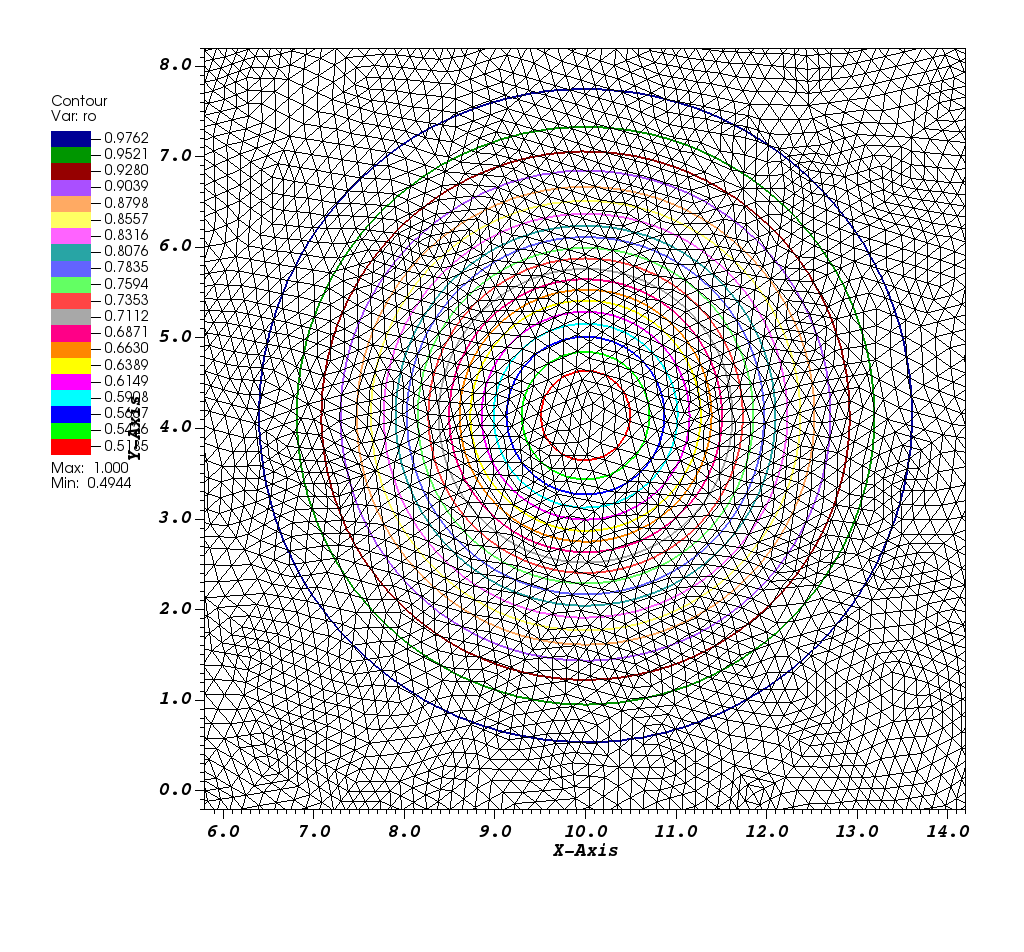}}\hspace*{0.05cm}
\subfigure[zoom of $\rho_\sigma$, cubic]{\includegraphics[trim=1.2cm 2.2cm 0.8cm 0.8cm,clip,width=6.0cm]{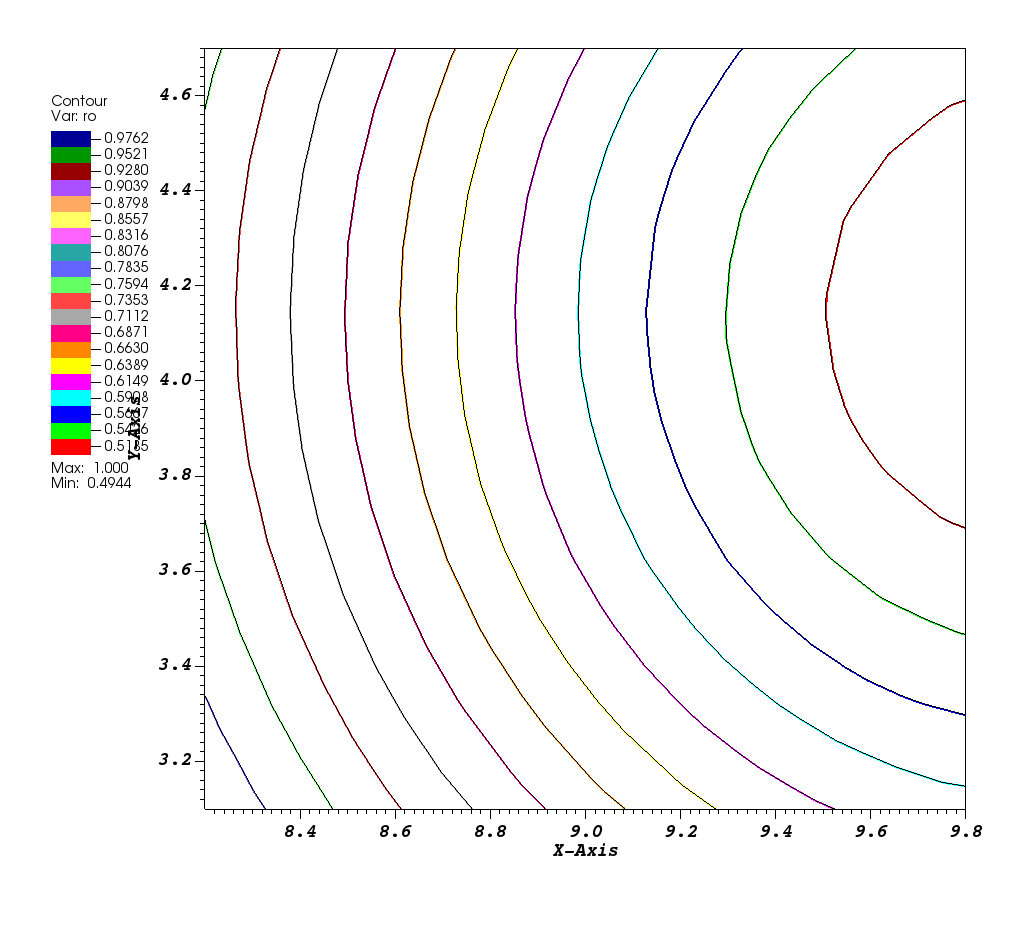}}}
\vskip5pt
\centerline{\subfigure[average values $\xbar{\rho}_E$, cubic]{\includegraphics[trim=1.2cm 2.2cm 0.8cm 0.8cm,clip,width=6.0cm]{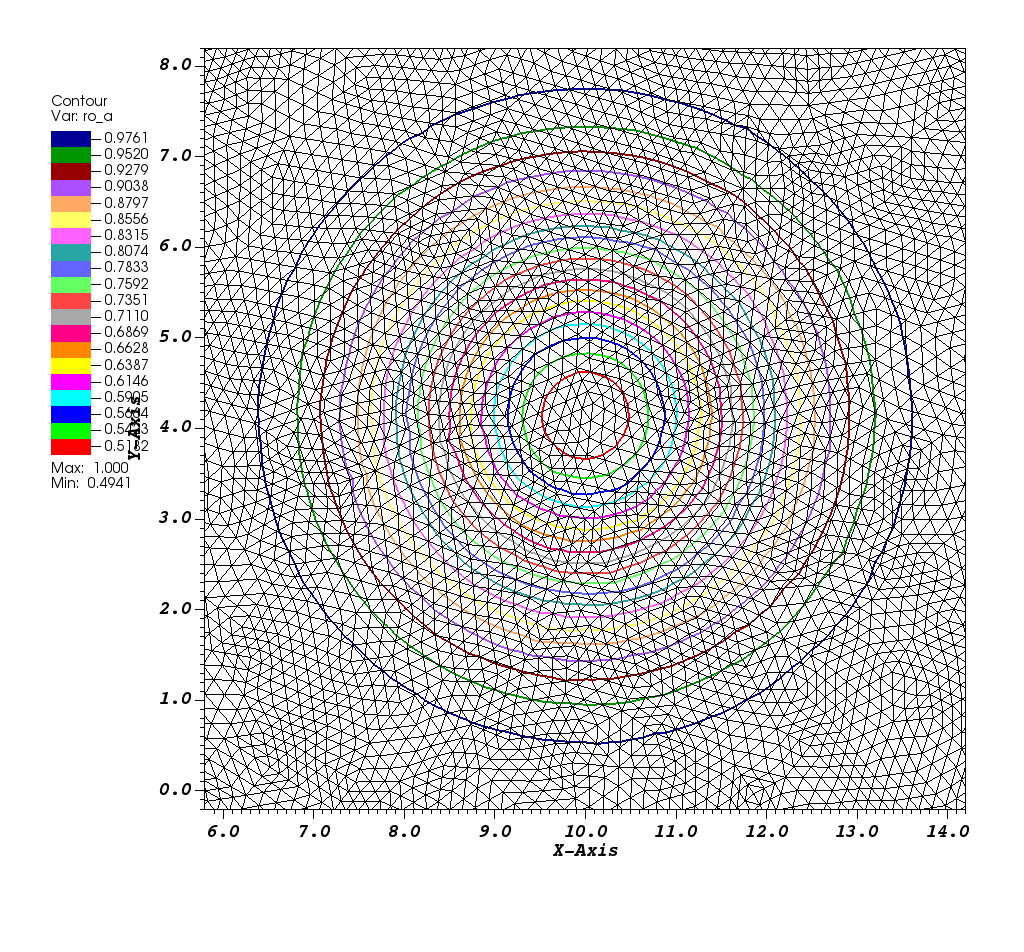}}\hspace*{0.05cm}
\subfigure[zoom of $\xbar{\rho}_E$, cubic]{\includegraphics[trim=1.2cm 2.2cm 0.8cm 0.8cm,clip,width=6.0cm]{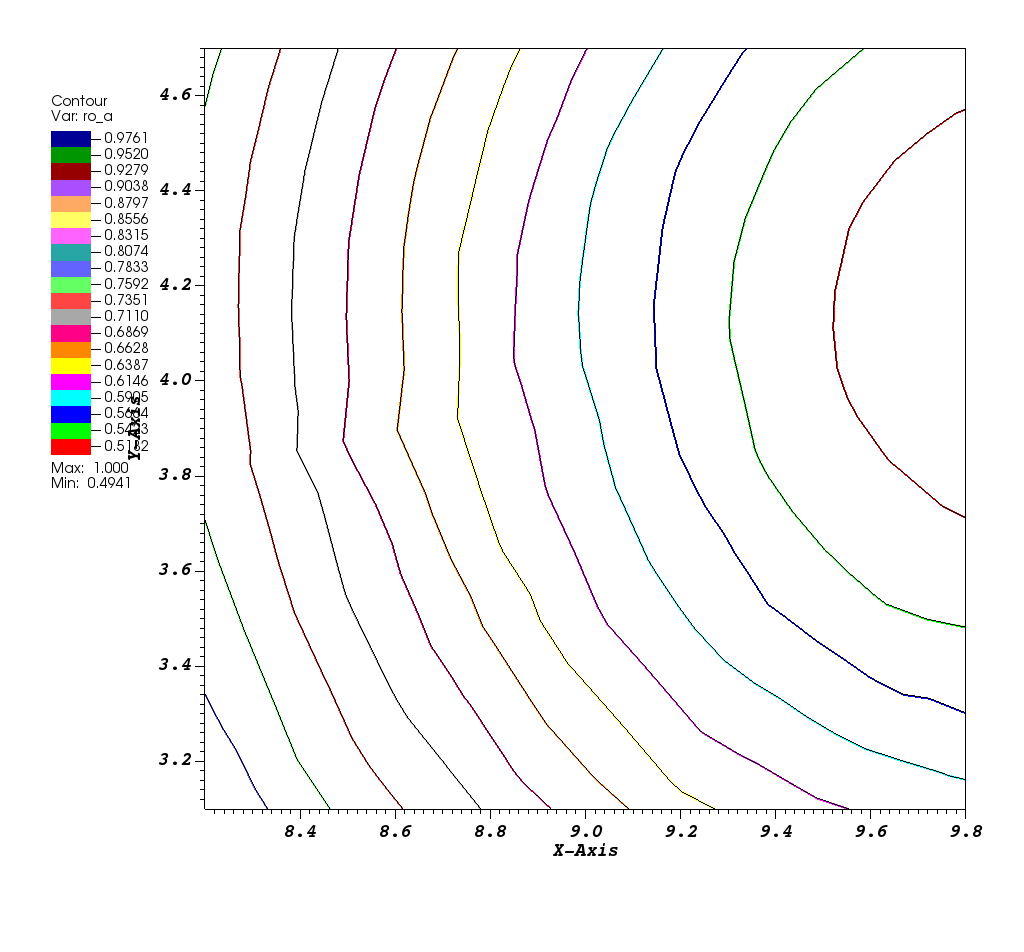}}}
\caption{\label{vortex:1_3} Moving vortex problem: density, zoom of the solution and the mesh. The exact solution is in black lines, the numerical solution with colored isolines. The mesh has 70762 DoFs, 15930 elements, fourth-order scheme.}
\end{figure}

In the second example of non-linear system cases, we consider the following initial condition:
\begin{equation}\label{ex2cond}
\rho=\big ( T_\infty+\Delta T\big )^{\frac{1}{\gamma-1}},~ u=u_\infty-{y_2}Me^{\frac{1-R}{2}},~ v=v_\infty+y_1Me^{\frac{1-R}{2}},~ p=\rho^\gamma,
\end{equation}
prescribed in the domain $[-20,20]^2$. In \eqref{ex2cond}, $(u_\infty,v_\infty)=(1,\frac{\sqrt{2}}{2})$, $T_\infty=1$, $M=\frac{5}{2\pi}$,  $\by=({\bf x-x_{0}})/2$ with $\bx_{0}=(-10,-10)$, $R=\Vert\by\Vert^2$, and $\Delta T=-\frac{\gamma-1}{2\gamma}M^2e^{{1-R}}$. 
This vortex is simply advected with the velocity $(u_\infty, v_{\infty})$ so that it is very easy to compute the exact solution. The third- and fourth-order results at the final time $T=20$, for the scheme \eqref{scheme:fv}-\eqref{scheme:pt}, are shown in Figures \ref{vortex:1} and \ref{vortex:1_3}. The mesh as a comparison with the exact solution is also shown. 
We have also displayed zooms of the solutions at $T=20$, superimposed with the exact one. For the fourth-order case, it is difficult to see a difference, at this scale, between the exact and numerical solutions. In the two cases, the mesh that is plotted is obtained from the actual computational mesh where we have sub-divided the elements using the Lagrange points: these meshes represent all the degrees of freedom that we actually use. It turns out that the plotting meshes (represented) are identical for the two types of solutions.

\subsubsection{2D Sod problem}
\begin{figure}[ht!]
\centerline{\subfigure[point values $\rho_\sigma$, quadratic]{\includegraphics[trim=1.2cm 2.2cm 1.2cm 1.2cm,clip,width=6.05cm]{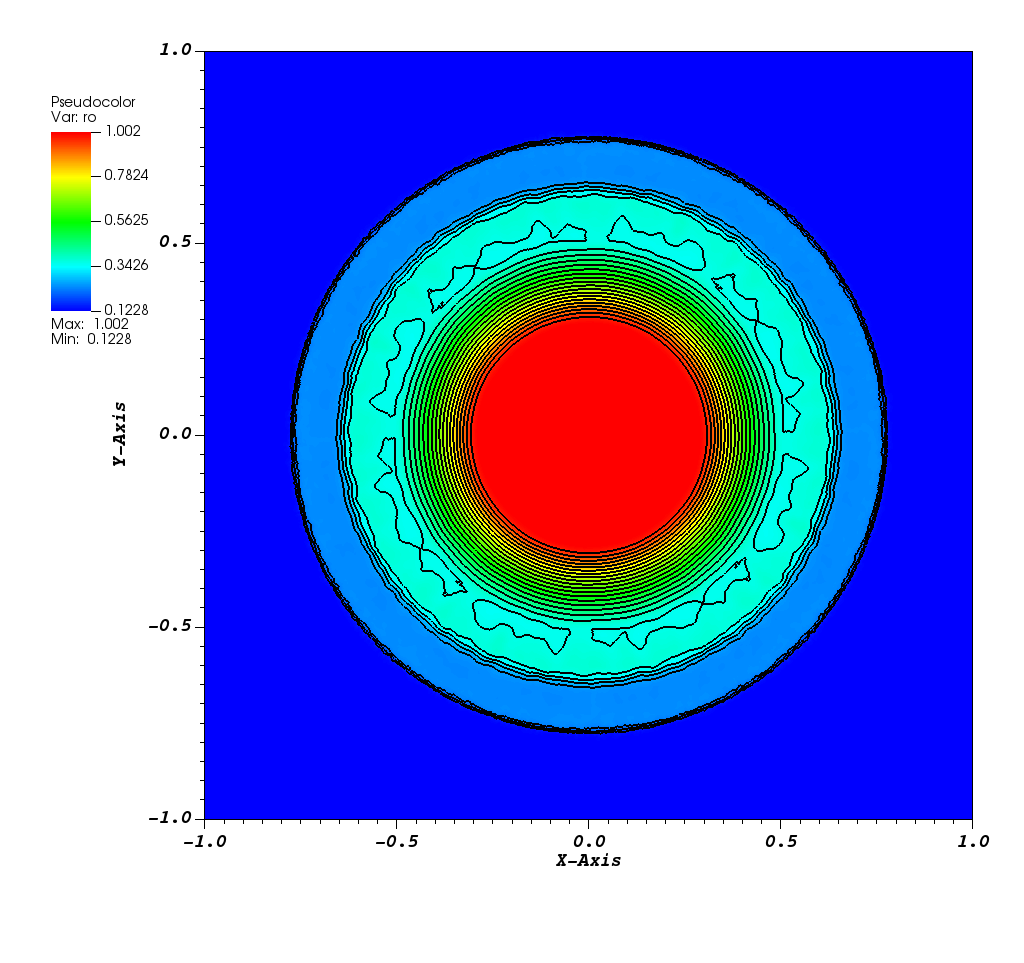}}\hspace*{0.05cm}
\subfigure[average values $\xbar{\rho}_E$, quadratic]{\includegraphics[trim=1.2cm 2.2cm 1.2cm 1.2cm,clip,width=6.05cm]{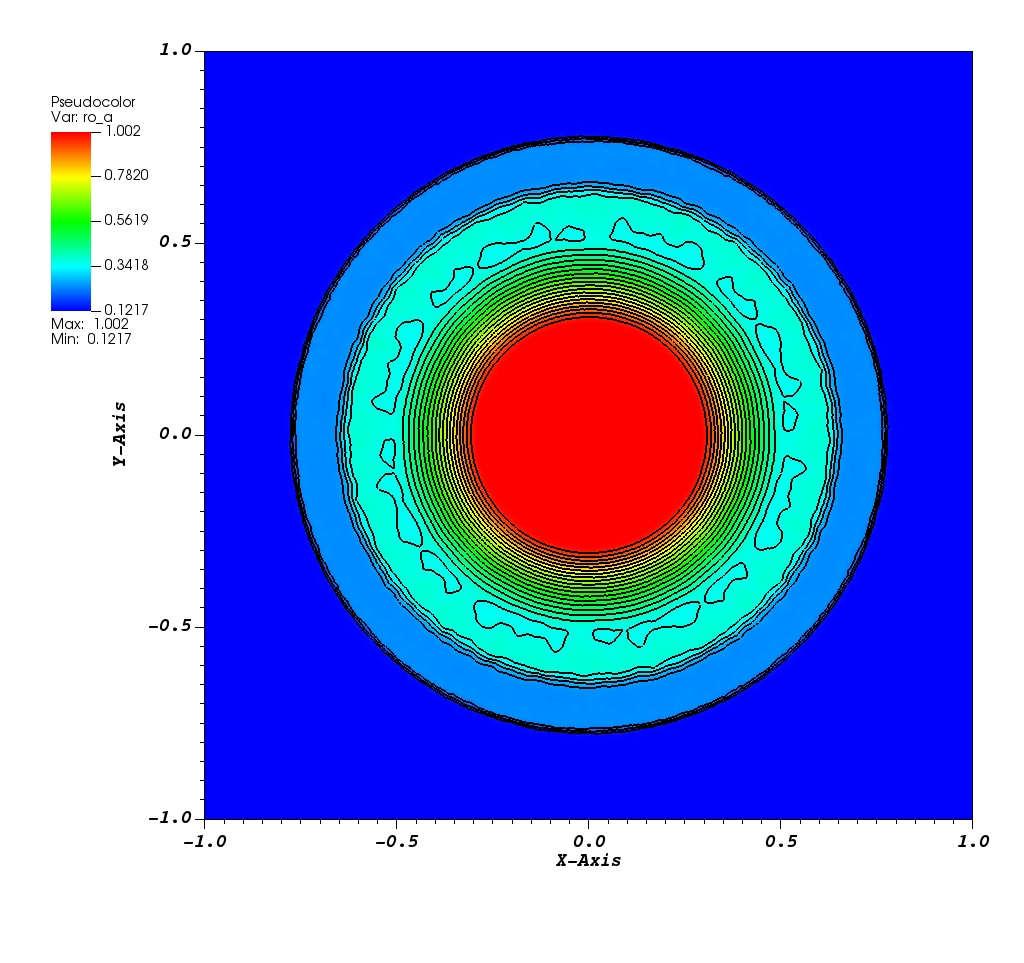}}}
\vskip5pt
\centerline{\subfigure[point values $\rho_\sigma$, cubic]{\includegraphics[trim=1.2cm 2.2cm 1.2cm 1.2cm,clip,width=6.05cm]{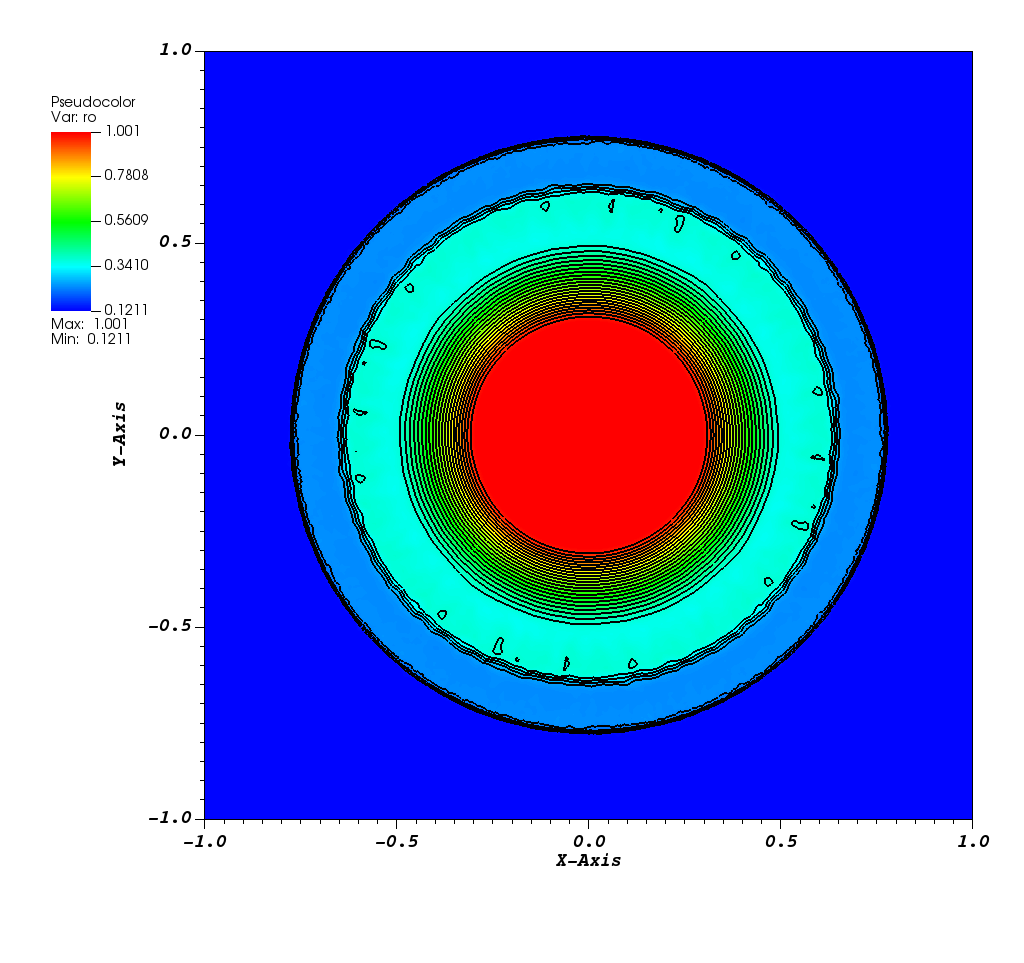}}\hspace*{0.05cm}
\subfigure[average values $\xbar{\rho}_E$, cubic]{\includegraphics[trim=1.2cm 2.2cm 1.2cm 1.2cm,clip,width=6.05cm]{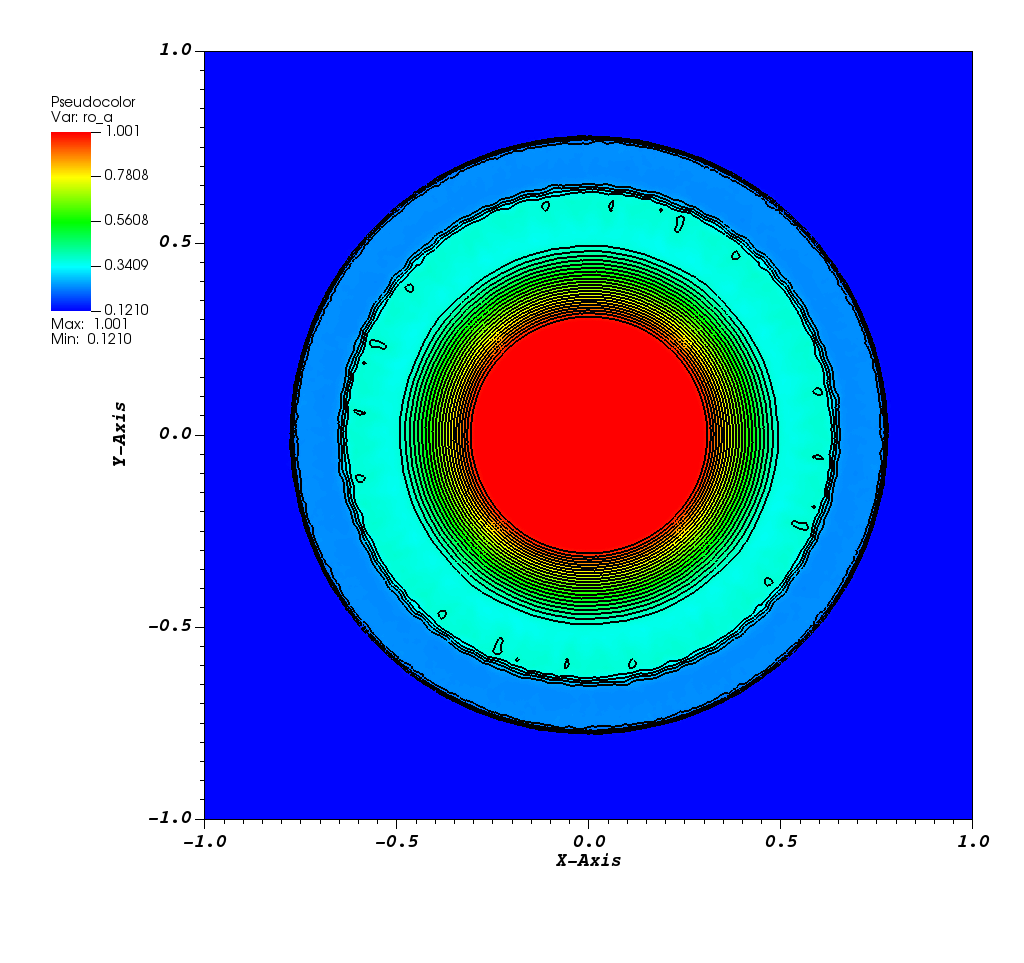}}}
\caption{Sod problem: profiles of density.\label{fig:sod}}
\end{figure}

\begin{figure}[ht!]
\centerline{\subfigure[flag on $\rho_\sigma$, quadratic]{\includegraphics[trim=1.2cm 2.2cm 1.2cm 1.2cm,clip,width=6.05cm]{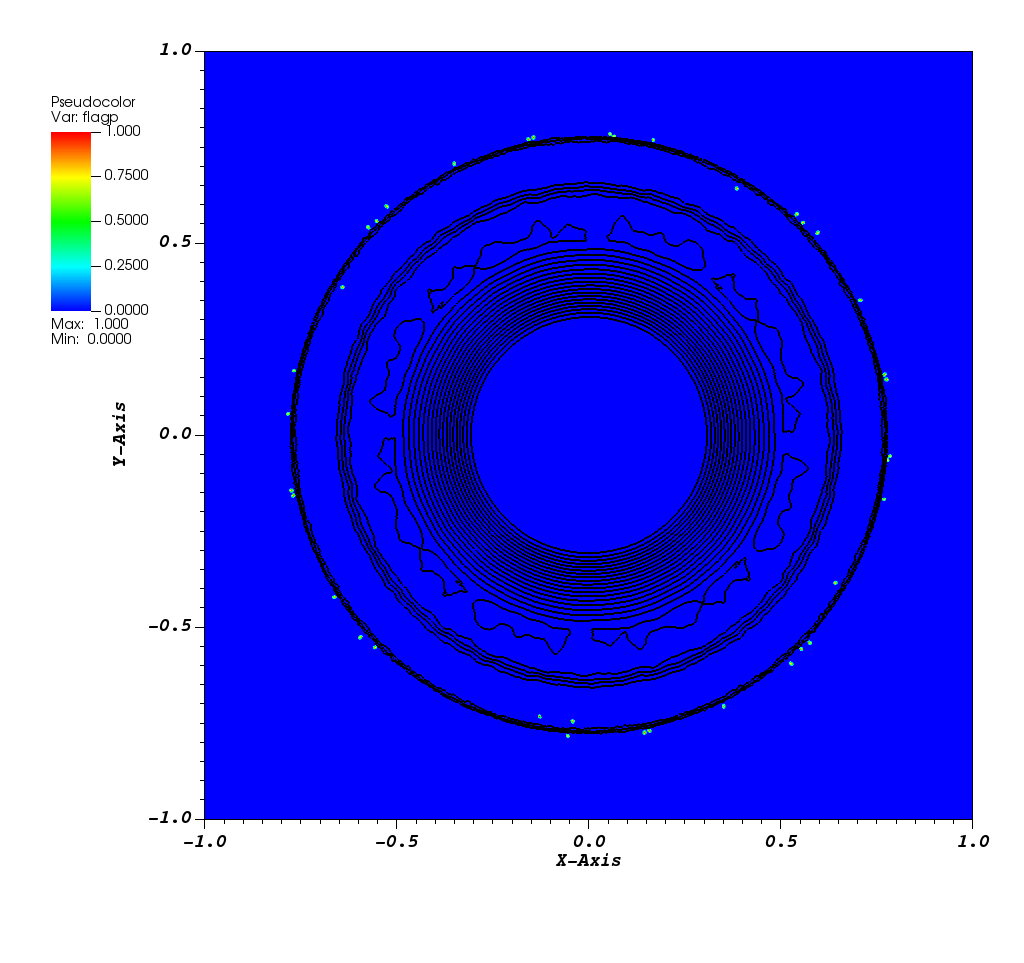}}\hspace*{0.05cm}
\subfigure[flag on $\xbar{\rho}_E$, quadratic]{\includegraphics[trim=1.2cm 2.2cm 1.2cm 1.2cm,clip,width=6.05cm]{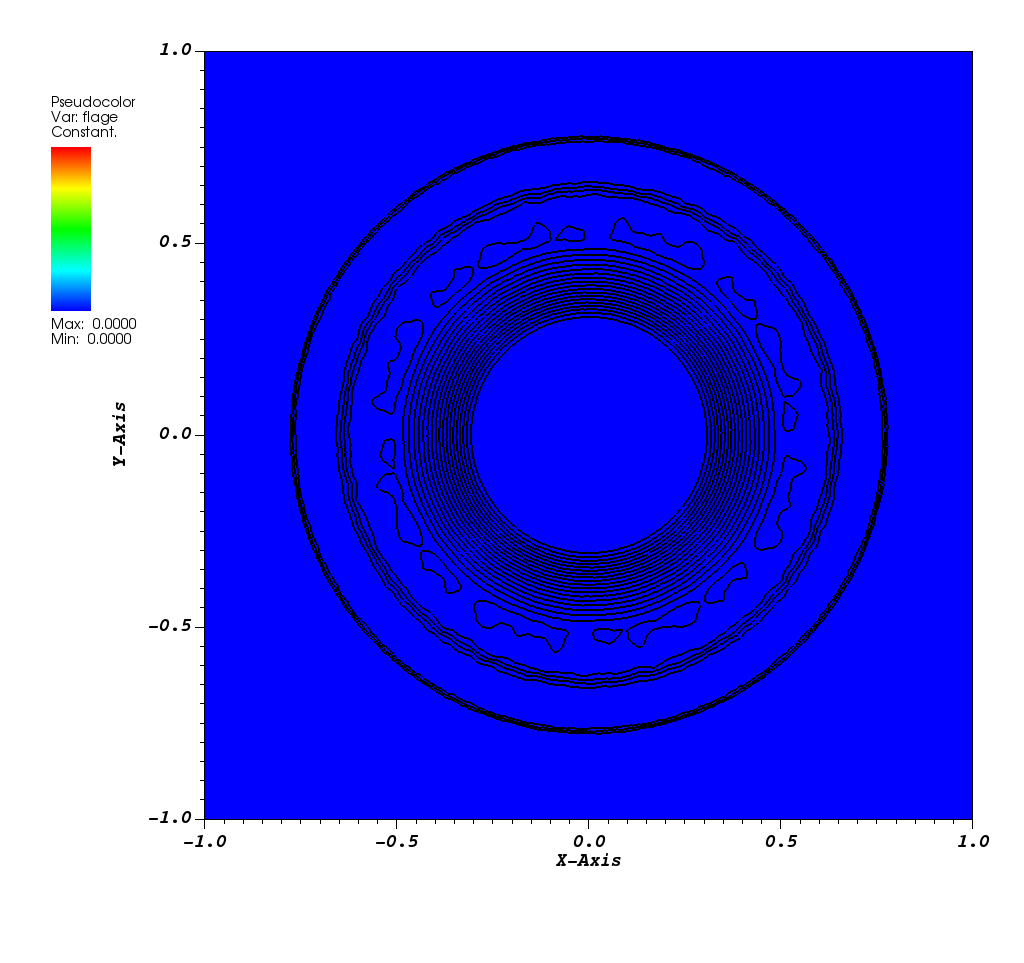}}}
\vskip5pt
\centerline{\subfigure[flag on $\rho_\sigma$, cubic]{\includegraphics[trim=1.2cm 2.2cm 1.2cm 1.2cm,clip,width=6.05cm]{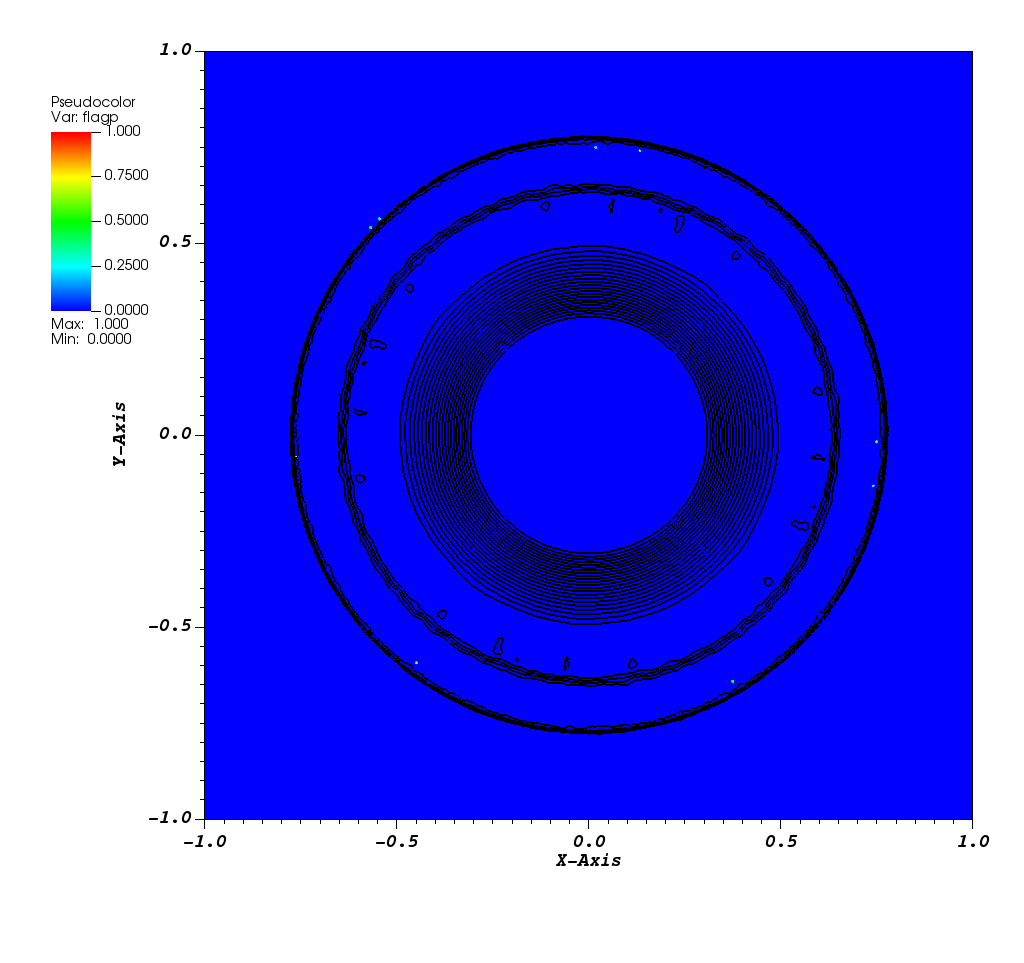}}\hspace*{0.05cm}
\subfigure[flag on $\xbar{\rho}_E$, cubic]{\includegraphics[trim=1.2cm 2.2cm 1.2cm 1.2cm,clip,width=6.05cm]{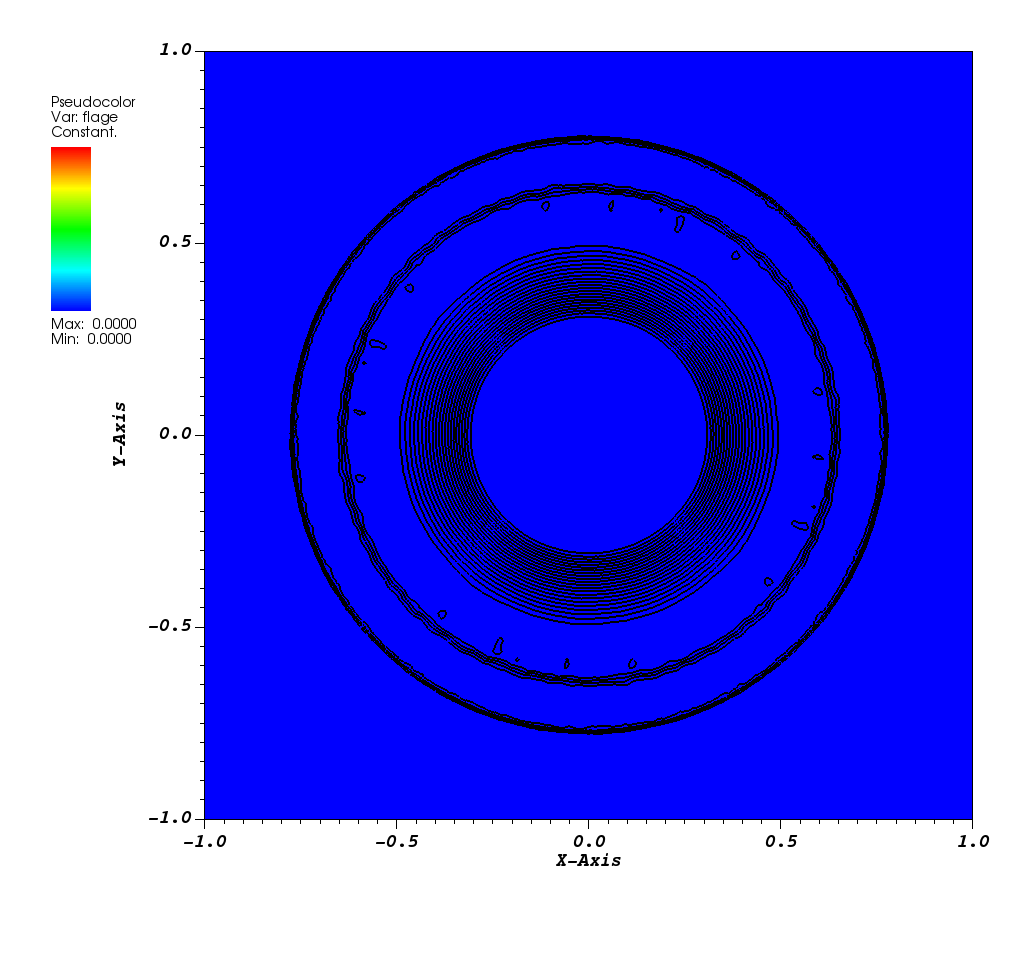}}}
\caption{Sod problem: flags at the positions where the first-order scheme is used.\label{fig:sod_flag}}
\end{figure}

In the third example of non-linear system cases, we test the proposed schemes on a well-known 2D Sod benchmark problem. The initial condition is given by
\begin{equation*}
    (\rho,u,v,p)=\left\{\begin{aligned}
         &(1,0,0,1.5) &&\mbox{if}~\|\mathbf{x}\|^2\leq0.25,\\
         &(0.125,0,0,0.1)&&\mbox{else},
    \end{aligned}\right.
\end{equation*}
the boundary condition is solid wall, and the final time we compute is $T=0.16$. The two approximations (quadratic and cubic) have been used to compute the solutions. We plot the results (density component) in Figure \ref{fig:sod}. As one can see, both two approximations produce similar results and cubic approximation is much more accurate as expected. In Figure \ref{fig:sod_flag}, we also plot the flags where the first-order scheme is activated at the final time. We can clearly see that only the first-order scheme for the point values is used and the cubic approximation leads to less triggers of the first-order scheme.  

\subsubsection{Liu-Lax problem}\label{LL}
\begin{figure}[ht!]
\centerline{\subfigure[point values $\rho_\sigma$, quadratic]{\includegraphics[trim=1.2cm 2.2cm 1.2cm 1.2cm,clip,width=6.05cm]{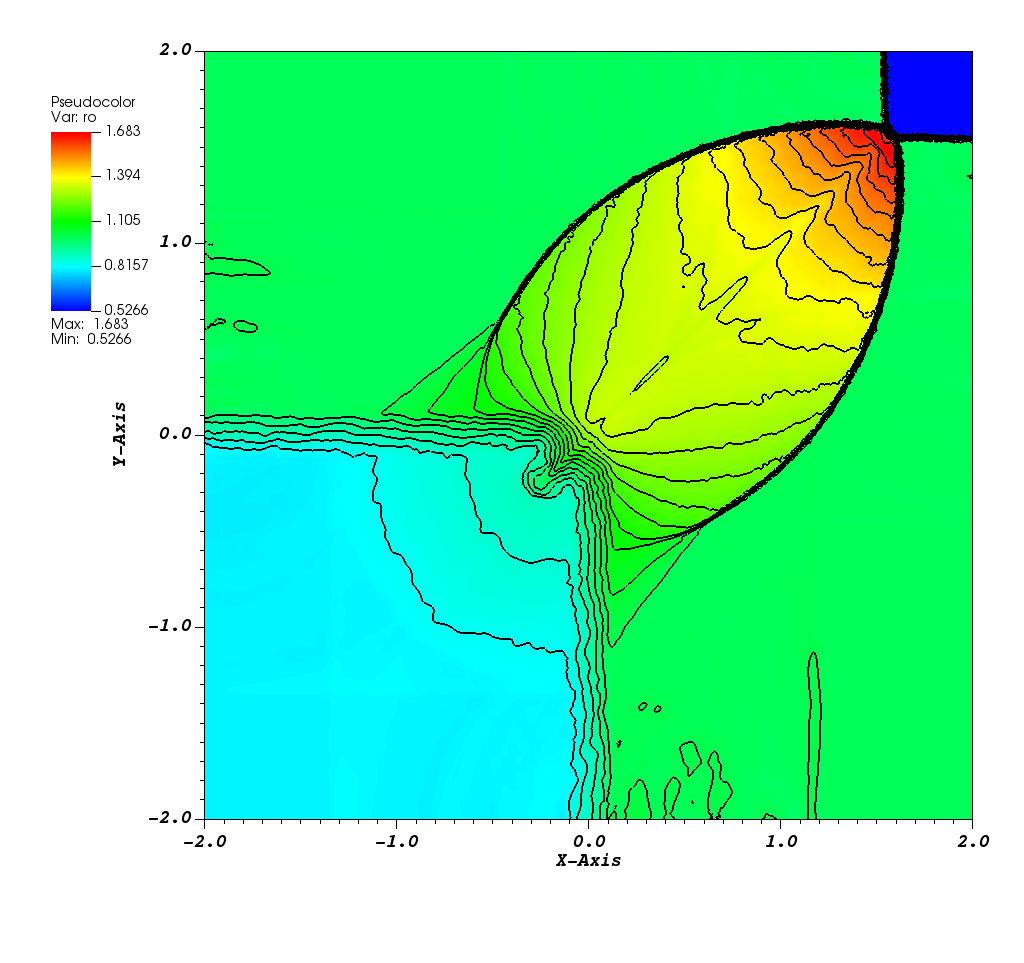}}\hspace*{0.05cm}
\subfigure[average values $\xbar{\rho}_E$, quadratic]{\includegraphics[trim=1.2cm 2.2cm 1.2cm 1.2cm,clip,width=6.05cm]{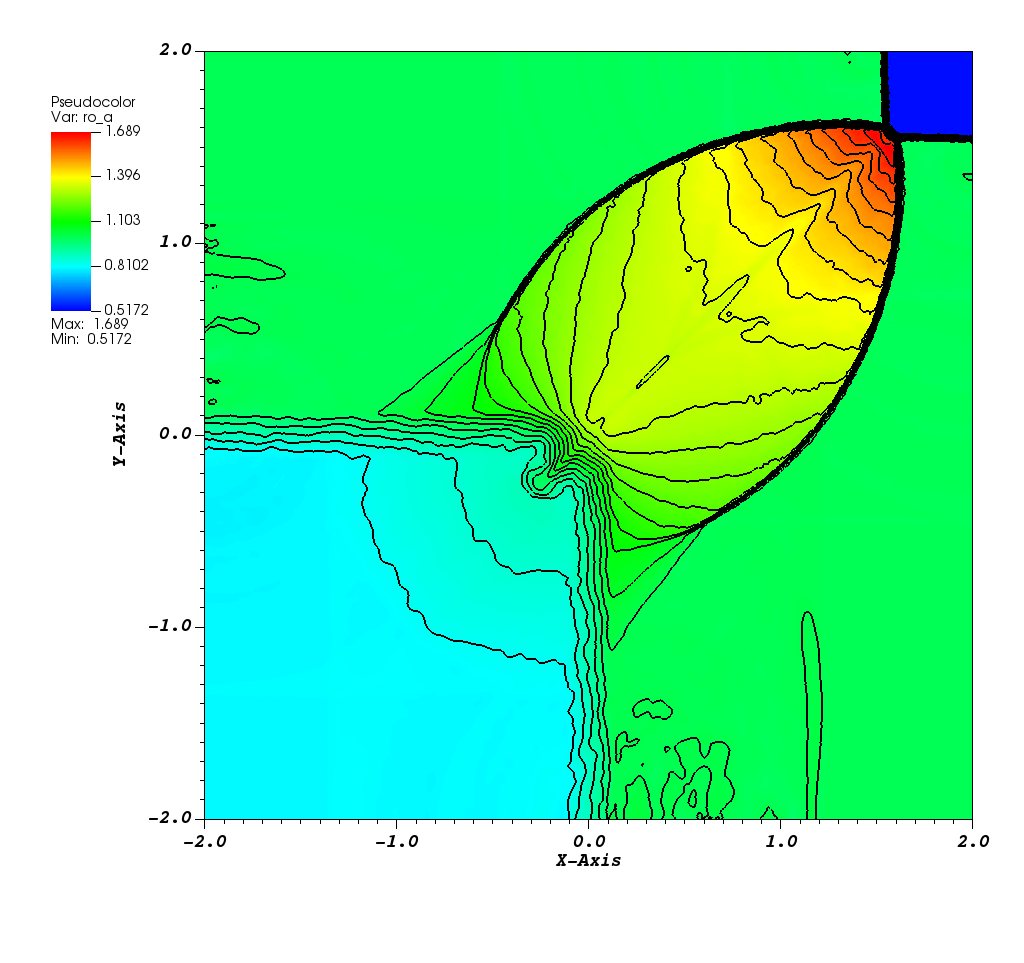}}}
\vskip5pt
\centerline{\subfigure[point values $\rho_\sigma$, cubic]{\includegraphics[trim=1.2cm 2.2cm 1.2cm 1.2cm,clip,width=6.05cm]{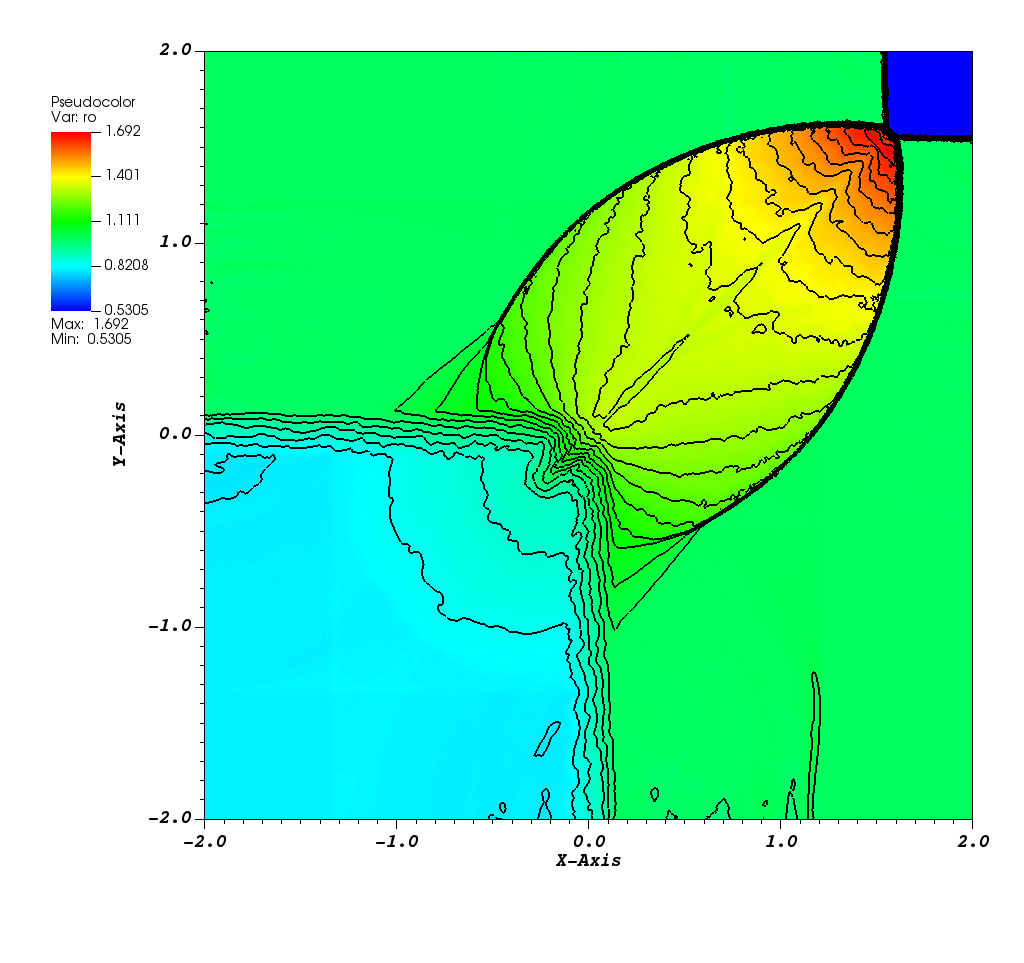}}\hspace*{0.05cm}
\subfigure[average values $\xbar{\rho}_E$, cubic]{\includegraphics[trim=1.2cm 2.2cm 1.2cm 1.2cm,clip,width=6.05cm]{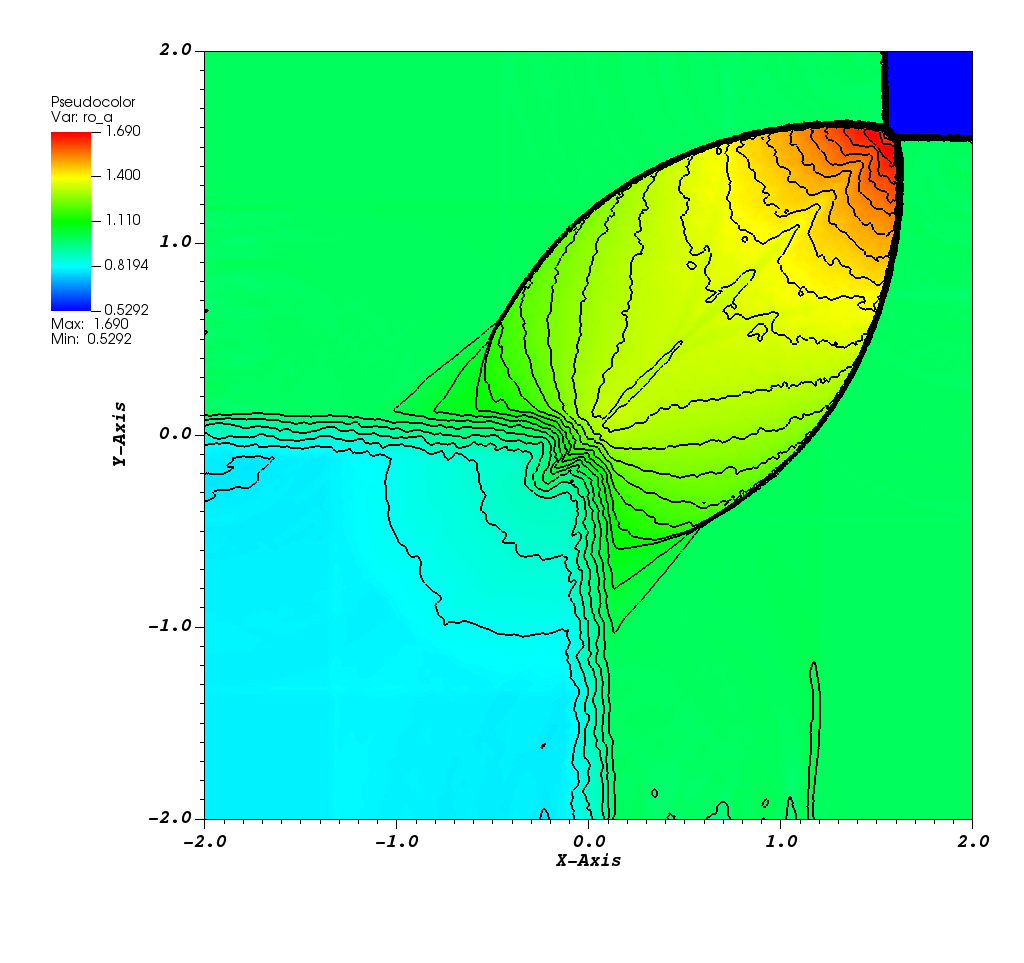}}}
\caption{Liu-Lax case: profiles of density.\label{LL13_ro}}
\end{figure}

In the fourth example of non-linear system cases, we consider the following initial condition
$$ (\rho, u,v,p)=\left \{\begin{array}{ll}
(\rho_1,u_1,v_1,p_1)=(0.5313,0,0,0.4)& \text{ if } x\geq0\text{ and } y\geq 0,\\
(\rho_2,u_2,v_2,p_2)=(1,0.7276,0,1) & \text{ if } x\leq 0 \text{ and } y\geq 0,\\
(\rho_3,u_3,v_3,p_3)=(0.8,0,0,1)&\text{ if } x\leq 0\text{ and }y\leq 0,\\
(\rho_4,u_4,v_4,p_4)=(1,0, 0.7276,1) &\text{ if } x\geq0\text{ and } y\leq 0,
\end{array}\right .
$$
prescribed in the computational domain $[-2,2]^2$. The states 1 and 2 are separated by a shock. The states 2 and 3 are separated by a slip line. The states 3 and 4 are separated by a steady slip-line. The states  4 and 1 are separated by a shock. The mesh has $53227$ vertices, $158878$ edges, $105652$ elements, and $212105$ (resp. $476635$) DoFs for the quadratic (resp. cubic) approximation. The obtained results of density field are plotted in Figure \ref{LL13_ro}.

\subsubsection{Kurganov-Tadmor problem} 
\begin{figure}[ht!]
\centerline{\subfigure[point values $\rho_\sigma$, quadratic]{\includegraphics[trim=1.2cm 2.2cm 1.2cm 1.2cm,clip,width=6.05cm]{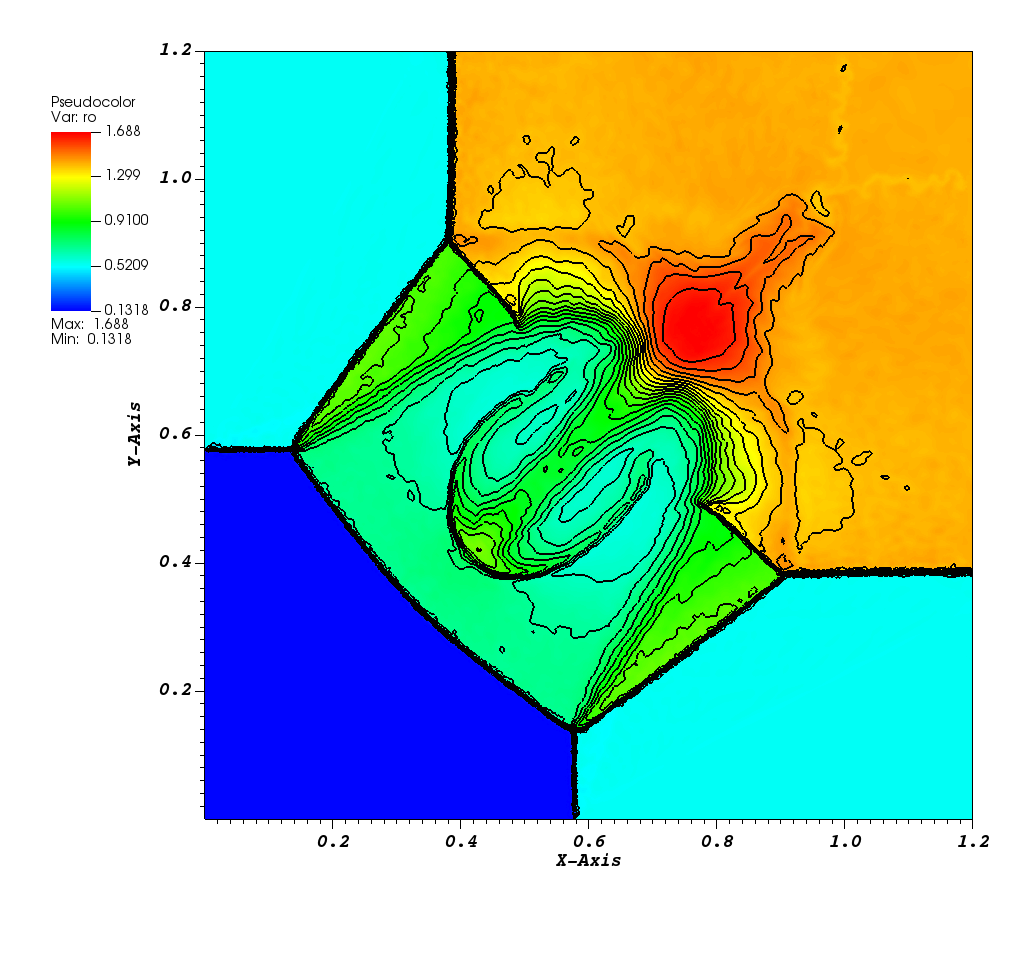}}\hspace*{0.05cm}
\subfigure[average values $\xbar{\rho}_E$, quadratic]{\includegraphics[trim=1.2cm 2.2cm 1.2cm 1.2cm,clip,width=6.05cm]{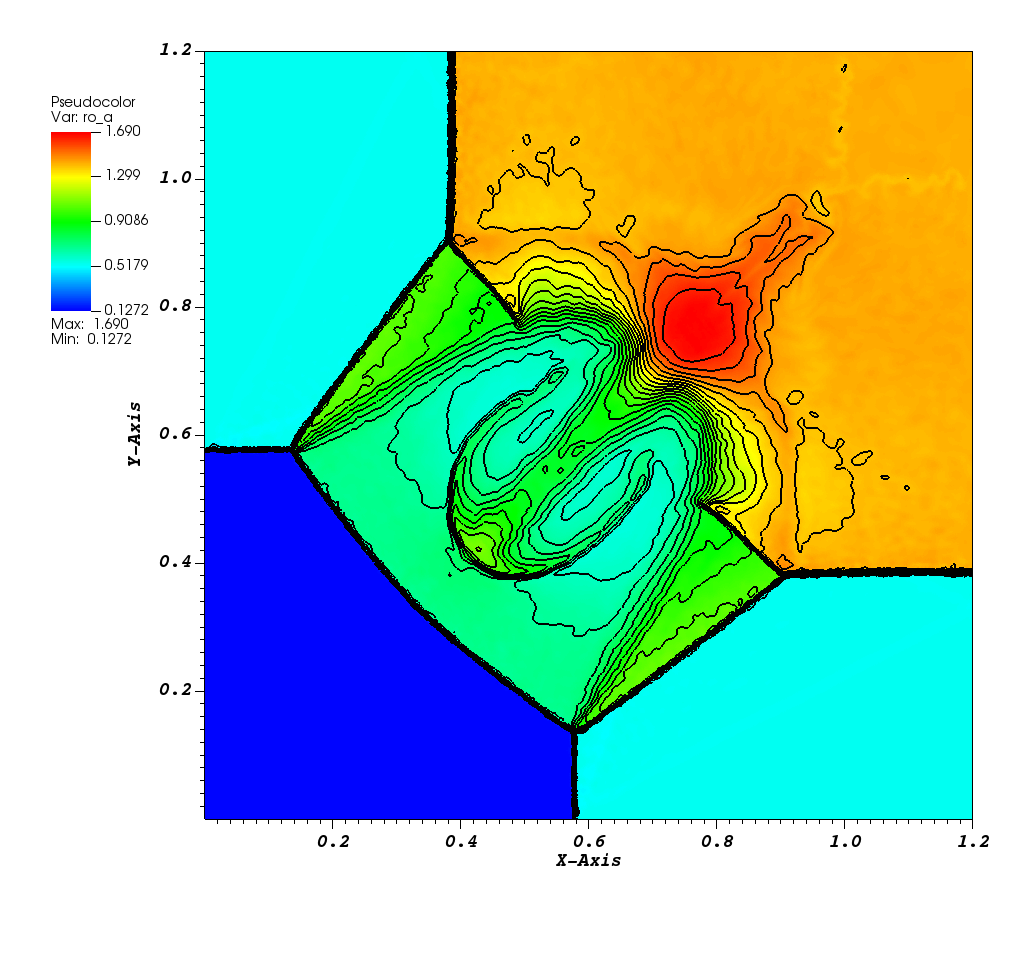}}}
\vskip5pt
\centerline{\subfigure[point values $\rho_\sigma$, cubic]{\includegraphics[trim=1.2cm 2.2cm 1.2cm 1.2cm,clip,width=6.05cm]{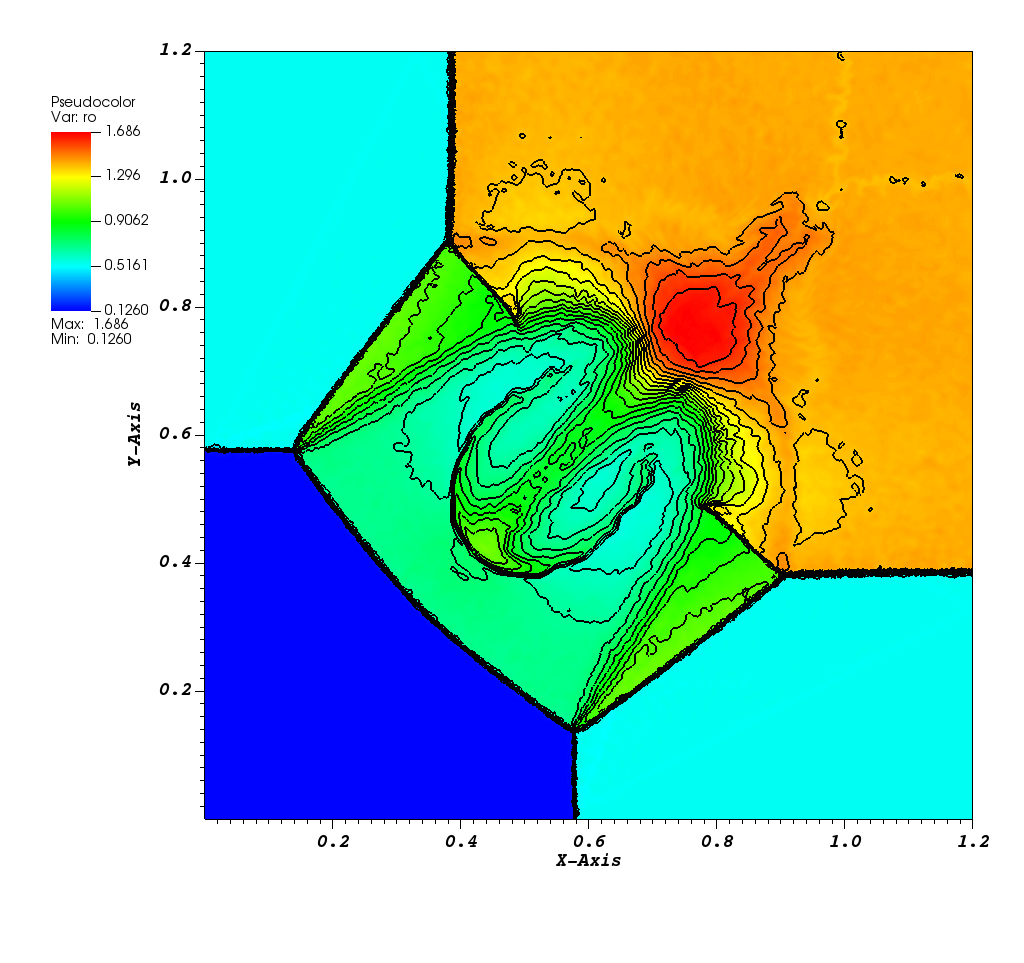}}\hspace*{0.05cm}
\subfigure[average values $\xbar{\rho}_E$, cubic]{\includegraphics[trim=1.2cm 2.2cm 1.2cm 1.2cm,clip,width=6.05cm]{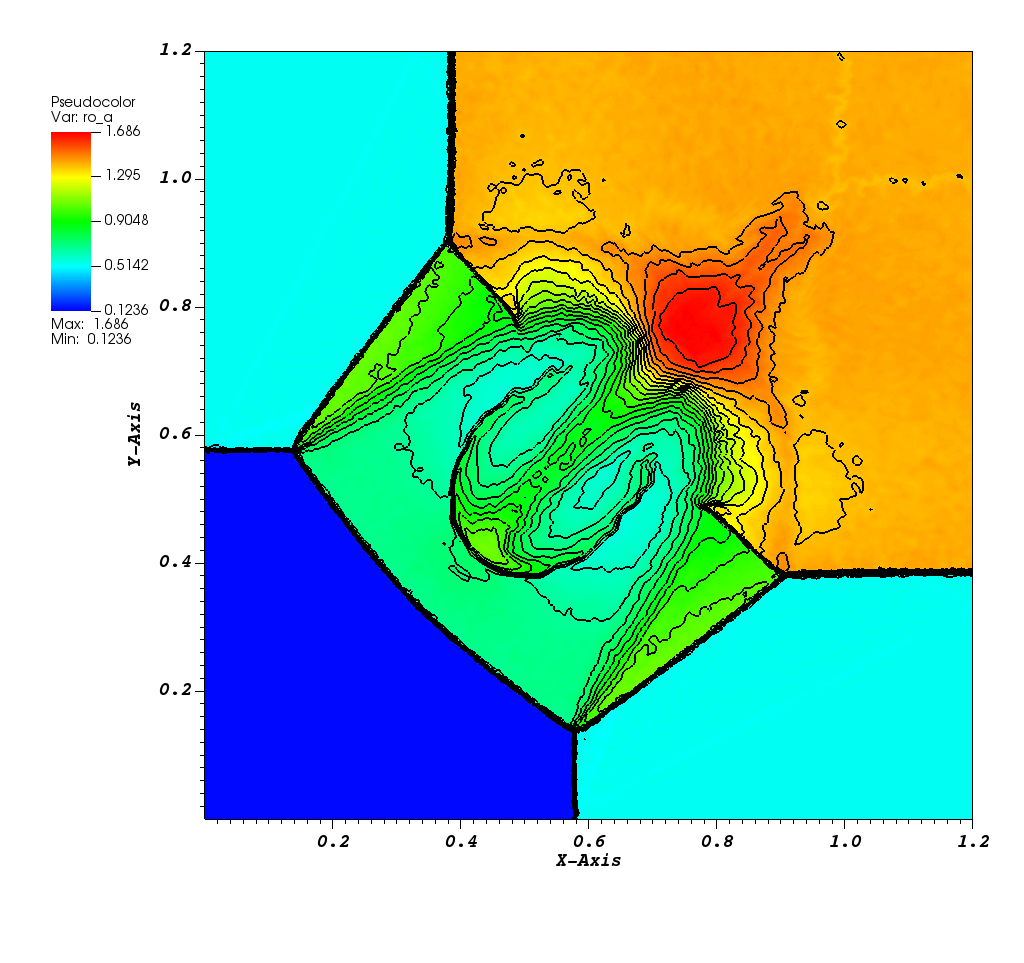}}}
\caption{Kurganov-Tadmor problem: profiles of density.\label{fig:KT_ro}}
\end{figure}

\begin{figure}[ht!]
\centerline{\subfigure[flag on $\rho_\sigma$, quadratic]{\includegraphics[trim=1.2cm 2.2cm 1.2cm 1.2cm,clip,width=6.05cm]{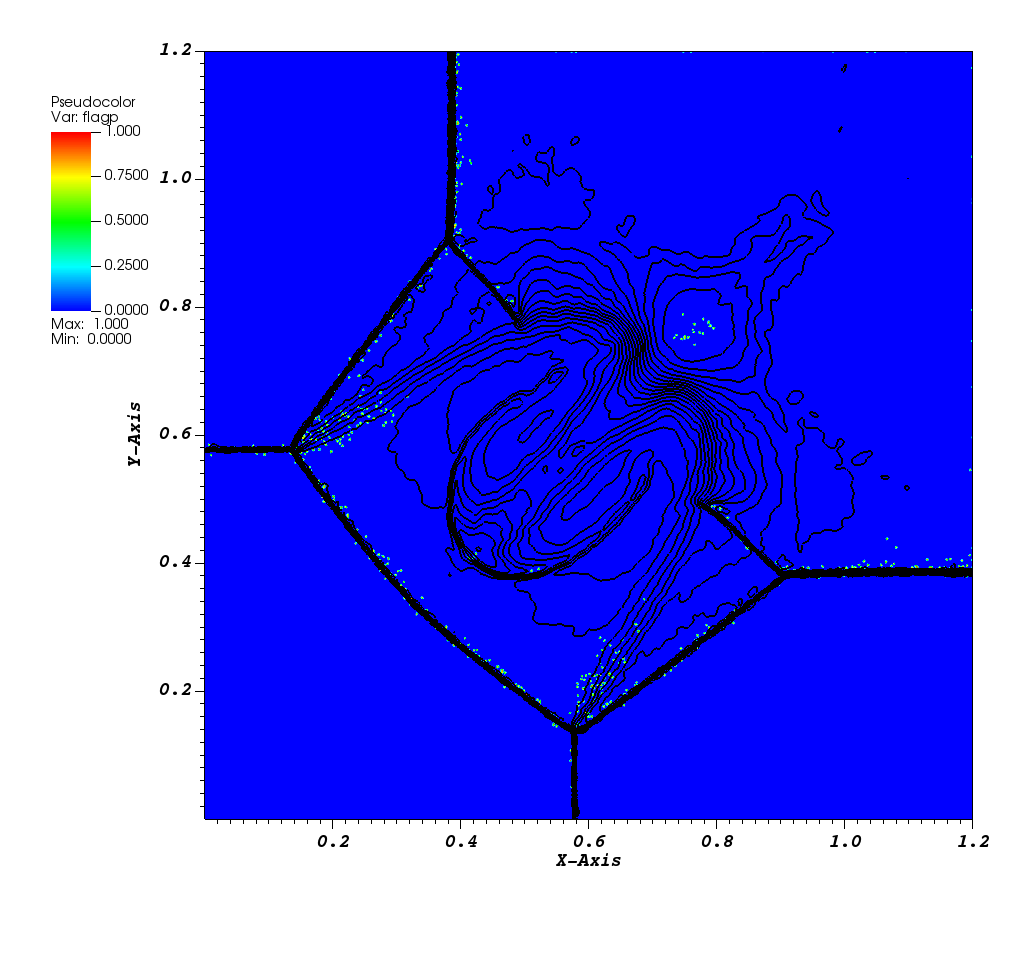}}\hspace*{0.05cm}
\subfigure[flag on $\xbar{\rho}_E$, quadratic]{\includegraphics[trim=1.2cm 2.2cm 1.2cm 1.2cm,clip,width=6.05cm]{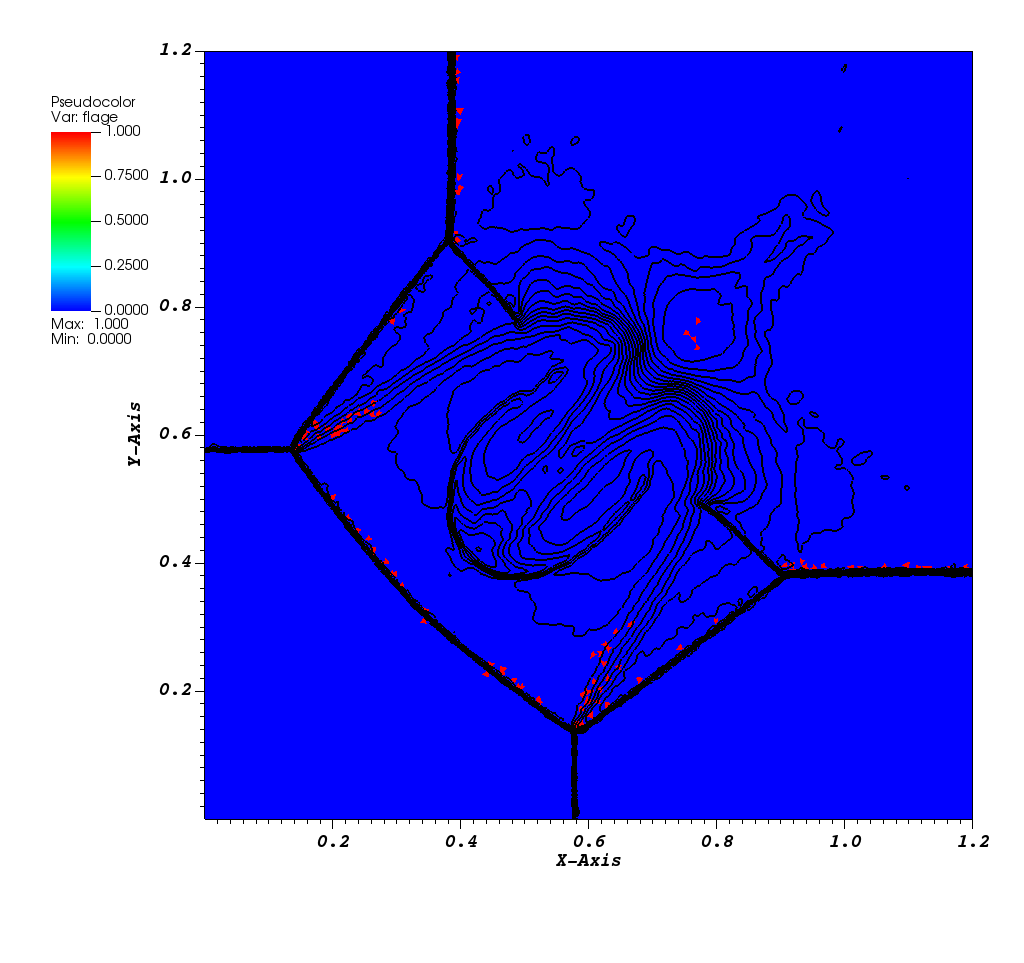}}}
\vskip5pt
\centerline{\subfigure[flag on $\rho_\sigma$, cubic]{\includegraphics[trim=1.2cm 2.2cm 1.2cm 1.2cm,clip,width=6.05cm]{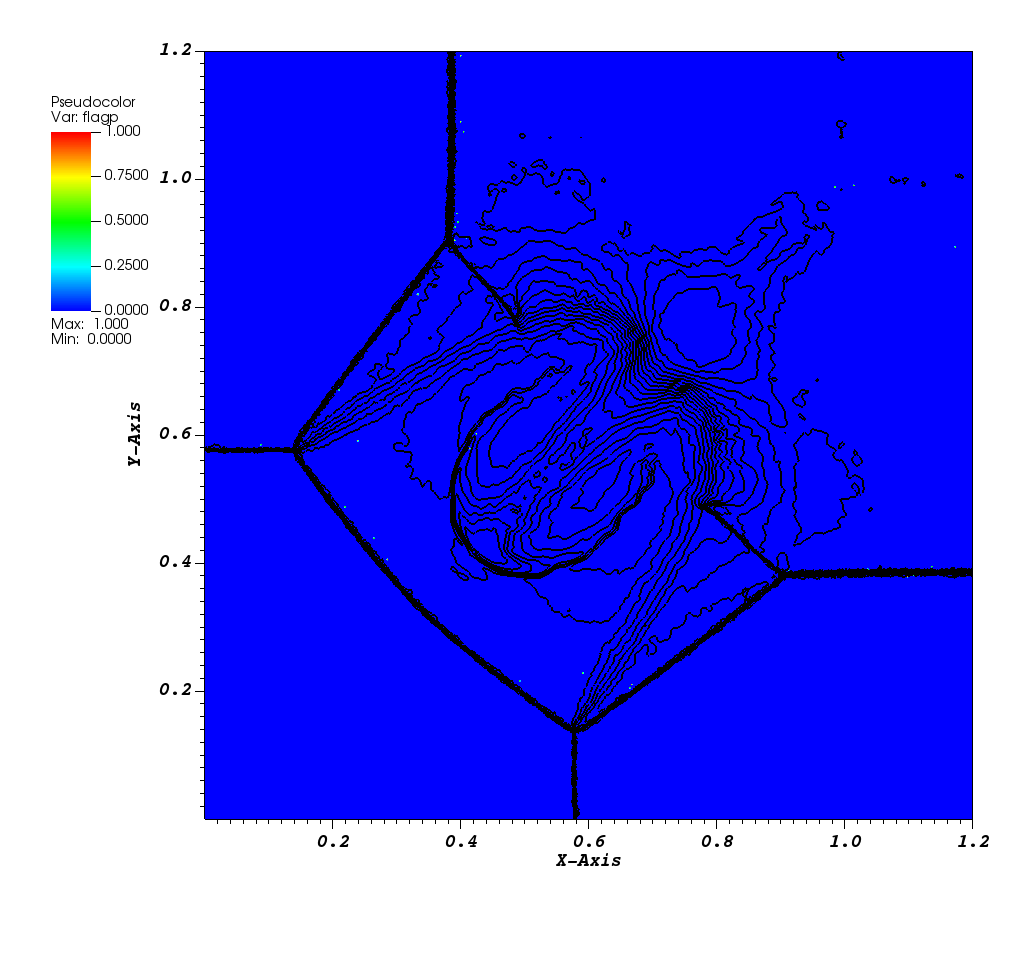}}\hspace*{0.05cm}
\subfigure[flag on $\xbar{\rho}_E$, cubic]{\includegraphics[trim=1.2cm 2.2cm 1.2cm 1.2cm,clip,width=6.05cm]{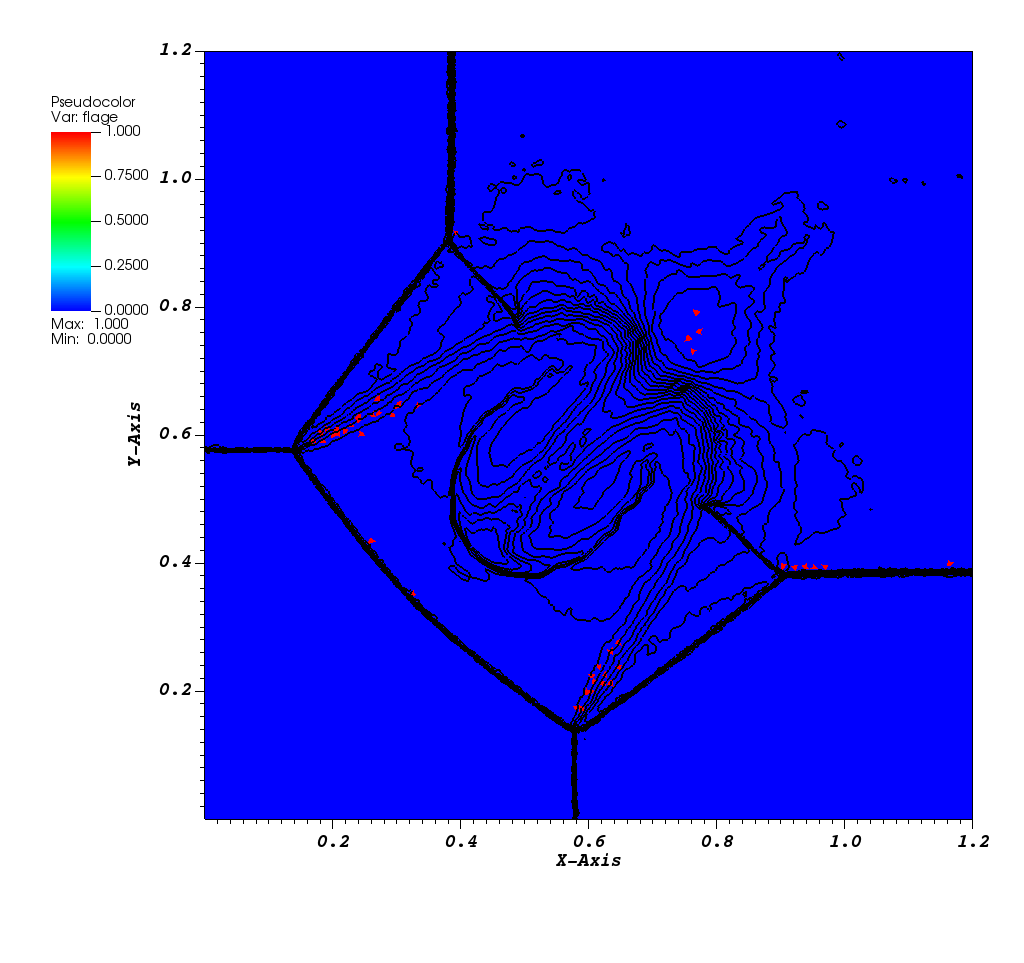}}}
\caption{Kurganov-Tadmor problem: flags at the positions where the first-order scheme is used.\label{fig:KT_flag}}
\end{figure}

In the fifth example of non-linear system cases, the initial condition is 
$$ (\rho, u,v,p)=\left \{\begin{array}{ll}
(\rho_1,u_1,v_1,p_1)=(1.5, 0, 0, 1.5)& \text{ if } x\geq1\text{ and } y\geq 1,\\
(\rho_2,u_2,v_2,p_2)=(0.5323, 1.206, 0, 0.3) & \text{ if } x\leq 1 \text{ and } y\geq 1,\\
(\rho_3,u_3,v_3,p_3)=(0.138, 1.206, 1.206, 0.029)&\text{ if } x\leq 1\text{ and }y\leq 1,\\
(\rho_4,u_4,v_4,p_4)=(0.5323, 0, 1.206, 0.3) &\text{ if } x\leq1\text{ and } y\leq 1.
\end{array}\right .
$$
Here the four states are separated by shocks. The domain is $[0,1.2]^2$. The solution at $T=1$ is displayed in Figure \ref{fig:KT_ro}.
The mesh has $53140$ vertices, $158617$ edges, $105478$ elements, and $211757$ (resp. $475852$) DoFs for the quadratic (resp. cubic) approximation. We see that, on the one hand, the solution looks very similar to what was obtained in the literature, see e.g. in \cite{KLZ, GKL, WDKL}. On the other hand, the first-order schemes (for average and point values) are activated on very few positions, see Figure \ref{fig:KT_flag}.

\subsubsection{Woodward-Collela problem} 
\begin{figure}[ht!]
\centerline{\subfigure[point values $\rho_\sigma$, quadratic]{\includegraphics[trim=1.2cm 2.2cm 1.2cm 1.2cm,clip,width=6.05cm]{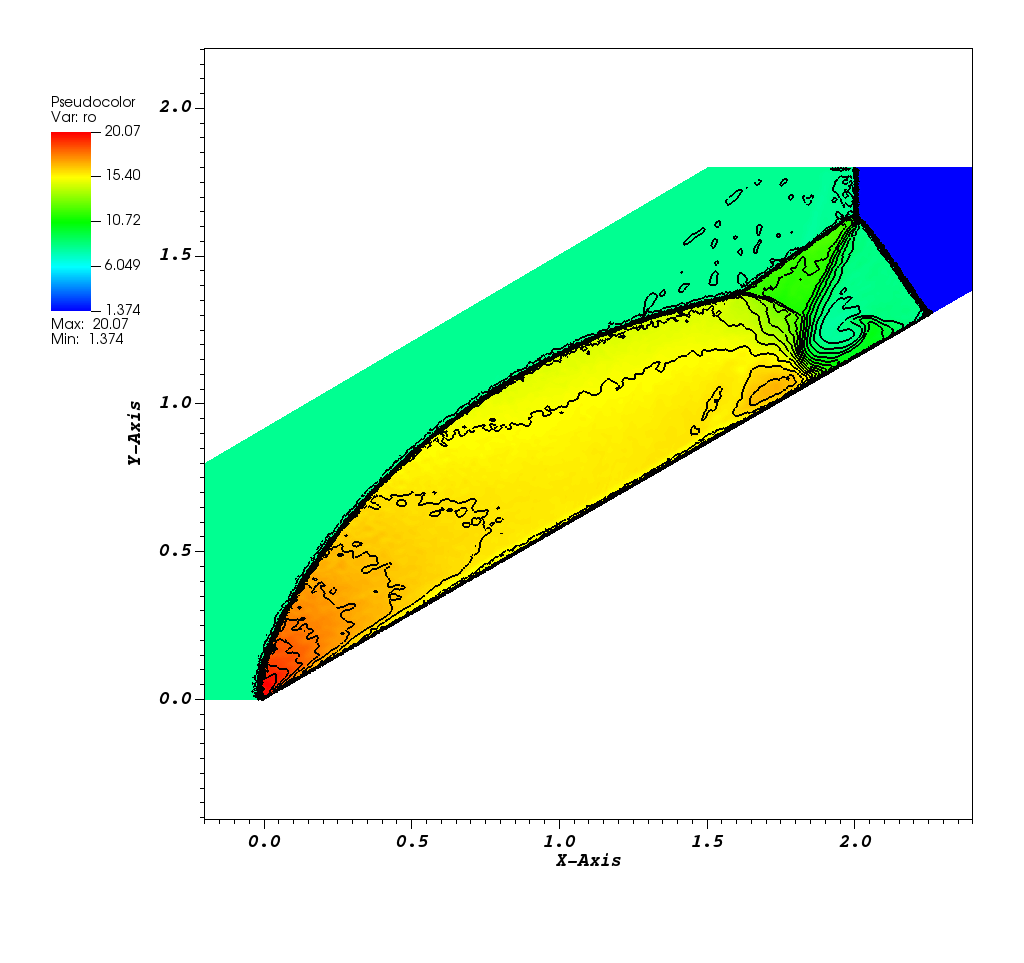}}\hspace*{0.05cm}
\subfigure[average values $\xbar{\rho}_E$, quadratic]{\includegraphics[trim=1.2cm 2.2cm 1.2cm 1.2cm,clip,width=6.05cm]{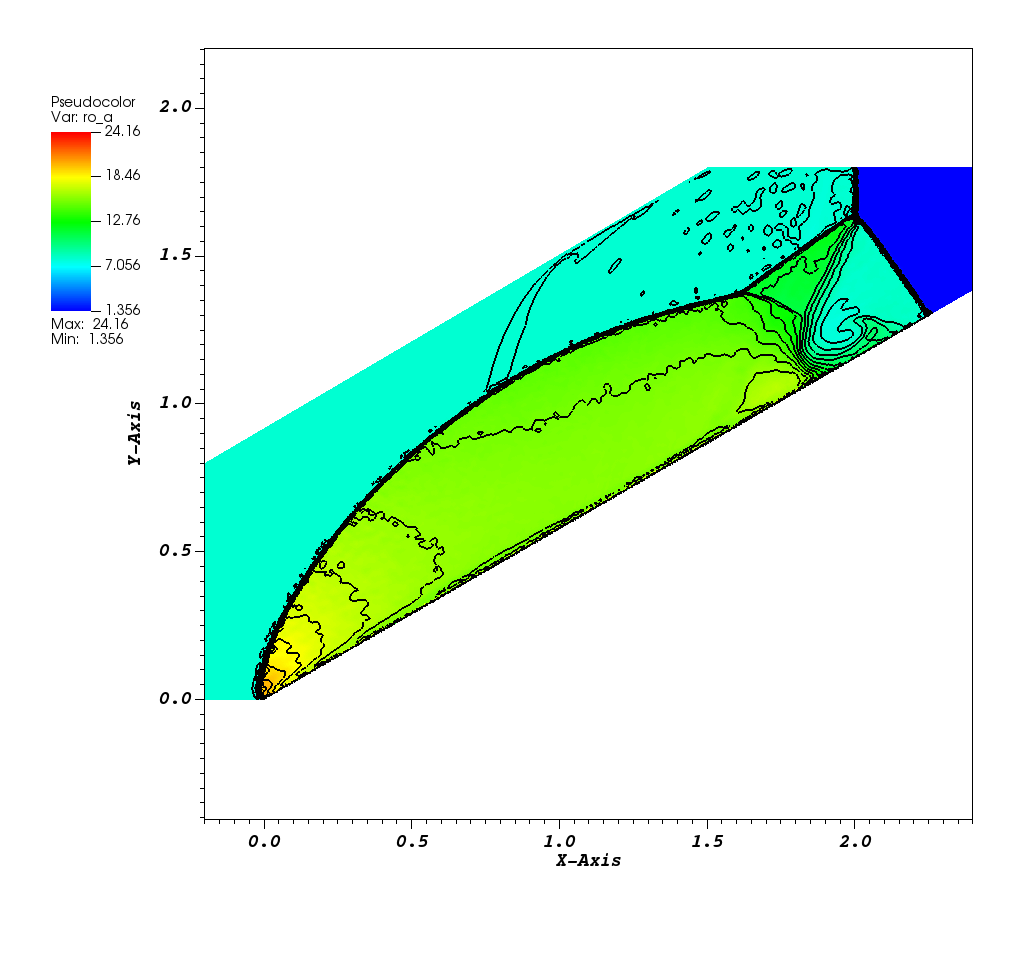}}}
\vskip5pt
\centerline{\subfigure[point values $\rho_\sigma$, cubic]{\includegraphics[trim=1.2cm 2.2cm 1.2cm 1.2cm,clip,width=6.05cm]{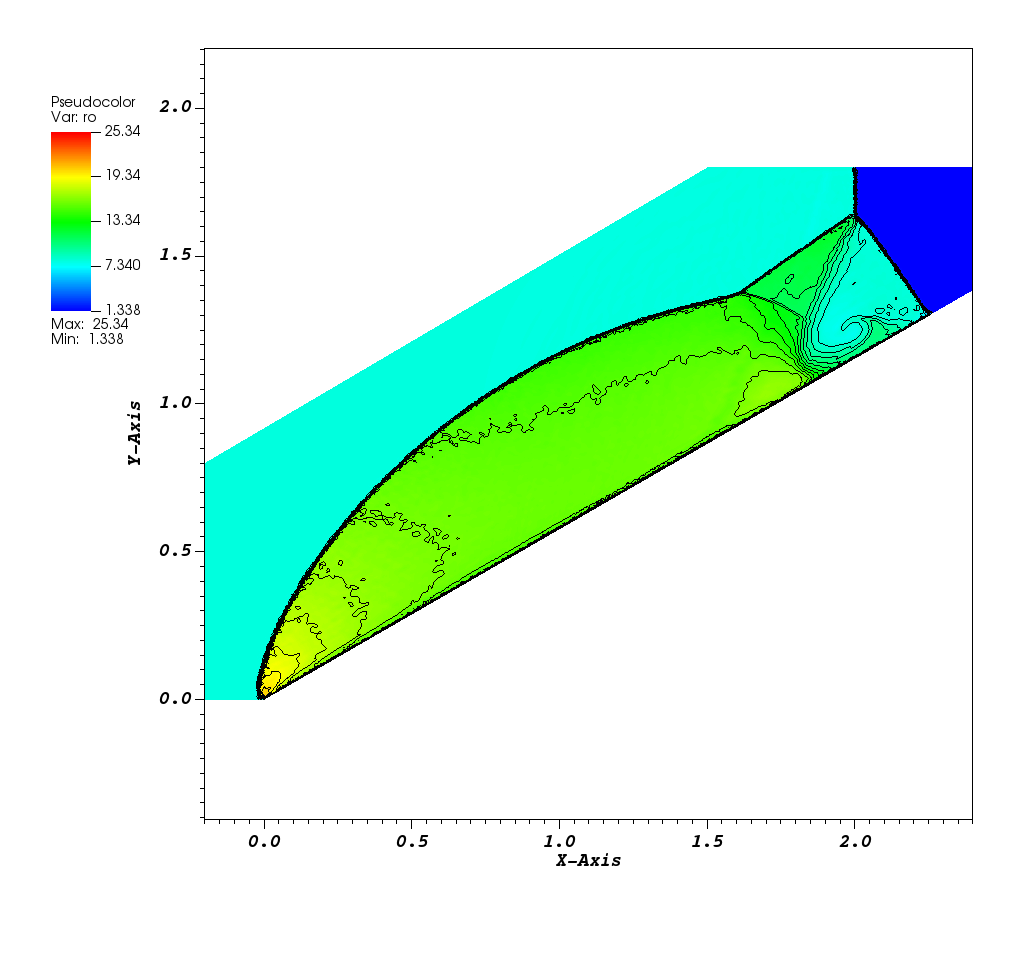}}\hspace*{0.05cm}
\subfigure[average values $\xbar{\rho}_E$, cubic]{\includegraphics[trim=1.2cm 2.2cm 1.2cm 1.2cm,clip,width=6.05cm]{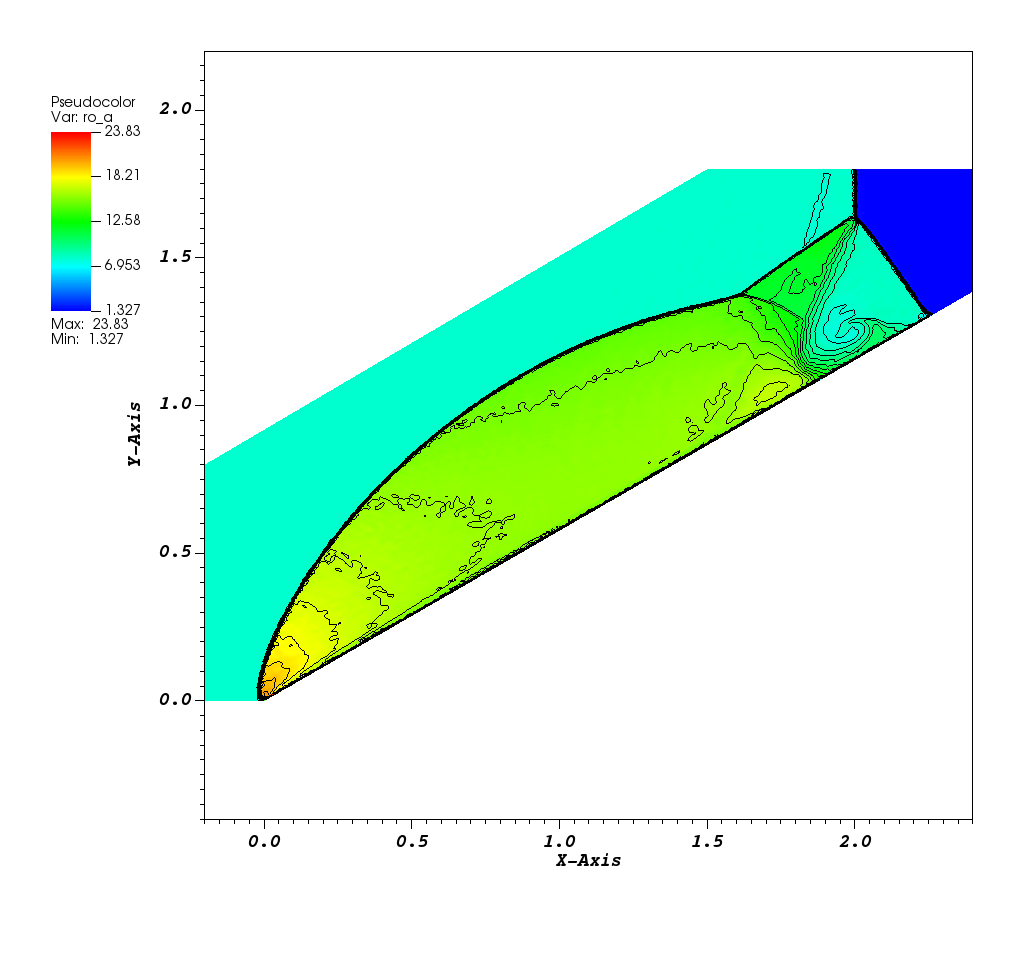}}}
\caption{Woodward-Collela problem: profiles of density.\label{fig:DMR_ro}}
\end{figure}

\begin{figure}[ht!]
\centerline{\subfigure[flag on $\rho_\sigma$, quadratic]{\includegraphics[trim=1.2cm 2.2cm 1.2cm 1.2cm,clip,width=6.05cm]{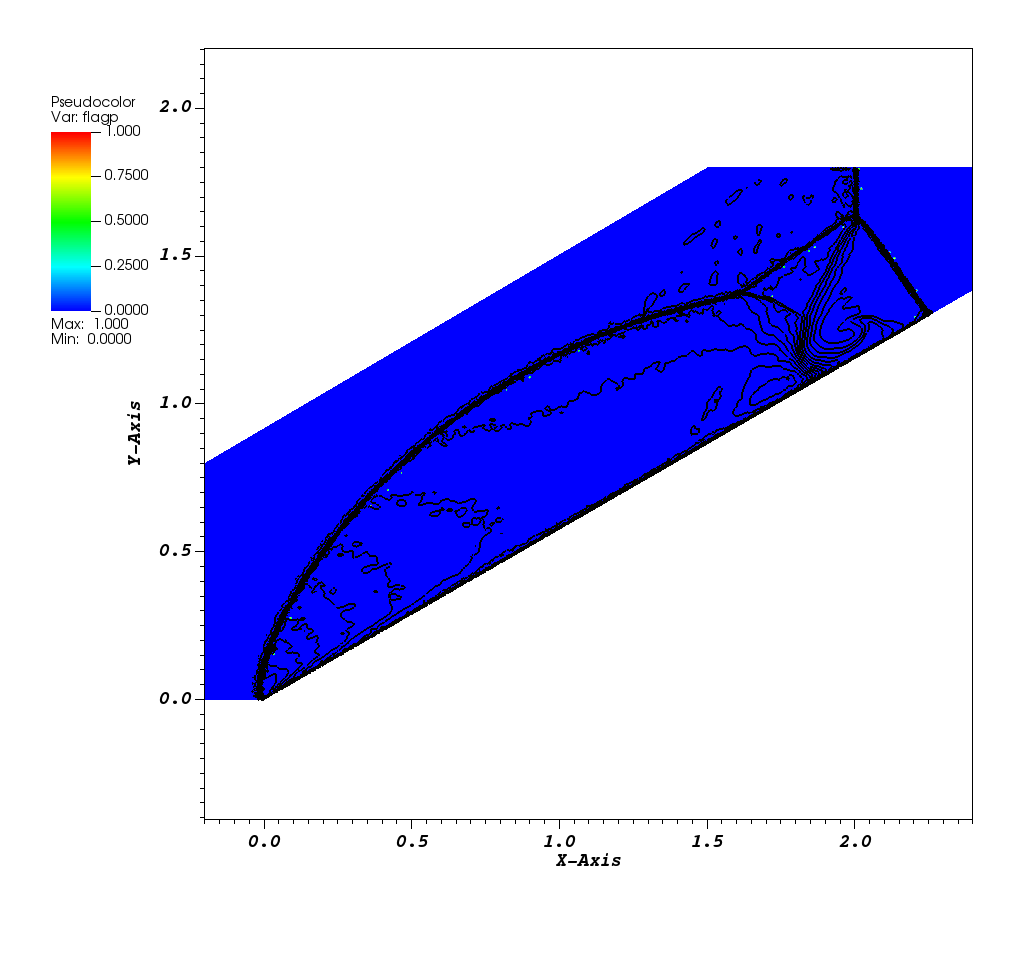}}\hspace*{0.05cm}
\subfigure[flag on $\xbar{\rho}_E$, quadratic]{\includegraphics[trim=1.2cm 2.2cm 1.2cm 1.2cm,clip,width=6.05cm]{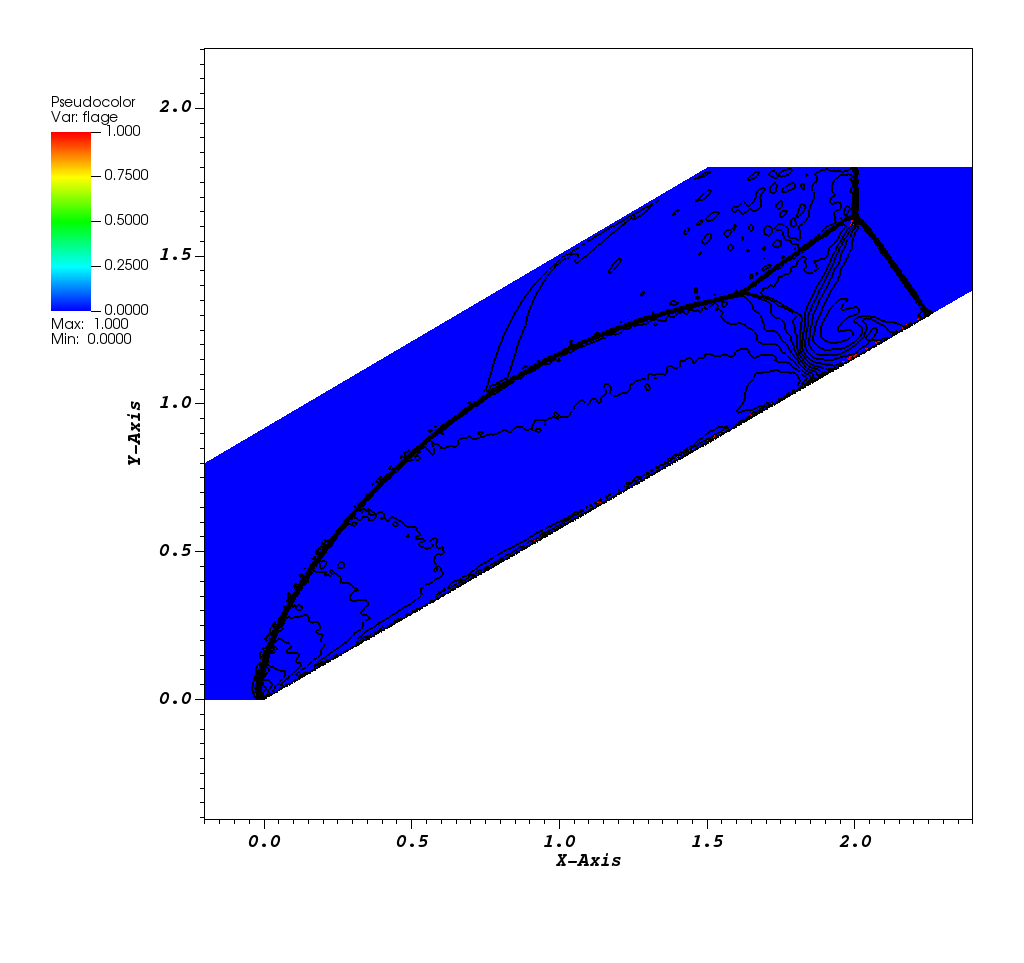}}}
\vskip5pt
\centerline{\subfigure[flag on $\rho_\sigma$, cubic]{\includegraphics[trim=1.2cm 2.2cm 1.2cm 1.2cm,clip,width=6.05cm]{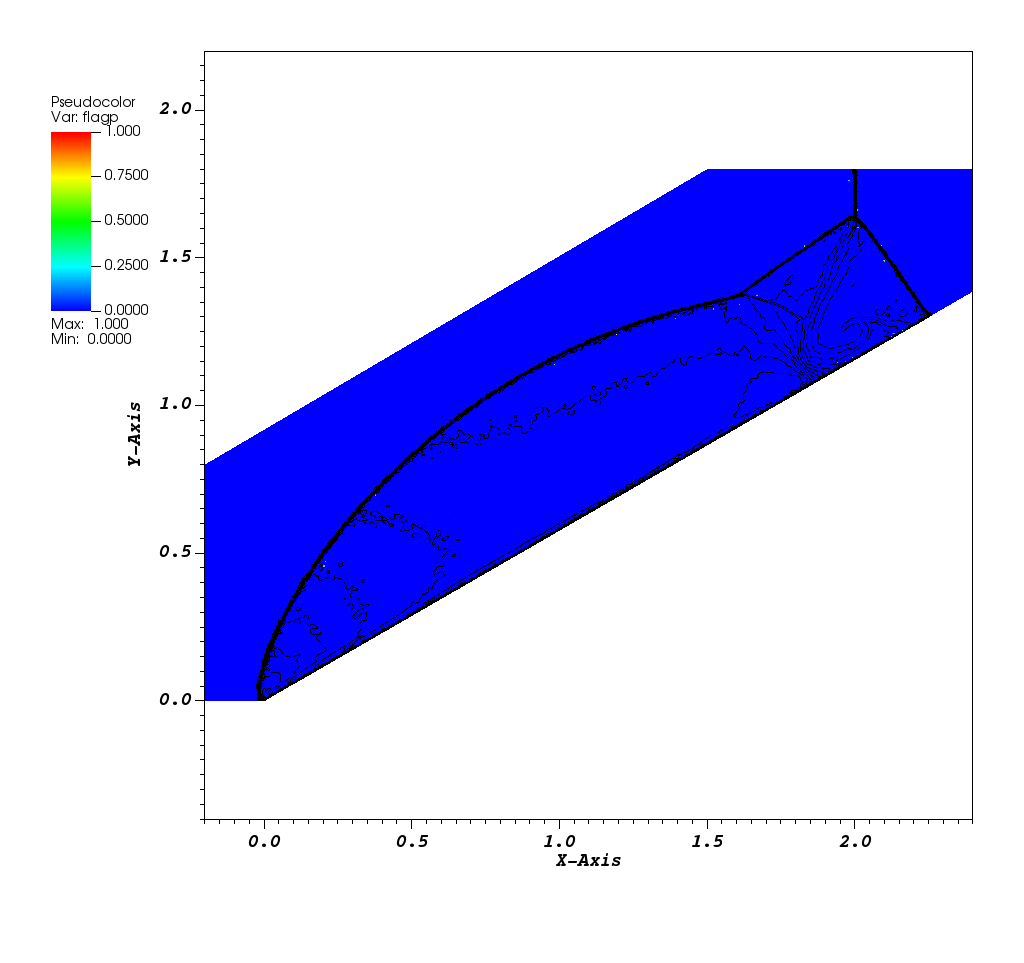}}\hspace*{0.05cm}
\subfigure[flag on $\xbar{\rho}_E$, cubic]{\includegraphics[trim=1.2cm 2.2cm 1.2cm 1.2cm,clip,width=6.05cm]{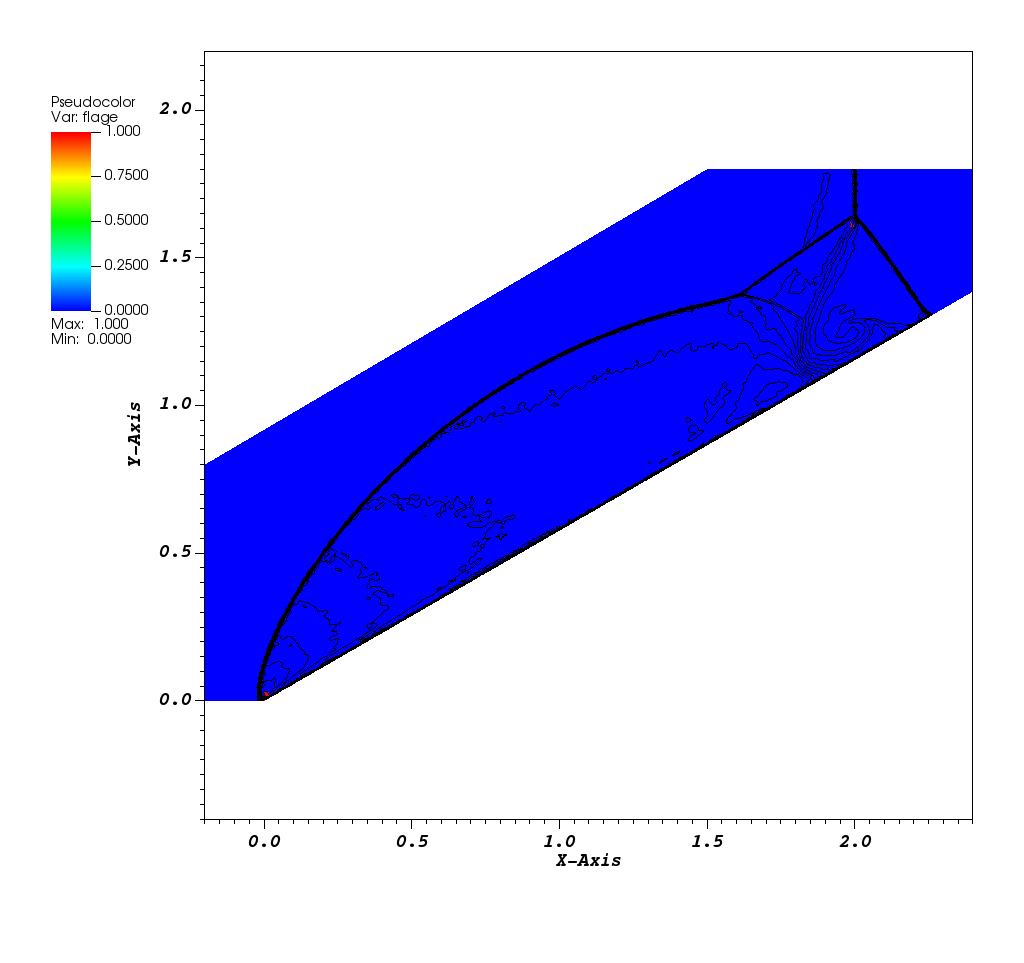}}}
\caption{Woodward-Collela problem: flags at the positions where the first-order scheme is used.\label{fig:DMR_flag}}
\end{figure}

In the final example of non-linear system cases, we consider a double Mack reflection problem. The domain is a ramp, and initially a shock at Mach 10 is set at $x=0$. 
The mesh has $29290$  vertices, $57870$ elements and $87159$ edges, i.e., $116449$ (reps. $261478$) DoFs for the quadratic (reps. cubic) approximation. The right condition is
$$(\rho_R, u_R, v_R, p_R)=(\gamma,0,0,1),$$ so that the speed of sound is $1$ ($\gamma=1.4$) and the left state is obtained assuming a shock at Mach 10, so that
$$(\rho_L,u_L,v_L,p_L)=(8,8.25,0,116.5).$$ The solution at $T=0.2$ is shown on Figure \ref{fig:DMR_ro}. {We observe slight oscillations in the numerical solution, indicating the need for carefully designing limiting strategies or adjusting the MOOD criteria. This will be a focus of our future research.} We also display the locations where the first-order scheme is activated at the final time, see Figure \ref{fig:DMR_flag}.

\section{Conclusions}
We have presented a computational method inspired by the so-called Active Flux method, following the ideas developed in 
\cite{Abgrall_AF} and \cite{BarzukowAbgrall}, and using triangular type meshes. The method is able to capture shocks. We have illustrated the behavior of this method on several cases with smooth and discontinuous problems and models ranging from scalar to systems. The behavior of the method is as expected.

\appendix

\section{Boundary conditions for the Euler equations}\label{BCs}

We have considered two type of boundary condition on $\partial \Omega$ where $\Omega\subset \R^2$: wall boundary conditions and inflow/outflow.
We describe how they are implemented  for the update of the average values and the point values.
\subsection{Average values}
Inside the domain, $\xbar\bu_E$ will satisfy:
$$\vert E\vert \dfrac{{\rm d}\xbar\bu_E}{{\rm d}t}+\oint_{\partial E} \hbbf_\bn(\bu^h, \bu^-)\; \bla{{\rm d}\ell}=\bla{\vert E\vert \dfrac{{\rm d}\xbar\bu_E}{{\rm d}t}+}\sum_{\Gamma, \text{edge of }E}\int_{\Gamma} \hbbf_\bn(\bu^h, \bu^-)\; {\rm d}\Gamma=0,$$
where the numerical flux is the continuous one for the high-order scheme and a standard numerical flux for the first-order case. Here, $\bu^h$ is the value of $\bu$ on the boundary of $\partial E$, $\bu^-$ is the value on the opposite side of the element, and $\bn$ is the local normal.

If one edge $\Gamma$ is contained in $\partial \Omega$, we modify 
$$\int_{\Gamma} \hbbf_\bn(\bu^h, \bu^-)\; {\rm d}\ell\approx \vert \Gamma\vert \sum_{\text{ quadrature points }}\omega_q \hbbf_\bn(\bu^h(\bx_q), \bu^-(\bx_q))$$ evaluated by a quadrature rule. We use a Gaussian quadrature with 3 points. The following cases are considered:
\begin{itemize}
\item for wall, by setting the state $\bu^-(\bx_q)$ equal to the state obtained from $\bu^h(\bx_q)$ by symmetry with respect to $\Gamma$, i.e., same density, same pressure, but the velocity $\bv$ is
$$\bv^-(\bx_q)=\bv^h(\bx_q)-2\frac{\langle \bv^h(\bx_q), \bn\rangle}{\Vert \bn\Vert^2} \bn.$$
\item for inflow/outflow, since we want to impose $\bu=\bu_\infty$ weakly, we take
$\bu^-=\bu^\infty$, and the flux is some upwind flux (the Roe one for the numerical experiments).
\item for supersonic inflow/outflows or if the solution stays constant in a neighborhood of the boundary, we will impose Neumann condition, i.e., we force the solution to stay constant.
\end{itemize}

\subsection{Point values}
It is a bit more involved.
Let us remind that the update of the point values is made by \eqref{scheme:pt:1}
\begin{equation}\label{scheme:bc0}
   \dfrac{{\rm d}\bu_\sigma}{{\rm d}t}+\sum\limits_{E, \sigma\in E} \Phi_\sigma^E(\bu)=0, 
\end{equation}
that we have to modify into
\begin{equation}\label{scheme:bc}
\dfrac{{\rm d}\bu_\sigma}{{\rm d}t}+\sum\limits_{E, \sigma\in E} \Phi_\sigma^E(\bu)+\sum_{\Gamma\subset \partial\Omega, \sigma\in\Gamma}\Phi_\sigma^\Gamma(\bu)=0,
\end{equation}
and we need to define $\Phi_\sigma^\Gamma(\bu)$.
We do it in two cases: 
\begin{itemize}
\item Neumann condition: again, we force the solution to stay the same, i.e.,
$$
\sum\limits_{E, \sigma\in E} \Phi_\sigma^E(\bu)+\sum_{\Gamma\subset \partial\Omega, \sigma\in\Gamma}\Phi_\sigma^\Gamma(\bu)=0.$$
In practice, once $\sum\limits_{E, \sigma\in E} \Phi_\sigma^E(\bu)$ in 
 \eqref{scheme:bc0} has been computed, we set this quantity to $0$.
\item wall. For each $\Gamma$ in \eqref{scheme:bc},  with normal $\bn$, there is one element $E$ that has $\Gamma$ as one of its edges. Then we consider $E^-$, the (virtual) element that is symmetric with respect to $\bn$. We consider $\bu^-$, the state symmetric to $\bu^h$ with respect to $\bn$, i.e., same density and pressure, and symmetric velocity. By this we define $\Phi^\Gamma(\bu^h)$ by $\Phi^\Gamma(\bu^h)=\Phi^{E^-}(\bu^-)$.
\end{itemize}

\bibliographystyle{plain}
\bibliography{iserle}

\end{document}